\documentclass[12pt]{article}
\usepackage{amsfonts}
\usepackage{amsmath,amssymb,color}
\usepackage{graphicx}
\usepackage{subcaption}
\usepackage{array,multirow}
\usepackage{float}
\usepackage{afterpage}
\usepackage{arydshln}
\usepackage{natbib}
\usepackage{titlesec}
\usepackage{blkarray}
\usepackage{mathtools}
\usepackage[toc,page]{appendix}
\usepackage{extarrows}

\newtheorem{theorem}{Theorem}
\newtheorem{lemma}{Lemma}

\newtheorem{proposition}{Proposition}
\newtheorem{definition}{Definition}

\newtheorem{example}{Example}
\newtheorem{remark}{Remark}

\makeatletter
\def\blfootnote{\xdef\@thefnmark{}\@footnotetext}
\makeatother

\usepackage{amsfonts}
\usepackage{amsmath,amssymb,color,xcolor}
\usepackage{graphicx}
\usepackage{subcaption}
\usepackage{array,multirow}
\usepackage{float}
\usepackage{afterpage}
\usepackage{arydshln}
\usepackage{mathtools}
\usepackage{natbib}
\usepackage[onehalfspacing]{setspace}
\usepackage{makecell}
\usepackage{blkarray}
\usepackage{float}

\makeatletter
\newcommand{\doublewidetilde}[1]{{%
  \mathpalette\double@widetilde{#1}%
}}
\newcommand{\double@widetilde}[2]{%
  \sbox\z@{$\m@th#1\widetilde{#2}$}%
  \ht\z@=.9\ht\z@
  \widetilde{\box\z@}%
}
\makeatother

\newcommand*{\Cdot}{\raisebox{-0.25ex}{\scalebox{1.5}{$\cdot$}}}
\newcommand{\aaa}{\boldsymbol \alpha}

\newcommand{\ttt}{\boldsymbol \theta}
\newcommand{\TT}{\boldsymbol \Theta}

\newcommand{\rr}{\boldsymbol r}

\newcommand{\mo}{\mathbf 1}
\newcommand{\mz}{\mathbf 0}
\newcommand{\bq}{\boldsymbol q}
\newcommand{\cs}{\boldsymbol s}

\newcommand{\ba}{\boldsymbol A}

\DeclarePairedDelimiter\ceil{\lceil}{\rceil}

\newcommand{\ma}{\mbox{$\mathcal A$}}
\newcommand{\bb}{\boldsymbol b}
\newcommand{\vv}{\boldsymbol v}

\newcommand{\cc}{\boldsymbol c}

\newcommand{\ee}{\boldsymbol e}

\newcommand{\one}{\mbox{$\mathbf 1$}}

\newcommand{\pp}{{\boldsymbol p}}

\newcommand{\qq}{\boldsymbol q}
\newcommand{\cg}{\mbox{$\boldsymbol g$}}

\newcommand{\RR}{\boldsymbol R}

\newcommand{\zero}{\boldsymbol 0}
\newcommand{\bo}{\boldsymbol }
\newcommand{\eeta}{\boldsymbol \eta}
\newcommand{\mt}{\mathcal }

\newcommand{\norm}[1]{\left\lVert#1\right\rVert}

\begin{document}

%\renewcommand{\baselinestretch}{2}

%\markright{ \hbox{\footnotesize\rm Statistica Sinica
%%{\footnotesize\bf 24} (201?), 000-000
%}\hfill\\[-13pt]
%\hbox{\footnotesize\rm
%%\href{http://dx.doi.org/10.5705/ss.20??.???}{doi:http://dx.doi.org/10.5705/ss.20??.???}
%}\hfill }

%\markboth{\hfill{\footnotesize\rm YUQI GU AND GONGJUN XU} \hfill}
%{\hfill {\footnotesize\rm } \hfill}

%\renewcommand{\thefootnote}{}
%$\ $\par

%%%%%%%%%%%%%%%%%%%%%%%%%%%%%%%%%%%%%%%%%%%%%%%%%%%%%%%%%%%%%%%%%%%%%%%%%%%%%%%%%%%%%%%%%%%%%%%%%%%%%%%%%%%%%%%%%%%%%%%%%%%%

%\fontsize{12}{14pt plus.8pt minus .6pt}\selectfont \vspace{0.8pc}
%\centerline{\large\bf Sufficient and Necessary Conditions for the}
%\vspace{2pt} \centerline{\large\bf Identifiability of the $Q$-matrix}
%\vspace{.4cm} \centerline{Yuqi Gu and Gongjun Xu} \vspace{.4cm} \centerline{\it
%University of Michigan} \vspace{.55cm} %\fontsize{9}{11.5pt plus.8pt minus .6pt}
%\selectfont
\title{Sufficient and Necessary Conditions for the Identifiability of the $Q$-matrix}
%\runtitle{Partial Identifiability}
%\tableofcontents
\author{Yuqi Gu and Gongjun Xu \\ {\normalsize Department of Statistics}\\ {\normalsize University of Michigan}}
\date{}
%\begin{titlepage}
\maketitle

%\blfootnote{This work was supported in part by National Science Foundation grants SES-1659328 and DMS-1712717.}

%%%%%%%%%%%%%%%%%%%%%%%%%%%%%%%%%%%%%%%%%%%%%%%%%%%%%%%%%%%%%%%%%%%%%%%%%%%%%%%%%%%%%%%%%%%%%%%%%%%%%%%%%%%%%%%%%%%%%%%%%%%%

\begin{abstract}
{ Restricted latent class models (RLCMs)  have recently gained prominence in educational assessment, psychiatric evaluation, and medical diagnosis. Different from conventional latent class models, restrictions on RLCM model parameters are   imposed by a design matrix   to respect practitioners' scientific assumptions. The design matrix,  called the $Q$-matrix in cognitive diagnosis literature, is usually constructed by practitioners and domain experts, yet it is subjective and could be misspecified.  To address this problem, researchers have proposed to estimate the design $Q$-matrix from the data. On the other hand, the fundamental   learnability issue of the $Q$-matrix and model parameters remains underexplored and existing studies often impose stronger than needed or   even impractical conditions.  This paper proposes the sufficient and necessary   conditions for the joint identifiability of the $Q$-matrix and  RLCM model parameters. The developed identifiability conditions only depend  on the design matrix   and therefore is easy to verify in practice.
} 


\end{abstract}

%\onehalfspacing
%\title{Sufficient and Necessary Conditions for the Identifiability of the $Q$-matrix}
%\author{}
%\date{}
%\maketitle

\def\thefigure{\arabic{figure}}
\def\thetable{\arabic{table}}

\renewcommand{\theequation}{\thesection.\arabic{equation}}

%\fontsize{12}{14pt plus.8pt minus .6pt}\selectfont

\setcounter{section}{0}
\setcounter{equation}{0}

\section{Introduction}
 %\section{Introduction}

Latent class models are widely used statistical tools in social and biological sciences to model the relationship between a set of observed responses and a set of discrete latent attributes of interest. This paper focuses on a family of \textit{restricted latent class models} (RLCMs), which play a key role in various fields, including cognitive diagnosis in educational assessments \citep[e.g.,][]{Junker, HensonTemplin09, Rupp, dela2011}, psychiatric evaluation \citep{Templin,jaeger2006distinguishing,de2017analysis},  online  testing and learning  \citep{wang2016hybrid,zhang2016smart,xu2016initial}, disease etiology diagnosis  and scientifically-structured clustering of patients \citep{wu2017nested, wu2018}.

 Different from conventional latent class models, the model parameters of   RLCMs are constrained through a design matrix, often called the $Q$-matrix in the cognitive diagnosis literature \citep{Rupp}.  
 The   $Q$-matrix  encodes practitioners'  understanding of how the responses depend on the underlying latent attributes. Various RLCMs have been developed with different cognitive diagnostic assumptions \citep[e.g.,][]{DiBello,dela,Templin,davier2008general,HensonTemplin09},  including the basic  {Deterministic Input Noisy output ``And" gate (DINA)}  model \citep{Junker}, which serves as a basic submodel for more general cognitive diagnostic models.
 %, which requires only two parameters for each diagnostic item regardless of the number of latent attributes measured by the item, and the more complex generalized-DINA (GDINA) model \citep{dela2011}, which includes many commonly used diagnostic models as special cases. 
 See Section 2 for a review of these models.

Despite the popularity of   RLCMs, the fundamental identifiability issue is challenging to address.   Identifiability   of RLCMs has   long been a concern in practice, as noted in the literature   \citep*{DiBello,MarisBechger,TatsuokaC09,deCarlo2011,davier2014dina}. 
Existing identifiability results of   unrestricted latent class models in statistics \citep{teicher1967identifiability,GOODMAN74,gyllenberg1994non,allman2009identifiability} 
 %The induced constraints  of the $Q$-matrix pose  additional   challenges to the study of  the identifiability the $Q$-restricted latent class models. the    identifiability results for the unrestricted latent class models 
 cannot be directly applied to RLCMs due to the structural  constraints   induced by the  $Q$-matrix here. 
Recently,  the identifiability of  RLCM model parameters has been   studied for the basic DINA model \citep{chen2015statistical,Xu15,dina} and  general RLCMs \citep{Xu2016,partial}, assuming that the $Q$-matrix is prespecified and correct. 
%In particular, \cite{dina} firstly established the necessary and sufficient identifiability conditions of the basic DINA model parameters under a correctly specified $Q$-matrix.
  
However, the $Q$-matrix, specified by scientific  experts upon construction of the diagnostic items,   can be misspecified.  Moreover, in an exploratory analysis of newly designed items, a large part or the whole $Q$-matrix may not be available. The misspecification of the $Q$-matrix could lead to a serious lack of fit of the model and consequently inaccurate inference on the latent attribute profiles of the individuals. Therefore, it is desirable to estimate the $Q$-matrix and the model parameters jointly from the response data    \citep*[e.g.,][]{dela2,decarlo2012recognizing,JLGXZY2012,de2015general,chen2018bayesian}.
  To achieve reliable and valid estimation and inference on the $Q$-matrix, a fundamental issue is to ensure  joint identifiability of the $Q$-matrix and the associated model parameters.
%However, identifiability and related statistical properties  of the $Q$-matrix have largely been an underexplored area in the  literature and it is still not clear when the  $Q$-matrix can be consistently estimated.
Such joint identifiability  has been recently studied  in \cite{JLGXZY2011} and \cite{chen2015statistical} under the DINA model and \cite{xu2018jasa} under   general RLCMs. Nevertheless, these
 existing works mostly focus  on developing  sufficient conditions for  joint identifiability,  so they often impose stronger than needed or sometimes impractical constraints on the experimental design of cognitive diagnosis.

It  remains an open problem  what would be the minimal requirements, i.e., the necessary and sufficient conditions, for  joint identifiability of the $Q$-matrix and model  parameters. 
This paper addresses this problem and has the following contributions.
 
 First, under the   DINA model, we derive the  necessary and sufficient conditions for joint identifiability of the $Q$-matrix and the associated DINA model parameters. 
 {Our necessary and sufficient conditions are succinctly and neatly written as three algebraic properties of the $Q$-matrix, which we summarize as \textit{completeness} (Condition $A$), \textit{distinctness} (Condition $B$), and \textit{repetition} (Condition $C$); please see Theorem \ref{thm-dina} for details.
 In plain words, these three conditions require the binary $Q$-matrix to be \textit{complete} by containing an identity submatrix, to have all columns \textit{distinct} other than the part of the identity submatrix, and %to contain at least three items \textit{repeatedly} measuring each latent attribute.
 to \textit{repeatedly} contain at least three entires of ``$1$" in each column.} 
 % end of blue color
%%%
The proposed conditions not only guarantee identifiability, but also give the minimal requirements for  the $Q$-matrix and DINA model parameters to be estimable from the observed responses. 
 The identifiability result can be directly applied to the {Deterministic Input Noisy output ``Or" gate (DINO)} model \citep{Templin}, due to  the duality of the DINA and DINO models \citep{chen2015statistical}. The derived identifiability conditions also serve as necessary requirements for joint  identifiability under general RLCMs that cover the DINA as a submodel.
 
 Second,  we study a weaker notation of identifiability, the so-called generic identifiability, and propose  sufficient  and necessary conditions for it under both the DINA  model and general RLCMs.
Generic identifiability  implies that those parameters for which identifiability does not hold live in a set of Lebesgue measure zero \citep{allman2009identifiability}. The motivation for studying generic identifiability is that strict identifiability conditions sometimes could be too restrictive in practice. 
 For instance,   it is known that   unrestricted latent class models  are not strictly identifiable \citep{gyllenberg1994non}, while they are generically identifiable under certain conditions \citep{allman2009identifiability}.
 However, as to  RLCMs,  the model parameters are forced by the $Q$-matrix-induced constraints to fall in a measure zero subset of the parameter space, and thus existing results for unrestricted models cannot be directly applied. It is   unknown what generic identifiability  conditions are needed to jointly identify the $Q$-matrix as well as the model parameters.
In this work, we propose sufficient and necessary conditions for generic identifiability, and explicitly characterize the non-identifiable measure-zero subset.
Our mild sufficient conditions for generic identifiability under general RLCMs can be summarized as the following properties of the $Q$-matrix, \textit{double generic completeness} (Condition $D$) and \textit{generic repetition} (Condition $E$); see Theorem \ref{thm-multi} for details. In plain words, these two conditions require the binary $Q$-matrix to contain two \textit{generically complete} square submatrices with all diagonal elements equal to ``$1$", and to additionally (\textit{repeatedly}) contain at least one entry of ``$1$" other than the part of these two submatrices.

%\begin{figure}[h!] 
%\centering
%\includegraphics[width=0.8\textwidth]{SS_flowchart}
%\caption{\blue{Flowchart of the main theoretical results. 
%The ``ID." is short for ``identifiability". See Theorem \ref{thm-dina} for Conditions $A$, $B$ and $C$, and   Theorem \ref{thm-multi} for   $D$ and $E$.}}
%\label{fig-flowmain}
%\end{figure}

The rest of the paper is organized as follows.
Section \ref{sec-review} gives an introduction to   RLCMs   and reviews some popular models in cognitive diagnosis. 
Section \ref{sec-def} introduces the definitions of strict and generic identifiability for RLCMs, and presents an illustrative example.
Sections \ref{sec-dina} and \ref{sec-general} contain  our main theoretical results for strict and generic identifiability, for DINA model and general RLCMs, respectively. Section \ref{sec-sum} gives some discussions. All the proofs of the theoretical results and additional simulation studies that verify the developed theory are included in the Supplementary Material.      The Matlab codes for checking all the proposed conditions are  available at \verb|https://github.com/yuqigu/Identify_Q|.

\section[Cognitive Diagnosis]{RLCMs for Cognitive Diagnosis}\label{sec-review}
RLCMs  are key statistical tools in cognitive diagnostic assessments with the aim to estimate individuals'  attribute profiles based on their response data in the assessment.
Specifically,   consider   a diagnostic test with $J$ items. A subject (such as an examinee or a patient) provides a $J$-dimensional binary response vector $\RR = (R_1,...,R_J)^\top$ to the $J$ items.  These responses are assumed to be dependent in a certain way on $K$ unobserved latent attributes.
Under  RLCMs, a complete set of $K$ latent attributes is known as a latent class or an attribute profile,
 denoted by a vector $\aaa= (\alpha_1,\ldots, \alpha_K)^\top$, where $\alpha_k \in \{0,1\}$ is  a binary   indicator of the absence or presence of the $k$th attribute, respectively.

RLCMs assume a two-step data generating process.
 The first step has a population model for the attribute profile vector.
We assume   that  the attribute profile        follows  a  categorical distribution with population proportions
$\pp:=(
p_{\aaa}: \aaa\in\{0,1\}^{K})^\top$
where $p_{\aaa} > 0$ for all $\aaa\in\{0,1\}^K$ and $\sum_{\aaa\in\{0,1\}^{K}}p_{\aaa}=1$.
% The distribution of $\aaa$ is thus characterized by the parameter vector $\pp$.

 The second step of the data generating process follows a   latent class model framework, incorporating  constraints based on the underlying cognitive processes.
Given a subject's attribute profile $\aaa$, his/her responses to the $J$ items $\{R_{j}: j=1,\cdots, J\}$, are
assumed conditionally independent, and each  $R_j$  follows a Bernoulli distribution with  parameter $\theta_{j,\aaa}= P(R_j=1 \mid \aaa) $. The $\theta_{j,\aaa}$ denotes the probability of  a positive response, and is also called an item parameter of item $j$.
 The collection of all the item parameters, denoted by the item parameter matrix $\TT = (\theta_{j,\aaa})_{J\times 2^K}$, are further constrained by the design matrix $Q$. The $Q$-matrix is the key structure that specifies the relationship between the $J$ items and the $K$ latent attributes.
Specifically, the $Q$-matrix is a $J\times K$ binary matrix, with entries $q_{j,k}\in\{1,0\}$ indicating whether or not   the $j$th item is linked to the $k$th latent attribute.
When $q_{j,k}=1$, we say attribute $k$ is required by item $j$.
The $j$th row vector  $\qq_j$ of $Q$ gives the full attribute requirements of item $j$.
Given an attribute profile $\aaa$ and a matrix $Q$, we 
write $\aaa\succeq \qq_j  \mbox{ if }  \alpha_k \geq q_{j,k} \mbox{ for all }  k \in\{1, \ldots, K\},$
and 
$\aaa\nsucceq \qq_j  \mbox{ if there exists $k$ such that } \alpha_k < q_{jk};$
similarly we define the operations $\preceq$ and $\npreceq$.

  If $\aaa\succeq \qq_j$, a subject having attribute pattern $\aaa$ possesses all the attributes required by item $j$ specified by the $Q$-matrix, and would be   ``capable'' of answering item $j$ correctly. On the other hand, if $\aaa'\nsucceq \qq_j$, the subject with $\aaa'$ misses some required attribute of item $j$ and is expected to have a smaller  positive response probability than those subjects with $\aaa\succeq \qq_j$.
That is, the RLCMs we consider in this paper assume
\begin{align}\label{eq-mono}
\theta_{j,\aaa}&
> \theta_{j,\aaa'} \text{ for any } \aaa\succeq \qq_j \text{ and } \aaa'\nsucceq \qq_j.
%\\ 
%\notag
%\theta_{j,\aaa}& = \theta_{j,\aaa'} \text{ for any }\aaa\otimes \qq_j=\aaa'\otimes \qq_j,
\end{align} 
{Such \textit{monotonicity assumption} in \eqref{eq-mono} is common to most RLCMs.}
% In some practical scenarios, a more stringent monotonicity assumption $\theta_{j,\aaa}>\theta_{j,\aaa'}$ for any $\aaa\odot \qq_j\succ\aaa'\odot \qq_j$ would be employed \citep{HensonTemplin09}, where ``$\odot$'' denotes the elementwise multiplication operator.
  Another common assumption of RLCMs is that mastering those non-required attributes of an item will not change the positive response probability to it, i.e., $\theta_{j,\aaa}=\theta_{j,\aaa'}$ if $\aaa\odot \qq_j=\aaa'\odot \qq_j$, where ``$\odot$'' denotes the elementwise multiplication operator \citep{HensonTemplin09}. Under the introduced setup, the response vector $\RR$ has probability mass function in the form 
\begin{equation}
	\label{eq-pmf}
\mathbb P(\RR=\rr\mid  Q, \TT,\pp)=\sum_{\aaa\in\{0,1\}^K}p_{\aaa}\prod_{j=1}^J \theta_{j,\aaa}^{r_j} (1-\theta_{j,\aaa})^{1-r_j},\quad  \rr\in\{0,1\}^J,
\end{equation}
where the constraints on the $\theta_{j,\aaa}$'s imposed by $Q$ are made implicit. 
%Note that \eqref{assumption2} implies that $\max_{\aaa\succeq\qq_j}\theta_{j,\aaa}=\min_{\aaa\succeq\qq_j}\theta_{j,\aaa},$ which is a key assumption for the   identifiability of the restricted latent class model parameters \citep{Xu2016}.
 
Next, we review some popular  cognitive diagnosis models   and illustrate how they fall into the family of RLCMs. 
\begin{example}[DINA model]\label{exp-dina}
\normalfont{
One of the basic cognitive diagnosis models is the DINA model \citep{Junker}. The DINA model assumes a conjunctive relationship among attributes, meaning that to be capable of providing a positive response to an item, it is necessary to possess all its required attributes indicated by the $Q$-matrix.
  For an item $j$  and  a subject with attribute profile $\aaa$, an   ideal response under the DINA model is defined as $\Gamma_{j,\aaa}^{DINA} = I(\aaa\succeq\bq_j)$, which indicates whether the subject is capable of item $j$.  
The uncertainty is incorporated at the item level with the slipping parameter $s_j  = P(R_j=0\mid \Gamma_{j,\aaa}  =1)$ denoting the probability that a capable subject slips the positive response, and the guessing parameter, $g_j  = P(R_j=1\mid \Gamma_{j,\aaa}  =0)$ denoting the probability that a non-capable subject coincidentally gives the positive response by guessing.
Then the positive response probability for item $j$ of class $\aaa$ is  
 $
\theta_{j,\aaa}^{DINA} = (1-s_j)^{\Gamma_{j,\aaa} } g_j^{1-\Gamma_{j,\aaa} }.
$ 
The DINA model has only two parameters $s_j$ and $g_j$ for each item regardless of the number of attributes required by the item. In the following discussion, we denote $\cs =(s_1,\ldots,s_J)^\top$ and $\cg =(s_1,\ldots,s_J)^\top$. Given the $Q$-matrix, the DINA model parameters $(\TT,\pp)$ can then be equivalently expressed by $  (\cs,\cg,\pp) $.
We further assume  $\mathbf{1}-\cs\succ \cg$ \citep{Xu15}, which makes DINA satisfy the monotonicity assumption \eqref{eq-mono}.  
%For such diagnostic models with each item   associated with exactly two item parameters, we call them  Two-Parameter  RLCMs.  
Identifiability results of the basic DINA model are presented in Section \ref{sec-dina}.
}
\end{example}

\begin{example}[GDINA   model and General RLCMs]\label{exp-gdina}
\normalfont{ \cite{dela2011} extended the DINA model to the
Generalized DINA (GDINA) model.
The  formulation of the GDINA model based on $\theta_{j,\aaa}$ can be decomposed into the sum of the effects {due to}
the presence of specific attributes and their interactions. Specifically, for an item $j$ with $\qq$-vector $\qq_j=(q_{j,k}: k=1, \cdots, K)$, the positive response probability is
\begin{eqnarray}\label{GGDINA}
\theta_{j,\aaa}^{GDINA}  &=&
\sum_{\mathcal S\subseteq \{1,\ldots,K\}}\beta_{j,\mathcal S}\prod_{k\in\mathcal S} q_{j,k}\prod_{k\in \mathcal S}\alpha_k.
%\quad \mathcal K_j = \{1\leq k\leq K:\,q_{j,k}=1\}.
%\beta_{j,0} 
%+ \sum_{k=1}^{K}\beta_{j,k}(q_{j,k}\alpha_{k})
%+ \sum_{k'=k+1}^{K}\sum_{k=1}^{K-1}\beta_{j,kk'}(q_{j,k}\alpha_{k})(q_{j,k'}\alpha_{k'})\\
%\notag
%&&+\cdots + \beta_{j,12\cdots K}\prod_{k}(q_{j,k}\alpha_{k}).
\end{eqnarray}
%%Each $\beta_{j,S}$ models the change of the positive response probability resulting from the mastery of the attributes in the set $S$, where $S$ ranges from every subset of  $\mathcal K_j$.
Note that  not all $\beta$-coefficients in the above equation are included in the model. For a subset $\mathcal S$ of the $K$ attributes $\{1,\ldots,K\}$, the $\beta_{j,\mathcal S}\neq 0$ only if $\prod_{k\in \mathcal S}q_{j,k}=1$.
%For instance, when $\qq_j \neq \mathbf 1^\top$, one does not need parameter $ \beta_{j,12\cdots K}$ since $\prod_{k}(q_{j,k}\alpha_{k})=0$.
The interpretation is that, $\beta_{j,\varnothing}$ denotes the probability of a positive response when none of the required attributes are present in $\aaa$;
 when $q_{j,k}=1$, $\beta_{j,\{k\}}$ is in the model, 
representing the change in the positive response probability resulting from the mastery of a single attribute $k$; 
when $q_{j,k}=q_{j,k'}=1$, $\beta_{j,\{k,k'\}}$ is in the model,
 representing the change in the positive response probability due to the interaction effect of mastery of both 
$k$ and ${k'}$.
%similarly, when $\qq_j=\mathbf 1^\top$, $\beta_{j,\{1,2,\cdots,K\}}$ is the change in the positive response probability due to the interaction effect of  mastery of all the required attributes.
 Under the GDINA model, each $\theta_{j,\aaa}$ models  the main effects and all the interaction effects  of the attributes measured by the item. For such diagnostic models, we call them  {\it general}  RLCMs. Another popular general RLCM is the Log-linear Cognitive Diagnosis Model \citep[LCDM;][]{HensonTemplin09} and the General Diagnostic Model \citep[GDM;][]{davier2008general}. Identifiability  results of   general  RLCMs are presented in Section \ref{sec-general}. 
}
\end{example}

%\section{Identifiability Results}\label{sec-id}
%
%\noindent{\bf 3.1~  Definitions of  identifiability and generic identifiability}
% \bigskip
%
% \noindent 
 
 \section[Definitions and Illustrations]{Definitions and Illustrations of strict and generic identifiability}\label{sec-def}
%This section studies the identifiability of the $Q$-matrix and the associated RLCM parameters $(\TT,\pp)$.
This section introduces the definitions of joint strict identifiability and joint generic identifiability of $(Q,\TT,\pp)$ for RLCMs, and gives an illustrative example. 
%That is, identifiability of $(Q,\TT,\pp)$.
%Note that $\TT$ depends on the $Q$-matrix as certain elements in the matrix $\TT$ are restricted to be equal under both the $Q$-matrix structure and specific diagnostic model assumptions; with a slight abuse of notation, here we just use $\TT$ to denote the constrained parameter matrix under the corresponding $Q$-matrix.

We would also like to point out that {the monotonicity assumption   stated in \eqref{eq-mono}}, is necessary  for the identifiability of the $Q$-matrix. Since otherwise any $Q\neq \mathbf 1_{J\times K}$ with parameters $(\TT,\pp)$ can not be distinguished from $\bar Q= \mathbf 1_{J\times K}$ with the same parameters $(\TT,\pp)$ under the general RLCM. The monotonicity constraints ensure  that the constraints induced by  $Q\neq \mathbf 1_{J\times K}$ and  $\bar Q= \mathbf 1_{J\times K}$ cannot be  the same and therefore $Q$ can be identified under additional conditions to be discussed in Sections \ref{sec-dina} and \ref{sec-general}. In the following we assume the monotonicity assumption introduced in Section 2 is satisfied.

Another common issue with identifiability of the $Q$-matrix is label swapping. In the setting of   RLCMs, arbitrarily reordering columns of a $Q$-matrix would not change the distribution of the responses. As a consequence, it is only possible to identify $Q$ up to column permutation, and we will write $\bar Q\sim  Q$ if $\bar Q$ and $Q$ have an identical set of column vectors, and write $(\bar Q, \bar\TT,\bar\pp)\sim (Q, \TT,\pp)$ if $\bar Q\sim  Q$ and $(\bar\TT,\bar\pp)=(\TT,\pp)$.

We first introduce the definition  of   identifiability of $Q$-matrix as well as the model parameters $(\TT,\pp)$ , which we term as \textit{joint strict identifiability}.

\begin{definition}[Joint Strict Identifiability]
Under an RLCM, the  design matrix $Q$ joint  with the model parameters $(\TT,\pp)$ are said to be strictly  identifiable  if for any $(Q, \TT,\pp)$, there is no $(\bar Q, \bar\TT,\bar\pp)\nsim (Q, \TT,\pp)$ such that
 \begin{equation}\label{eq-orig}
 \mathbb P(\RR=\rr\mid  Q, \TT,\pp) =\mathbb  P(\RR=\rr\mid \bar Q, \bar \TT,\bar\pp) ~\mbox{ for all } ~ \rr\in \{0,1\}^J. 
 \end{equation}
 In the following discussion, we will write \eqref{eq-orig} simply as $\mathbb P(\RR\mid  Q, \TT,\pp) =\mathbb  P(\RR\mid \bar Q, \bar \TT,\bar\pp)$.
 \end{definition}

 Despite being the most stringent criterion for identifiability, strict identifiability could be too restrictive, ruling out many cases where the $(Q,\TT,\pp)$ are ``almost surely" identifiable. 
In the literature of unrestricted latent class models, \cite{allman} proposed and studied the so-called \textit{generic identifiability}. Here we also introduce the concept of generic identifiability for RLCMs as follows.

\begin{definition}[Joint Generic Identifiability] \label{def-gen}
Consider an RLCM with  parameter space $\boldsymbol\vartheta_Q$, which is of full dimension in $\mathbb R^m$ {with $m$ corresponding to the number of free parameters in the model}. The matrix $Q$ joint with the model parameters $(\TT,\pp)$ are said to be generically identifiable, 
%if $(Q,\TT,\pp)$ are strictly identifiable excluding a zero measure set of $\mathbb R^m$.
if the following set has Lebesgue measure zero in $\mathbb R^m$:
$
\boldsymbol\vartheta_{non} =\{
  (\TT,\pp):\,\exists  (\bar Q, \bar\TT,\bar\pp) \nsim (Q, \TT,\pp)~\text{such that}~ 
  \mathbb P(\RR\mid Q,\TT,\pp) = \mathbb P(\RR\mid\bar Q,\bar\TT,\bar\pp)\}.
$
\end{definition}
%In Definition \ref{def-gen}, $m=2^K-1+2J$ for the DINA model

\subsection{Illustration of Generic Identifiability Phenomenon with $Q_{4\times 2}$}
\begin{figure}%[H]
\centering
\begin{subfigure}{0.5\textwidth}
\includegraphics[width=\linewidth]{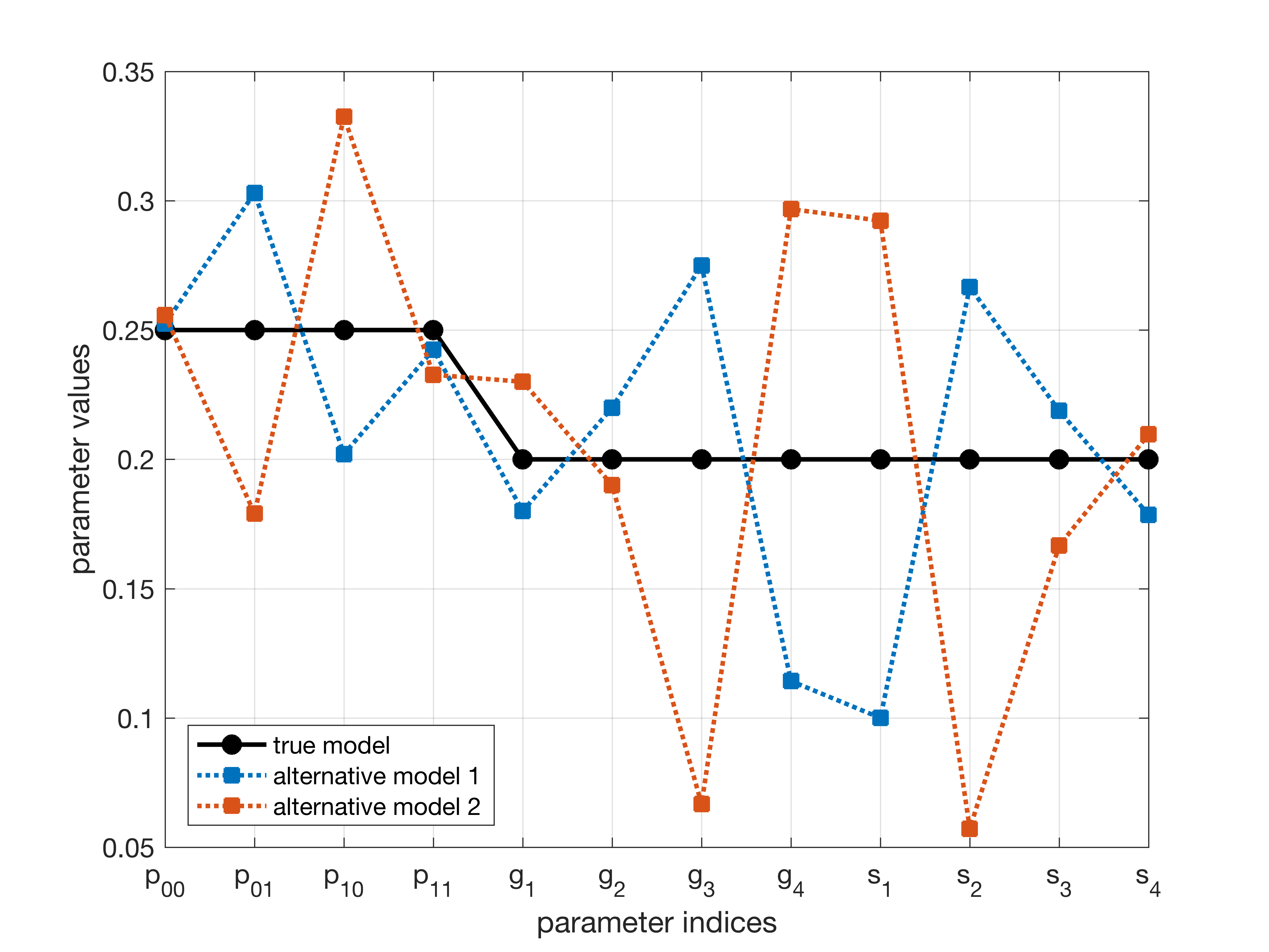}
\caption{} \label{fig:a24}
\end{subfigure}\hspace*{\fill}
\begin{subfigure}{0.5\textwidth}
\includegraphics[width=\linewidth]{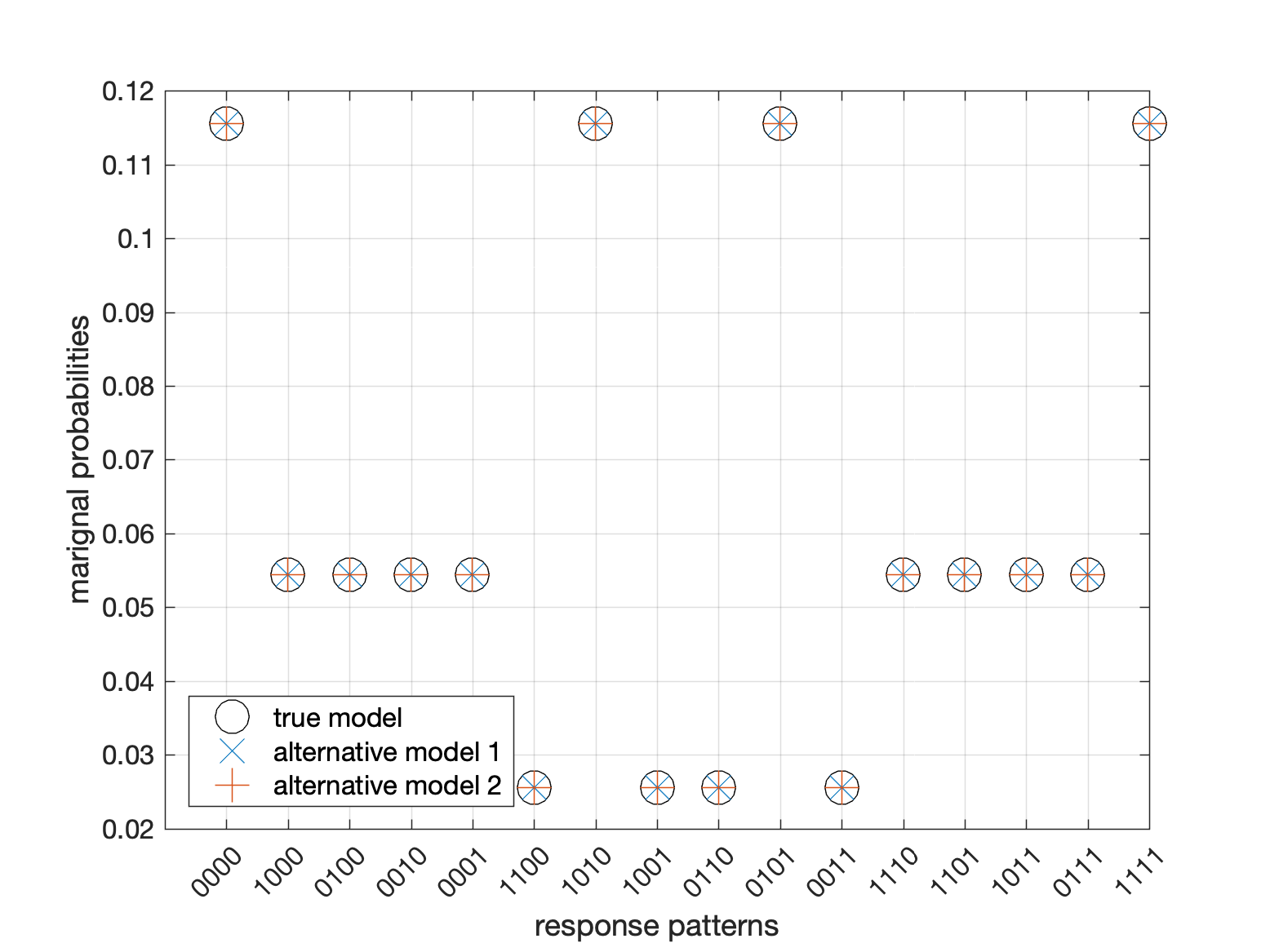}
\caption{} \label{fig:b24}
\end{subfigure}
\caption{Illustration of non-identifiability under   $Q_{4\times 2}$ in scenario (a).}
\label{fig-Q24-a}
\end{figure}

We use an example to show the difference between generic identifiability and strict identifiability. Consider the $Q$-matrix $Q_{4\times 2}$ in \eqref{eq-q24}. Under the DINA model, it will be proved  that this   $Q$-matrix joint with the associated   model parameters $(\cs,\cg,\pp)$ are   generically identifiable (by part (b.2) of Theorem \ref{thm-dina-gen}), but not strictly identifiable (by Theorem \ref{thm-dina}). 
\begin{equation}\label{eq-q24}
Q_{4\times 2} = \begin{pmatrix}
1 & 0 & 1& 0\\
0 & 1 & 0 & 1\\
%1 & 0 \\
%0 & 1 \\
\end{pmatrix}^\top.
\end{equation}
In particular, as long as the true proportions $\pp=(p_{(00)},p_{(01)},p_{(10)},p_{(11)})$ satisfy the following inequality constraint, $(Q_{4\times 2},\cs,\cg,\pp)$ are identifiable (see proof of Theorem \ref{thm-dina-gen} (b.2) for reason):
\begin{equation}\label{eq-q24-id}
p_{(01)}p_{(10)}\neq p_{(00)}p_{(11)}.
\end{equation}
%%%
%We call constraints like \eqref{eq-q24-id}  \textit{generic identifiability constraints}, since for RLCM parameters, these constraints  hold almost surely with respect to the Lebuesgue measure in the parameter space. 
On the other hand, when $p_{(01)}p_{(10)}= p_{(00)}p_{(11)}$, the model parameters are not identifiable and there exist infinitely many sets of parameters providing the same distribution of the observed response vector.
Here the parameter space $\boldsymbol\vartheta_Q=\{(\cs,\cg,\pp):\, \mathbf1 - \cs\succ \cg,~ \pp\succ\mathbf0,~ \sum_{\aaa} p_{\aaa}=1 \}$ is of full dimension in $\mathbb R^{11}$, while the non-identifiable subset $\boldsymbol\vartheta_{non} = \{(\cs,\cg,\pp):\, p_{(01)}p_{(10)}= p_{(00)}p_{(11)}\}$ has Lebesgue measure zero in $\mathbb R^{11}$. 
We use a simulation study to illustrate the generic identifiability phenomenon. 
Under the $Q_{4\times 2}$ in \eqref{eq-q24}, consider the following two simulation scenarios, 
\begin{itemize}
\item[(a)] the true model parameters are set to be $g_j=s_j=0.2$ for $j=1,2,3,4$ and $p_{(00)}=p_{(01)}=p_{(10)}=p_{(11)}=0.25$, which violates \eqref{eq-q24-id}; 
\item[(b)] the true model parameters are randomly generated, which almost always satisfy \eqref{eq-q24-id}.
 Specifically, we randomly generate a total number of 100 true parameter sets $(\cs,\cg,\pp)$, with the following generating mechanism,
$ s_j\sim \mathcal U(0.1,0.3)$,  $g_j\sim \mathcal U(0.1,0.3)$ for  $j=1,2,3,4$
and $\pp  \sim \text{Dirichlet}(3,3,3,3)$.
 Here $\mathcal U(0.1,0.3)$ denotes the uniform distribution on   $[0.1,0.3]$, and Dirichlet$(3,3,3,3)$ denotes the Dirichlet distribution with parameter vector $(3,3,3,3)$. 
\end{itemize}
We show numerically that in scenario (a), there exist multiple different sets of valid DINA parameters that give the same distribution of $\RR$; while in scenario (b), the model $(Q,\cs,\cg,\pp)$ are almost surely identifiable and estimable.
 In particular, corresponding to scenario (a), Figure \ref{fig-Q24-a} (a) plots the true model parameters as well as the other two sets of valid DINA model parameters (constructed based on the derivations in the proof of Theorem \ref{thm-dina-gen} (b.2)), and Figure \ref{fig-Q24-a} (b) plots the marginal probabilities of all the $2^4=16$ response patterns under these three different sets of model parameters. We can see that despite these three sets of parameters are very different, they give the identical distribution of the $4$-dimensional binary response vector. 
 
%We point out that the alternative models 1 and 2 are precisely constructed based on our proof of generic identifiability of $(Q_{4\times 2},\cs,\cg,\pp)$ in Theorem \ref{thm-dina-gen}. In particular, the parameter specification in scenario (a) falls within a measure zero set $\boldsymbol{\vartheta}_{non}$ of the entire parameter space where non-identifiability can occur, and we utilize this property to construct the alternative model parameters. For details of construction, see ...

\begin{figure}[h!]
\centering
\begin{subfigure}{0.32\textwidth}
\includegraphics[width=\linewidth]{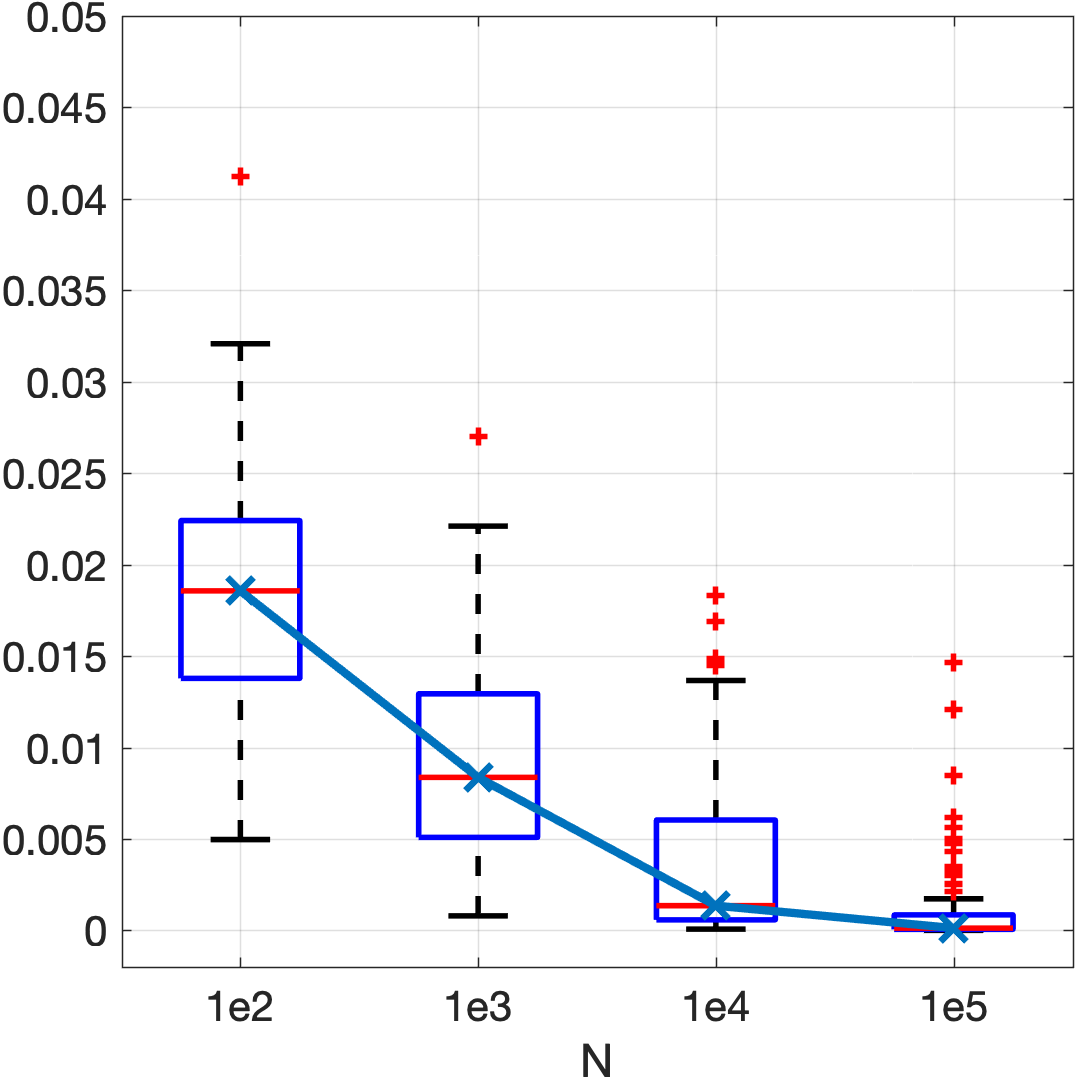}
\caption{MSE of $\pp$} \label{fig:a24}
\end{subfigure}
\hspace{\fill}
\begin{subfigure}{0.32\textwidth}
\includegraphics[width=\linewidth]{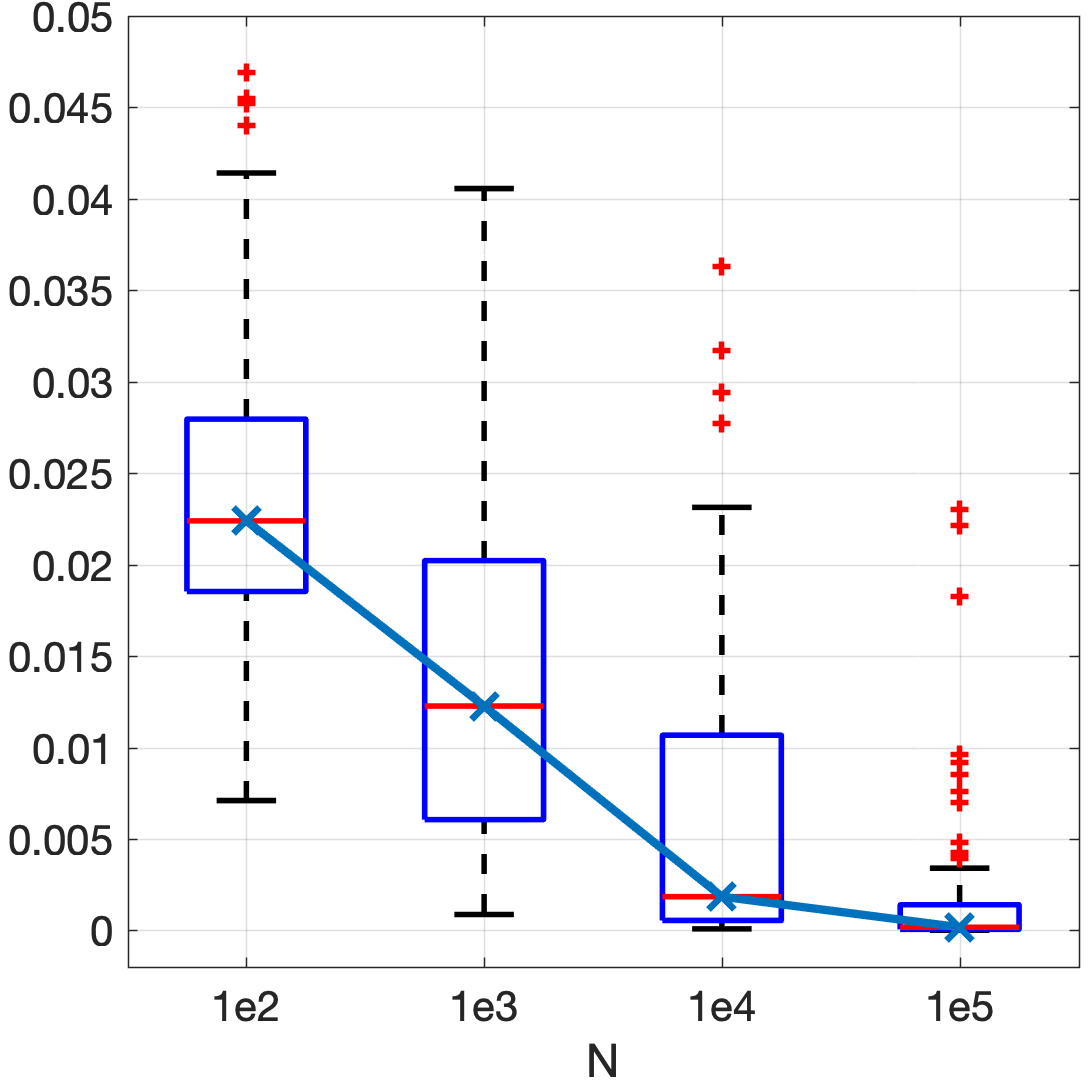}
\caption{MSE of $\cs$} \label{fig:b24}
\end{subfigure}
\hspace{\fill}
\begin{subfigure}{0.32\textwidth}
\includegraphics[width=\linewidth]{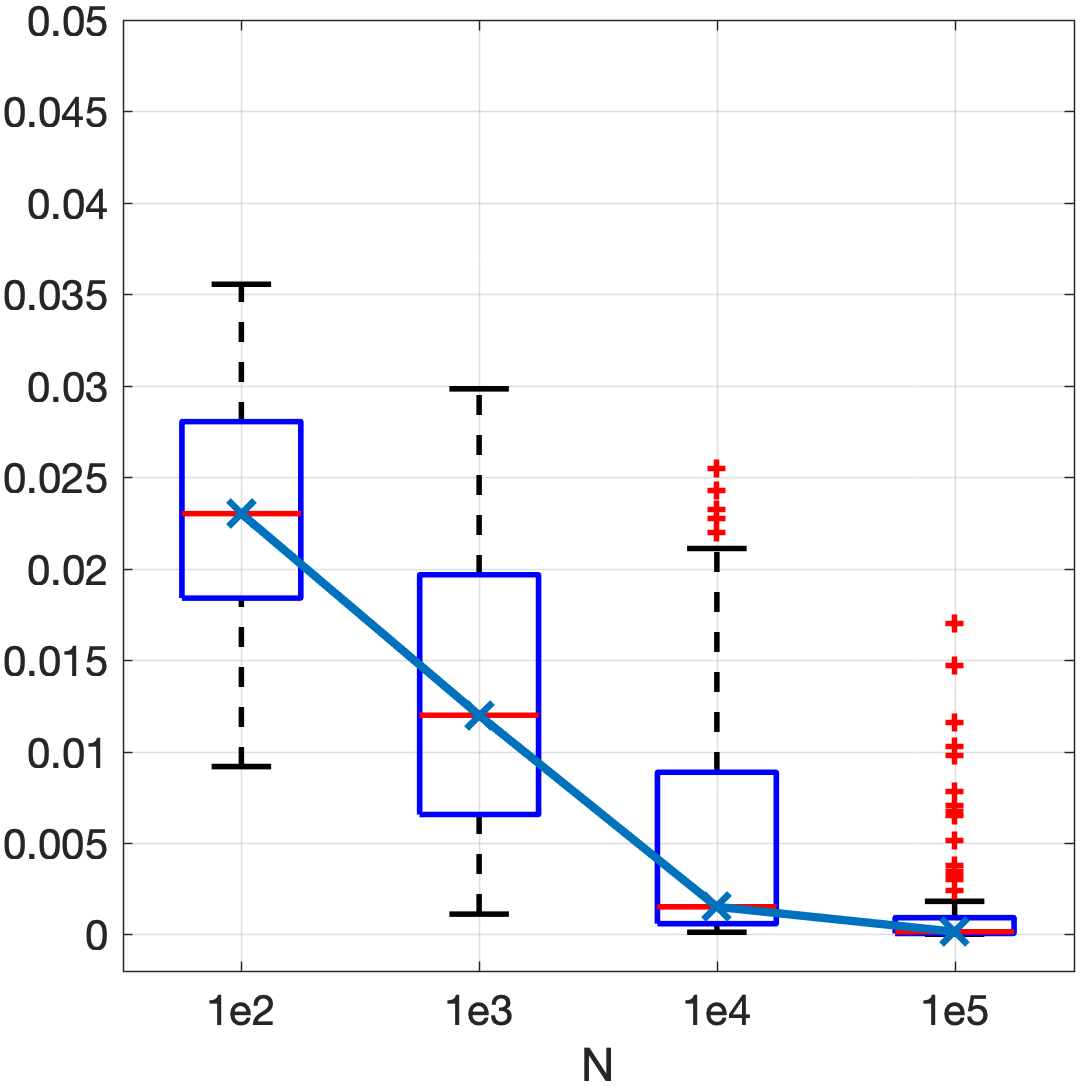}
\caption{MSE of $\cg$} \label{fig:c24}
\end{subfigure}
\caption{{Illustration of generic identifiability under $Q_{4\times 2}$, which corresponds to simulation scenario (b).}}
\label{fig-Q24-b}
\end{figure}

Corresponding to scenario (b), 
%Figure \ref{fig-Q24-b} shows the Mean Square Errors (MSEs) of the randomly generated true parameters generally decrease as the sample size $N$ increases. 
  we randomly generate $B = 100$ sets of true parameters $(\cs^i,\cg^i,\pp^i)$ for $i=1,\ldots,100$.
Then for each $(\cs^i,\cg^i,\pp^i)$,  we generate 200 independent datasets of   size $N$ with   $N=10^2$, $10^3$, $10^4$ and $10^5$, and then compute the Mean Square Error (MSE) of  the maximum likelihood estimators (MLE) of the slipping, guessing and proportion parameters, respectively. {To compute the MLE of model parameters for each simulated dataset, we run the EM algorithm with 10 random initializations and choose the estimators achieving the largest log-likelihood value out of the 10 runs. }%end of blue
Figure \ref{fig-Q24-b} shows the boxplots of Mean Square Errors (MSEs) associated with the $B = 100$ true parameter sets for each sample size $N$. As $N$ increases, we observe that the MSEs decrease  to zero, indicating the (generic) identifiability of these randomly generated parameters. 

%%%%%%%
 % normal font ends here
 
\begin{figure}[h!]
\centering
\includegraphics[width=0.65\linewidth]{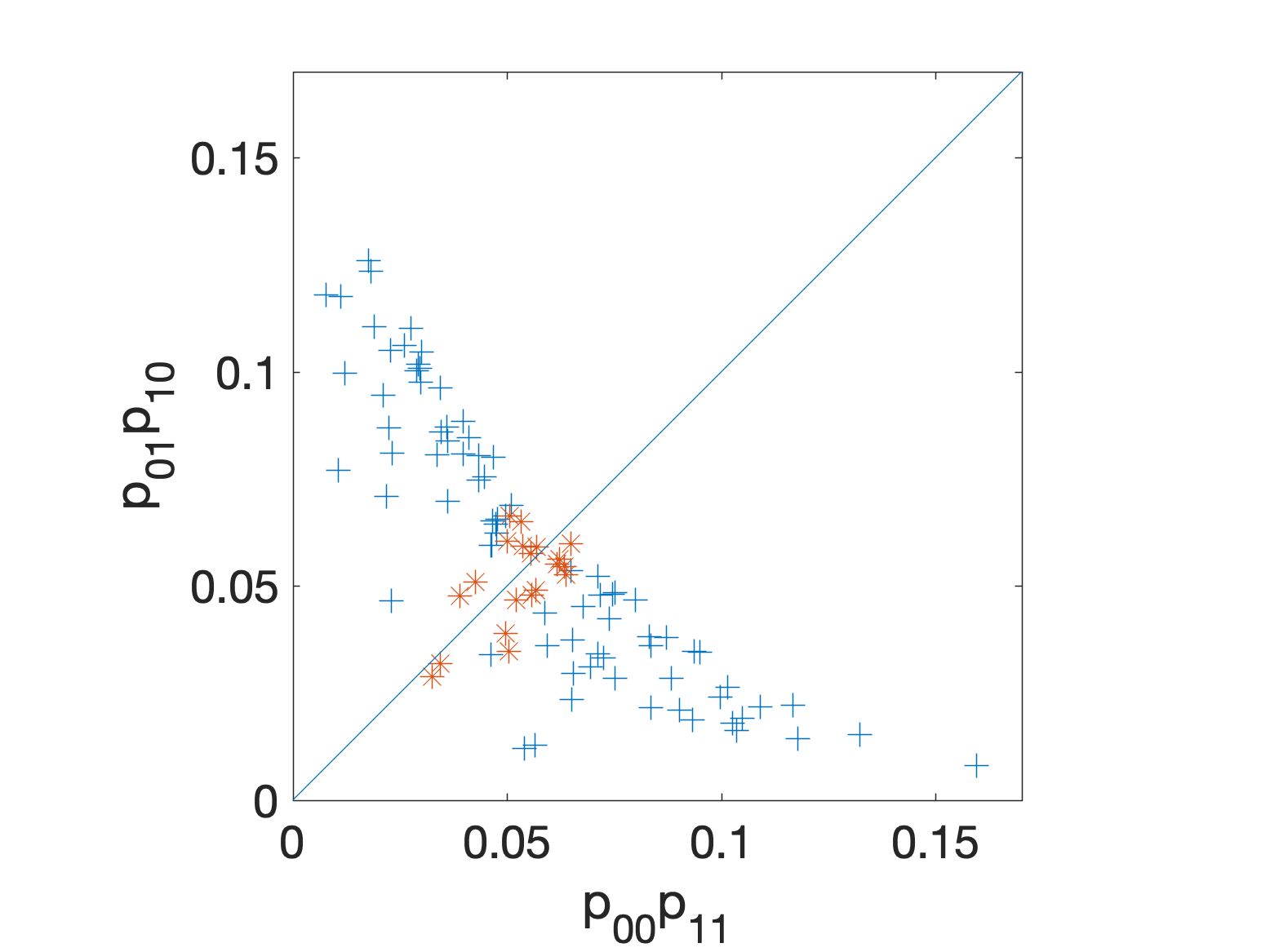}
\caption{Illustration of impact of the generic identifiability constraint \eqref{eq-q24-id}. Red ``$\ast$"s represent parameter sets with the $20\%$ largest MSEs in Figure \ref{fig-Q24-b}(a) with $N=10^5$; blue ``$+$"s represent the remaining parameter sets.}
\label{fig-Q24-ratio}
\end{figure}

On the other hand, Figure \ref{fig-Q24-b} also shows that there do exist several parameter sets whose MSEs are ``outliers" in the boxplots and converge to 0 much slower than others  as $N$ increases. This happens basically because these sets of parameters fall near the non-identifiability  set, $\mathcal V_{non}=\{(\cs,\cg,\pp):~p_{(01)}p_{(10)}- p_{(00)}p_{(11)}=0\}$, and it becomes more difficult to identify them  than others.
To illustrate this point, we consider the scenario corresponding to the rightmost boxplot in Figure \ref{fig-Q24-b}(a) with sample size $N=10^5$. 
For each one of the 100 sets of true parameters $(\cs^i,\cg^i,\pp^i)$, in Figure \ref{fig-Q24-ratio} we plot $p^i_{(00)}\cdot p^i_{(11)}$ and $p^i_{(01)}\cdot p^i_{(01)}$ as the $x$-axis and $y$-axis coordinates, respectively. Then each point represents one set of true parameters used to generate the data. Specifically, we plot those parameter sets with red ``$\ast$"s if their corresponding MSEs are the $20\%$ largest outliers in the rightmost boxplot in Figure \ref{fig-Q24-b}(a); and plot the remaining $80\%$ parameter sets with blue ``$+$"s.
One can clearly see that the closer the true parameters lie to the non-identifiability set $\mathcal V_{non}=\{(\cs,\cg,\pp):~p_{(01)}p_{(10)}-p_{(00)}p_{(11)}=0\}$ (represented by the straight reference line drawn from $(0,0)$ to $(0.17,0.17)$), the larger the MSEs are, and the slower the convergence rate of the MLEs is.
This indicates the phenomenon under generic identifiability that when the true model is close to the non-identifiable set, the convergence of their MLEs becomes slow.
%and larger than usual sample sizes are required to achieve reasonably small MSEs

Interestingly, the generic identifiability constraint \eqref{eq-q24-id} is equivalent to the statement that the two latent attributes are \textit{not independent} of each other. To see this, view each subject's 2-dimensional attribute profile as a random vector taking values in a $2\times 2$ contingency table. Then \eqref{eq-q24-id} states that the $2\times 2$ matrix of joint probabilities of attributes mastery
$$
\begin{pmatrix}
p_{(00)} & p_{(01)} \\
p_{(10)} & p_{(11)} \\
\end{pmatrix}
$$
has full rank with nonzero determinant $p_{(00)}p_{(11)}-p_{(01)}p_{(10)}$. This means one row (resp. column) of the matrix can not be a multiple of the other row (resp. column), and hence the two binary attributes can not be independent.
%the two attributes are not independent.
Intuitively, this implies that the DINA model essentially requires each attribute to be measured by at least three times for identifiability (as shown in Condition $B$ in Theorem \ref{thm-dina}). In particular, consider those attributes that are measured by only two items in the $Q$-matrix. If these attributes are independent, then intuitively they provide independent source of information, in which case the model is not identifiable. However, if these attributes are  dependent, then the dependency instead helps with the identification of the model structure.

\bigskip
Before stating  the strict and generic  identifiability results on $(Q,\TT,\pp)$, we   show in the next proposition that any all-zero row vector in the $Q$-matrix can be dropped without impacting the identifiability conclusion.

\begin{proposition}\label{prop-0vector}
Suppose the $Q$-matrix of size $J\times K$ takes the form $Q=((Q')^\top,\mathbf0^\top)^\top$,
%\[Q=\begin{pmatrix}
%Q'\\
%\mathbf0\\
%\end{pmatrix},\]
where $Q'$ is a $J'\times K$ submatrix containing the $J'$ nonzero $\bq$-vectors, and $\mathbf0$ denotes a $(J-J')\times K$ submatrix containing those zero $\qq$-vectors. Let $\TT'$ be the submatrix of $\TT$ containing its first $J'$ rows. Then for any RLCM, $(Q,\TT,\pp)$ are jointly strictly (generically) identifiable if and only if   $(Q',\TT',\pp)$ are jointly strictly (generically) identifiable.
\end{proposition}
 
 Therefore, without loss of generality, from now on we only consider $Q$-matrices without any zero $\qq$-vectors when studying joint identifiability. %of $(Q,\TT,\pp)$.
 We study various RLCMs that are popular in  cognitive diagnosis assessment. In particular, we present  in Section \ref{sec-dina} the sufficient and necessary   conditions for  strict and generic identifiability of $(Q,\TT,\pp)$ under the basic DINA model.  
 These identifiability results are also applicable to the DINO model \citep{Templin}, thanks to the duality between these two models  \citep{chen2015statistical}.
 Section \ref{sec-general}  presents the sufficient and necessary   conditions for   generic identifiability of $(Q,\TT,\pp)$ under general RLCMs, which include the popular GDINA and LCDM models.

%\noindent{\bf 3.2~ Identifiability of   $(Q, \TT, \pp)$ under the DINA model}
% \bigskip
%
%\noindent 

\section[DINA Model]{Identifiability  of   $(Q, \TT, \pp)$ under the DINA model}\label{sec-dina}

Under the DINA model,  
\cite{JLGXZY2011} first studied  identifiability of the $Q$-matrix under the assumption that the guessing parameters $\cg$ are known.
\cite{chen2015statistical} and \cite{xu2018jasa}  further proposed a set of sufficient conditions without assuming known item parameters.  
An important requirement  in these identifiability studies is  the completeness of the $Q$-matrix \citep{Chiu}.
 %A $Q$-matrix is said to be complete if  it can differentiate all latent attribute profiles, in the sense that under the $Q$-matrix, different attribute profiles have different response distributions. 
 Under the DINA model, the $Q$-matrix is said to be complete if
 %$\{\ee_k^\top:k=1,\ldots,K\}\subseteq \{\qq_j:j=1,\ldots,J\}$;
 %, equivalently, for each attribute there is some item which requires that and solely requires that attribute.  
it contains a $K\times K$ identity submatrix $I_K$ up to column permutation.
The previous studies in \cite{chen2015statistical} and \cite{xu2018jasa}  require  $Q$ to contain at least two complete submatrices $I_K$ for identifiability.

 However, it has been an open problem  what would be the minimal requirements on the $Q$-matrix for identifiability.
In the next theorem,  we solve this problem by providing the necessary and sufficient condition for identifiability of $(Q, \cs,\cg,\pp)$, under the earlier assumption that $p_{\aaa}>0$ for all $\aaa\in\{0,1\}^K$ \citep{Xu15,dina}. 

\begin{theorem}\label{thm-dina}
Under the DINA model, the following
Conditions $A$, $B$ and $C$ combined are necessary and sufficient for strict identifiability of $(Q, \cs,\cg,\pp)$.
\begin{itemize}
\item[$A$.] The true $Q$-matrix is complete. Without loss of generality, assume the   $Q$-matrix takes the following form
\begin{equation}\label{eq-true-q}
Q = \left(
\begin{array}{c}
I_K\\
Q^\star
\end{array}
\right).
\end{equation}

\item[$B$.] 
The column vectors of the sub-matrix $Q^\star$ in \eqref{eq-true-q} are distinct.

\item[$C$.] 
Each column in $Q$ contains at least three entries of ``~$1$''.
\end{itemize}
\end{theorem}

{In the Supplementary Material, we perform simulations to verify Theorem \ref{thm-dina}. In particular, see simulation study I for the verification of the sufficiency of the conditions $A$, $B$ and $C$ for joint identifiability; also see simulation studies III and IV regarding the necessity of the proposed conditions.}
We provide  comparisons of our Theorem 1 with some existing results.
First, although the same set of conditions $A$, $B$ and $C$ were also proposed in \cite{dina}, they assumed {\it a known $Q$} and  studied  identifiability of parameters $(\cs,\cg,\pp)$; on the contrary, Theorem 1 studies the joint identifiability of $(Q, \cs,\cg,\pp)$, which is theoretically much more challenging due to the unknown $Q$-matrix and therefore provides  a much stronger result than that in \cite{dina}. 
 In terms of estimation, Theorem \ref{thm-dina}   implies that one can consistently estimate both $Q$ and $(\cs,\cg,\pp)$, without worrying about a wrong $Q$-matrix would be indistinguishable from the true $Q$.
Second, Theorem \ref{thm-dina} also has much weaker requirements than the celebrated identifiability conditions resulting from three-way tensor decomposition \citep{kruskal1977three, allman}. Specifically, these classical results require the number of items $J\geq 2K+1$ for  (generic) identifiability.  
In contrast, conditions in Theorem \ref{thm-dina}  imply that we need  the number of items $J$ to be at least $K+\ceil*{\log_2(K)}+1$ under the DINA model. This is because other than the identity submatrix $I_K$, in order to satisfy Condition $B$ of \textit{distinctness}, the $Q$-matrix only needs to contain another $\log_2(K)$ items whose $K$-dimensional  $\qq$-vectors form a matrix with $K$ distinct columns. For example, for $K=8$, conditions in \cite{allman} require at least $2K+1=17$ items while our Theorem \ref{thm-dina} guarantees that the following $Q$ with $K+\log_2(K)+1=12$ items suffices for strict identifiability of $(Q,\cs,\cg,\pp)$ under DINA.
$$
Q=\begin{pmatrix}
& & &  I_8 & & &  &\\
0 & 0 & 1 & 1 & 1 & 0 & 1 & 1 \\
0 & 1 & 0 & 1 & 0 & 1 & 1 & 1 \\
1 & 0 & 0 & 0 & 1 & 1 & 1 & 1 \\
1 & 1 & 1 & 1 & 1 & 1 & 0 & 1
\end{pmatrix}.
$$

 Conditions $A$, $B$ and $C$ are the minimal requirements for  joint strict identifiability. 
When the true $Q$ fails to satisfy any of them, Theorem 1 implies that there must exist  $(Q,\cs,\cg,\pp)\nsim (\bar Q,\bar\cs,\bar\cg,\bar\pp)$ such that \eqref{eq-orig} holds.
In this scenario,   there are still cases where the model is ``almost surely" identifiable though not strictly identifiable,   as illustrated by the example under $Q_{4\times 2}$ in \eqref{eq-q24}; and on the other hand, there are also cases where the entire model is never identifiable, as shown in simulation studies III and IV in the Supplementary Material.
%, like Example \ref{exp-q} to appear later.
It is therefore desirable to study conditions that guarantee the former case, i.e.,  generic identifiability of $(Q,\cs,\cg,\pp)$.

In the following, we discuss necessity of  Conditions $A$, $B$, $C$ under the weaker notion of  generic identifiability. 
First,   Condition $A$ is necessary for  joint generic identifiability of  $(Q, \TT, \pp)$. If the true $Q$-matrix does not satisfy Condition $A$, then under DINA model, certain latent classes would be equivalent given $Q$, and their separate proportion parameters can never be identified, not even generically \citep{partial}. In certain scenarios where Condition $A$ fails, one can find a different $\bar Q$ that is not distinguishable from $Q$. {See simulation study IV in the Supplementary Material that illustrates the necessity of Condition $A$.}

Second, Condition $B$ is also difficult to relax and it serves as a necessary condition for generic identifiability when $K=2$.
Specifically, as shown in   \cite{dina}, when $K=2$,
  the only possible structure of the $Q$-matrix violating Condition $B$ while satisfying Conditions $A$ and $C$ is
$$
Q = \begin{pmatrix}
1 & 0  & 1 & \cdots & 1\\
0 & 1  & 1 & \cdots & 1\\
\end{pmatrix}^\top.
$$
And it is proved in \cite{dina} that for \textit{any} valid DINA parameters associated with this $Q$, there exist infinitely many different sets of DINA parameters that lead to the same distribution of the responses. Therefore the model is not generically identifiable.

Third,  differently from   Conditions $A$ and $B$, for generic identifiability, Condition $C$ can be  relaxed to certain extent. %that some attribute is required by only two items instead of three, but can not be further relaxed to the case of being required by only one item. 
The next theorem  characterizes how the $Q$-matrix structure in this case impacts generic identifiability. {For empirical verification of Theorem \ref{thm-dina-gen}, see simulation study II in the Supplementary Material.}

\begin{theorem}\label{thm-dina-gen}
Under the DINA model,  $(Q,\cs,\cg,\pp)$ are  not    generically identifiable if some attribute is required by only one item.

\noindent
If some attribute is required by only two items, suppose the  $Q$-matrix takes the following form after some column and row permutations, 
\begin{equation}\label{eq-prop1}
Q=\begin{pmatrix}
1 & \mathbf 0^\top\\
1 & \vv^\top\\
%1 & \mathbf 0^\top\\
\mathbf 0 & Q^\star
\end{pmatrix},
\end{equation}
where $\vv$ is a vector of length $K-1$ and $Q^\star$ is a $(J-2)\times(K-1)$ submatrix.
\begin{itemize}
\item[(a)]
If $\vv=\mathbf 1$,   $(Q,\cs,\cg,\pp)$ are  not  locally generically identifiable.
\item[(b)] If $\vv=\mathbf 0$,   $(Q,\cs,\cg,\pp)$ are  globally  generically identifiable if   
either 
\begin{itemize}
\item[(b.1)] the submatrix $Q^\star$ satisfies Conditions $A$, $B$ and $C$ in Theorem \ref{thm-dina}; or

\item[(b.2)] the submatrix $Q^\star$ has two submatrices $ I_{K-1}$.
\end{itemize}

\item[(c)]
If $\vv\neq \mathbf 0, \mathbf 1$,    $(Q,\cs,\cg,\pp)$ are  locally  generically identifiable if  $Q^\star$ satisfies Conditions $A$, $B$ and $C$ in Theorem \ref{thm-dina}.
\end{itemize}
\end{theorem}
 
\begin{remark}\label{rmk-local}
We say $(Q,\cs,\cg,\pp)$ are locally   identifiable, if in a neighborhood of the true parameters, there does not exist a different set of parameters that gives the same distribution of the responses. Local generic identifiability is a weaker notion than (global) generic identifiability, so the statement in part (a) of Theorem \ref{thm-dina-gen} also implies $(Q,\cs,\cg,\pp)$ are not globally generically identifiable.
\end{remark}

\begin{remark}\label{rmk-cons1}
 
In scenario (b.1) of Theorem \ref{thm-dina-gen}, the identifiable subset of the parameter space is 
%\begin{align*}%\label{eq-s1s2}
%
$\big\{(\cs,\cg,\pp): 
 \exists \aaa^{1}=(0,\alpha^1_2,\ldots,\alpha^1_K),\aaa^{2}=(0,\alpha^2_2,\ldots,\alpha^2_K)\in 
\{0\}\times\{0,1\}^{K-1},$ 
 such that   ${p_{\aaa^{1}}} {p_{\aaa^2+\ee_1}}\neq p_{\aaa^{2}}  p_{\aaa^1+\ee_1} \big\},
$
%\end{align*} 
where $\ee_j$ denotes the $J$-dimensional unit vector with the $j$th element being one and all the others   being zero.
In scenario (b.2) of Theorem \ref{thm-dina-gen}, we can write $Q=(I_K, I_K,(Q^{\star\star})^\top)^\top$ and the identifiable subset is 
%\begin{align*}%\label{eq-ratiok}
$\big\{(\cs,\cg,\pp):  
 \forall k\in\{1,\ldots,K\}, 
\exists \aaa^{k,1},\aaa^{k,2}\in \{0,1\}^{k-1}\times\{0\}\times \{0,1\}^{K-k-1}, 
$ such that  ${p_{\aaa^{k,1}}p_{\aaa^{k,2}+\ee_k}}\neq {p_{\aaa^{k,2}}p_{\aaa^{k,1}+\ee_k}}\big\}.$
%\end{align*}
The complements of these identifiable subsets in the parameter space give the non-identifiable subsets, which are both of measure zero in the   DINA model parameter space. 
 
\end{remark}

%We also point out that for generic identifiability, Condition $B$ \textit{can not be further relaxed} to that some attribute is required by only one attribute. In particular, if there exists some attribute only required by one item, for an arbitrary set of valid DINA parameters with $Q$, one can construct infinitely many sets of parameters associated with the same $Q$ that give the same distribution of the responses.
Next we give some discussions on generic identifiability of DINA model in the special case of $K=2$. 
%The only case not covered by Theorem \ref{thm-dina} and Theorem \ref{thm-dina-gen} is when the true $Q$-matrix takes the form of $Q_{4\times 2}$ as follows:
%Q_2 = \begin{pmatrix}
%1 & 0 \\
%0 & 1 \\
%1 & 0 \\
%0 & 1 \\
%0 & 1 \\
%\end{pmatrix}.
%\end{equation}
We have the following proposition.
\begin{proposition}\label{prop-dina-K2}
Under the DINA model with $K=2$ attributes,
%\begin{itemize}
%\item[(a)] When $J=4$,  the only $Q$-matrix structure not covered by Theorem \ref{thm-dina} and \ref{thm-dina-gen} takes the form 
%\begin{equation}\label{eq-q24}
%Q_{4\times 2} = \begin{pmatrix}
%1 & 0 \\
%0 & 1 \\
%1 & 0 \\
%0 & 1 \\
%\end{pmatrix},
%\end{equation}
%and $(Q_{4\times 2},\cc,\cs,\pp)$ are generically identifiable.
%\item[(b)] When $J\geq 5$, 
$(Q,\cs,\cg,\pp)$ are generically identifiable if and only if the conditions   in Theorem \ref{thm-dina} or \ref{thm-dina-gen}(b) hold.
%\end{itemize}
\end{proposition}
Proposition \ref{prop-dina-K2} gives a full characterization of joint generic identifiability when $K=2$, showing that the proposed generic identifiability conditions   are necessary and sufficient in this case.
%\end{remark}
The following example discusses all the possible $Q$-matrices with $K=2$ such that  $(Q,\cs,\cg,\pp)$  are not strictly identifiable, which proves Proposition \ref{prop-dina-K2} automatically.

\begin{example}
When $K=2$, the discussions on Conditions $A$ and $B$ before Theorem \ref{thm-dina-gen} show that $(Q,\cs,\cg,\pp)$ are not generically identifiable when  $A$ or $B$ is violated. So we only need to focus on the cases where Condition $C$  is violated while Conditions $A$ and $B$ are satisfied. Specifically, when $J\leq 5$, the $Q$-matrix could only take the following forms up to column and row permutations,
$$
%Q_1 = \begin{pmatrix}
%1 & 0\\
%0 & 1\\
%1 & 1\\
%\end{pmatrix},\quad
Q_1 = \begin{pmatrix}
1 & 0\\
0 & 1\\
1 & 1\\
0 & 1\\
\end{pmatrix},\qquad
Q_2 = \begin{pmatrix}
1 & 0\\
0 & 1\\
1 & 0\\
0 & 1\\
\end{pmatrix},\qquad
Q_3 = \begin{pmatrix}
1 & 0\\
0 & 1\\
1 & 0\\
0 & 1\\
0 & 1
\end{pmatrix}.
$$
%According to Theorem \ref{thm-dina},  none of   $(Q_i,\cs,\cg,\pp)$ for  $i=1,2,3$ are strictly identifiable. 
By Theorem \ref{thm-dina-gen}, $Q_1$ falls in scenario (a), so $(Q_1,\cs,\cg,\pp)$ are not   locally generically identifiable, i.e., even in a small neighborhood of the true parameters there exist infinitely many different sets of  parameters that give the same distribution of the responses. On the other hand, $Q_2$ falls in scenario (b.2) %with the submatrix $Q^\star$ containing two $I_{K-1}$ 
and $Q_3$ falls in  scenario (b.1), %with   $Q^\star$ satisfying Conditions $A$, $B$ and $C$, 
so   $(Q_2,\cs,\cg,\pp)$ and  $(Q_3,\cs,\cg,\pp)$ are both generically identifiable. 
In the case of $J>5$,   any $Q$ satisfying $A$ and $B$ while violating $C$ must contain one of the above $Q_i$ as a submatrix  and have some additional row vectors of $(0,1)$. By Theorem \ref{thm-dina-gen}, any such $Q$ extended from $Q_1$ is still not   locally generically identifiable, and any such $Q$ extended from $Q_2$ or $Q_3$ is globally generically identifiable. 
%This therefore proves the   result of Proposition \ref{prop-dina-K2}.
\end{example}

%\noindent{\bf 3.3~ Identifiability of $(Q, \TT, \pp)$ under general RLCMs}
% \bigskip

\section[General RLCMs]{Identifiability of $(Q, \TT, \pp)$ under general RLCMs}\label{sec-general}
 
%\blue{Add a discussion on strict identifiability first}
\noindent 
%Under general RLCMs, \cite{xu2018jasa} proposed a set of sufficient Conditions for strict identifiability of $(Q,\TT,\pp)$, and the Condition requires $Q$ to contain two identity submatrices $\mathcal I_K$.
%This requirement is strong and it can not be applied to identifying any $Q$ that lacks some single-attribute items. 
Since DINA is a submodel of the general RLCMs, Conditions $A$, $B$ and $C$ in Theorem \ref{thm-dina} are also necessary for strict identifiability of general RLCMs. 
For instance, our proposed Conditions $A$, $B$ and $C$ are weaker than  the sufficient conditions proposed by \cite{xu2018jasa} for strict identifiability of $(Q,\TT,\pp)$ under general RLCMs; and if their conditions are satisfied, the current conditions $A$, $B$ and $C$ are also satisfied. 
%and the condition requires $Q$ to contain two identity submatrices $\mathcal I_K$.
However, these necessary requirements may be strong in practice and they can not be applied to identifying any $Q$ that lacks some single-attribute items (i.e., lacks some unit vector as a row vector). 
A natural question is whether  Conditions $A$, $B$ and $C$ can be  relaxed under the  weaker notation of of generic identifiability.
This section addresses this question.
%Additionally, it is unknown what structural constraint of the $Q$-matrix is necessary for identifiability. In this section we address the problems with the goal of generic identifiability.

 Under general RLCMs, the next theorem shows that Condition $C$ (each attribute is required by at least three items)  is necessary for generic identifiability of $(Q,\TT,\pp)$, contrary to the results for the   DINA model where Conditions $A$ and $B$ can not be relaxed while Condition $C$ can. See simulation studies VI and VII in the Supplementary Material for the verification of Theorem \ref{prop-gdina-nece}.

\begin{theorem}\label{prop-gdina-nece}
Under a general RLCM, Condition $C$ in Theorem 1 is  necessary  for generic identifiability of $(Q,\TT,\pp)$.
%\begin{itemize}
%\item[($B$)] Each attribute is required by at least three items.
%\end{itemize}
Specifically, when the true $Q$-matrix violates $C$,  for any   model parameters $(\TT,\pp)$ associated with $Q$, there exist infinitely many   sets of     $(\bar Q, \bar \TT,\bar \pp)\nsim (Q,\TT,\pp)$  such that equation \eqref{eq-orig} holds, making $(Q,\TT,\pp)$ not generically identifiable.

\end{theorem}

While Condition $C$ is necessary, we next show that the other two conditions $A$ and $B$ can be further relaxed for generic identifiability of general RLCMs. Before stating the result, we first introduce a new concept about the $Q$-matrix, the \textit{generic completeness}.

\begin{definition}[Generic Completeness]
A $Q$-matrix with $K$ attributes is said to be generically complete, if after some column and row permutations, it has a $K\times K$ submatrix with all diagonal entries being ``1". %We call such a submatrix a generically complete component.
\end{definition}
Generic completeness is a relaxation  of the concept of completeness. In particular, a $Q$-matrix is generically complete, if up to column and row permutations, it contains a submatrix as follows:
$$
\begin{pmatrix}
    1 & * & \dots  & * \\
    * & 1 & \dots  & * \\
    \vdots & \vdots & \ddots & \vdots \\
    * & * & \dots  & 1
\end{pmatrix},
$$ where the off-diagonal   entries ``$*$"  are left unspecified.  
Note that any complete $Q$-matrix is also generically complete, while a generically complete $Q$-matrix may not have any single attribute item.

With the concept  of generic completeness,   the next theorem gives sufficient conditions for joint generic identifiability, and shows  that under general RLCMs, the necessary conditions $A$ and $B$ for strict identifiability are not necessary any more in the current setting. %for generic identifiability. %and they can be weakened to the following Conditions $D$ and $E$. 
\begin{theorem}\label{thm-multi}
Under a general RLCM, if the true $Q$-matrix satisfies the following Conditions D and E, then  $(Q,\TT,\pp)$ are {generically} identifiable.

\begin{enumerate}
\item[D.]  The $Q$-matrix has two  nonoverlapping  generically complete $K\times K$ submatrices $Q_1$ and $Q_2$. Without loss of generality, assume the $Q$-matrix is in the following form
\begin{equation}\label{eq-diag}
Q = \left(\begin{array}{c}
Q_1 \\
Q_2 \\
Q^\star
\end{array}\right)_{J\times K}.
%\quad
%Q_i =
%\left(\begin{array}{cccc}
%    1 & * & \dots  & * \\
%    * & 1 & \dots  & * \\
%    \vdots & \vdots & \ddots & \vdots \\
%    * & * & \dots  & 1
%\end{array}\right)_{K\times K} \mbox{for } i=1,2,
\end{equation}
\item[E.] Each column of the submatrix  $Q^\star$ in \eqref{eq-diag} contains at least one entry of ``$1$''.
\end{enumerate}
\end{theorem}

\begin{remark}\label{rmk-cons2}
Under Theorem \ref{thm-multi}, the identifiable subset of the parameter space is
$\{(\TT,\pp): 
  \det (T(Q_1,\TT_{Q_1}))\neq 0$, $\det (T(Q_2,\TT_{Q_2}))\neq 0,$  and $T(Q^\star,\TT_{Q^\star})\cdot\text{Diag}(\pp)$ has distinct column vectors\}.  Its complement is the non-identifiable subset and it has measure zero in the parameter space $\boldsymbol \vartheta_Q$, when $Q$ satisfies Conditions $D$ and $E$.
Please see the supplementary materials for the definition of the $T$-matrices ($T(Q_1,\TT_{Q_1})$, etc.).
\end{remark}

\begin{remark}\label{rmk-thm4}
The proof of   Theorem \ref{thm-multi} is based on the proof of Theorem 7 in \cite{partial}, which  proposed the same Conditions $D$ and $E$ as sufficient conditions for generic identifiability of model parameters given a known $Q$. 
We point out that though $D$ and $E$ serve as sufficient  conditions for generic identifiability both when  $Q$ is known and when $Q$ is unknown, the generic identifiability results in these two scenarios are different.
In particular, Theorem 8 in \cite{partial} shows that when $Q$ is known, some attribute can be required by only two items for generic identifiability (i.e., Condition $C$   can be relaxed); while our current Theorem \ref{prop-gdina-nece} shows that when $Q$ is unknown, Condition $C$ indeed becomes necessary. 
%Specifically, if $Q$ is unknown and Condition $C$ is violated, there exist multiple different $Q$-matrices leading to identical distribution of the responses. 
\end{remark}

The proposed sufficient Conditions $D$ and $E$ weaken the strong requirement of Conditions $A$ and $B$, especially the identity submatrix requirement that may be difficult to satisfy in practice.{See simulation study V in the Supplementary Material for the verification of Theorem \ref{thm-multi}.}
 Note that  Conditions $D$ and $E$ imply the necessary  Condition $C$ that each attribute is required by at least three items.  
%In addition, Condition $D$   requires $Q$ to contain two generically complete submatrices. 

We next discuss the necessity of Conditions $D$ and $E$.  As shown in Section 3.2, under DINA,  the completeness of $Q$ is necessary for joint strict identifiability of $(Q,\cs,\cg,\pp)$.
 For general RLCMs, we have an analogous conclusion that the generic completeness of $Q$, which is part of Condition $D$, is necessary for joint generic identifiability of $(Q,\TT,\pp)$. This is stated in the next theorem.
\begin{theorem}\label{prop-gencom}
Under a general RLCM, generic completeness of the $Q$-matrix is necessary for joint generic identifiability of $(Q,\TT,\pp)$.
\end{theorem}

Furthermore, we show  that Conditions $D$ and $E$ themselves are in fact necessary when $K=2$, indicating the difficulty of further relaxing them.
\begin{proposition}\label{prop-gdina-K2}
For a general RLCM with  $K=2$,   Conditions $D$ and $E$  are   necessary and sufficient  for  generic identifiability of $(Q,\TT,\pp)$.
\end{proposition}
%\blue{give some specific examples for K=2}

We use the following example to illustrate the   result of Proposition \ref{prop-gdina-K2}, which also gives a natural proof of the proposition.
\begin{example}
When $K=2$,  %(because a $Q$ with some all-zero rows is equivalent to the submatrix of it containing those not-all-zero rows), 
 a $Q$-matrix which satisfies the necessary Condition $C$ but   not   Conditions  $D$ or $E$ can only take   the following form  $Q_1$ or $Q_2$, up to row permutations,
$$
Q_1 = \begin{pmatrix}
1 & 1 \\
1 & 1 \\
1 & 1
\end{pmatrix},\qquad
Q_2=\begin{pmatrix}
1 & * \\
* & 1 \\
1 & 1 \\
1 & 1 \\
\end{pmatrix};
\qquad
\bar Q_2=\begin{pmatrix}
1 & 1 \\
1 & 1 \\
1 & 1 \\
1 & 1 \\
\end{pmatrix}.
$$
The ``$*$"s in $Q_2$ are unspecified values, either 0 or 1.
For $Q_1$ with $J=3$, $K=2$ and any     parameters $(\TT,\pp)$, there are $2^J=8$ constraints in \eqref{eq-orig} for solving $(\bar\TT,\bar\pp)$ under $Q_1$ itself, while the number of free parameters of $(\bar\TT,\bar\pp)$ is 
 $
 \lvert \{p_{\aaa}:\aaa\in\{0,1\}^2\}\cup\{\theta_{j,\aaa}:j\in\{1,2\},\aaa\in\{0,1\}^2\}  \rvert = 2^K + 2^K\times J = 16>8.
$  
For $Q_2$ with $J=4$, $K=2$ and any associated $(\TT,\pp)$, there are $2^J=16$ constraints in \eqref{eq-orig} for solving $(\bar\TT,\bar\pp)$, while the number of free parameters of $(\bar\TT,\bar\pp)$ under the alternative $\bar Q_2$ is $2^K+J\times2^K = 20>2^J=16$. In both cases there are infinitely many sets of solutions of \eqref{eq-orig} as alternative model parameters, so neither $(Q_1,\TT,\pp)$ nor  $(Q_2,\TT,\pp)$ are generically identifiable.
\end{example}

\section[Discussion]{Discussion}\label{sec-sum}
In this work, we study the identifiability issue of RLCMs with unknown   $Q$-matrices.
For the basic DINA model, we derive the necessary and sufficient conditions for  strict joint identifiability of the $Q$-matrix and the associated model parameters. 
We also study a slightly weaker identifiability notion, generic identifiability, and propose sufficient and necessary conditions for it under the DINA model and general RLCMs, respectively.

%\begin{figure}[h!] 
%\centering
%\resizebox{11cm}{!}{%
%\begin{tikzpicture}[FlowChart, start chain = A going below]
%
%\node [input0] (input) 
%{$Q$-matrix with $K=2$};                
%   
%   
%\node [output2] (two-p) 
%[node distance=1.5cm and 0cm,below left=of input] 
%{\textbf{Proposition \ref{prop-dina-K2}}\\
%necessary \& sufficient\\ for Generic Identifiability
%};
%
%\node [output2]  (multi-p) 
%[node distance=1.5cm and 0cm, below right=of input] 
%{\textbf{Proposition \ref{prop-gdina-K2}}\\
%necessary \& sufficient\\ for Generic Identifiability};
%   
%\draw[->] (input) -- ($(input.south)+(0,-0.75)$) -| (two-p) node[above,pos=0.25] {two-param. model} ;
%\draw[->] (input) -- ($(input.south)+(0,-0.75)$) -| (multi-p) node[above,pos=0.25] {multi-param. model};
%  
%\end{tikzpicture}
%}
%\caption{\blue{Full characterization of identifiability for the special case  $K=2$.}}
%\label{fig-flow2}
%\end{figure}

\paragraph{Statistical consequences of identifiability.}
In the setting of RLCMs, identifiability naturally leads to estimability, in different senses under strict and generic identifiability. If the $Q$-matrix and the associated model parameters are strictly identifiable, then $Q$ and model parameters can be jointly estimated from data consistently.
If the $Q$-matrix and model parameters are generically identifiable, then for true parameters ranging almost everywhere in the parameter space with respect to the Lebesgue measure, the $Q$-matrix and model parameters can be jointly estimated from data consistently.

As pointed out by one reviewer, the analysis of identifiability is under an ideal situation with an infinite sample size. Indeed, general identification problems assume the hypothetical exact knowledge of the distribution of the observed variables, and ask under what conditions one can recover the underlying parameters \citep{allman2009identifiability}.
%Identifiability is also very related to the estimability of model parameters.
Next we discuss the finite sample estimation issue under the proposed identifiability conditions for strict identifiability, following a similar argument as Proposition 1 in \cite{xu2018jasa}.
Denote the true $Q$-matrix and model parameters by $Q^0$ and $\bo\eta^0 = (\TT^0,\pp^0)$. 
Consider a sample with $N$ i.i.d. response vectors $\RR_1,\RR_2,\ldots,\RR_N$,
and denote the log-likelihood of the sample by $\ell(\TT,\pp)=\sum_{i=1}^N \log \mathbb P(\RR_i\mid Q,\TT,\pp)$.
Under a specified RLCM, a $Q$-matrix determines the  structure of the item parameter matrix $\TT$, by specifying which entries of it are equal. For a given $\TT$, we can define an equivalent formulation of it, a sparse matrix $\bo B$ having the same size as $\TT$, in the following way. Under a general RLCM such as the GDINA model in Example \ref{exp-gdina}, the item parameters can be parameterized as  $\theta_{j,\aaa}=\sum_{\mathcal S\subseteq \{1,\ldots,K\}}\beta_{j,\mathcal S}\prod_{k\in \mathcal S}\alpha_k$.
 Based on this, we define the $j$th row of $\bo B$ as a $2^K$-dimensional vector collecting all these $\beta$-coefficients; that is, $\bo B_j=(\beta_{j,0},\beta_{j,1},\ldots, \beta_{j,K},\ldots, \beta_{j,12\cdots K})$.  Then as long as the $\bq$-vector $\qq_j\neq \mathbf 1_K$, the vector $\bo B_j$ would be ``sparse" and so is the matrix $\bo B$.
%Namely, the  $(j,\aaa)$-th entry of matrix $\bo B$ are the difference between $\theta_{j,\aaa}$ and the highest level item parameter of item $j$.
%Then equivalently, the $Q$-matrix determines a sparsity pattern of the matrix $\bo B$. 
 For the true $Q^0$, we denote the corresponding $\bo B$-matrix by $\bo B^0$.
Under a specified RLCM (such as DINA or GDINA), the identification of $Q^0$ is then implied by the identification of the indices of nonzero elements of $\bo B^0$. 
Denote  the support of the true $\bo B^0$ and any candidate $\bo B$ by $S_0$ and $S$, respectively.
Define 
$
	C_{\min}(\eeta^0) = \allowbreak \inf_{\{S\neq S_0,\, |S|\leq |S_0|\}} \allowbreak (|S_0\setminus S|)^{-1} h^2(\bo\eta^0,\bo\eta),
$
where $h^2(\bo\eta^0,\bo\eta)$ denotes the Hellinger distance between the two distributions of $\bo R$ indexed by parameters $\bo\eta^0$ under the true $\bo B^0$  and  $\bo\eta$ under the candidate $\bo B$.
Denote the $Q$-matrix and the model parameters that maximize the log-likelihood $\ell(\TT,\pp)$ subject to the $L_0$ constraint $|S|\leq |S_0|$ by $\widehat{\bo\eta}=(\widehat\TT,\widehat\pp)$, and denote the ``oracle" MLE of model parameters obtained assuming $Q^0$ is known by $\widehat\eeta^0=(\widehat\TT^0,\widehat\pp^0)$.
Then we have the following finite sample error bound for the estimated $Q$-matrix and model parameters.

\begin{proposition}\label{prop-finite}
Suppose $Q^0$ satisfies the proposed sufficient conditions for joint strict identifiability, %and $(\TT^0,\pp^0)$ do not fall in the measure zero non-identifiable subset of the parameter space. 
then $C_{\min}(\TT^0,\pp^0)\geq c_0$ for some positive constant $c_0$. Furthermore,
	\begin{equation}
	\label{eq-bound}
	\mathbb P(\widehat Q \not\sim Q^0)
	\leq \mathbb P(\widehat {\bo \eta}\neq \widehat\eeta^0)
	\leq c_2\exp\{-c_1 N C_{\min}(\TT^0,\pp^0) 
	\},
    \end{equation}
where $c_1,c_2>0$ are some constants. Namely, when joint strict identifiability conditions hold, the finite sample estimation error has an exponential bound.
%\begin{itemize}
%	\item[(a)] Suppose
%
%    \item[(b)] Suppose 
%\end{itemize}
\end{proposition}
Proposition \ref{prop-finite} shows that the estimation error decreases exponentially in $N$ if the model is   identifiable. 
On the other hand, when the   identifiability conditions fail to hold, there exist alternative models that are close to the true model in terms of the Hellinger distance. This would make the $C_{\min}(\TT^0,\pp^0)$  in \eqref{eq-bound} equal to zero, instead of bounded away from zero as shown in Proposition \ref{prop-finite}. %The upper bound for the probability of falsely estimating the MLE becomes larger accordingly in this scenario. 
Therefore, the finite sample error bound in \eqref{eq-bound} becomes $O(1)$ in this non-identifiable scenario. 
In particular, in the case where generic identifiability conditions are satisfied, 
  $C_{\min}(\TT^0,\pp^0)$ depends on the distance between the true parameters and the  non-identifiable measure-zero subset of the parameter space; as the true parameters become  closer to this measure-zero set, $C_{\min}(\TT^0,\pp^0)$ decreases to   zero and  
  a larger sample size therefore may be needed to achieve a prespecified   level of estimation accuracy.
%This also partly explains the phenomenon in Figure \ref{fig-Q24-ratio}.

%Under strict identifiability, the convergence rate will be impacted by how close the model parameters lie to the boundary of the parameter space; while under generic identifiability, the convergence rate will also be impacted by how close the model parameters lie to the measure-zero non-identifiable subset of the parameter space.

\paragraph{Potential extensions to other latent variable models.}
%The techniques used to establish identifiability theory in this paper might be extended to develop identifiability conditions for other latent variable models, 
We briefly discuss the potential extensions of the developed theory to some other latent variable models, such as restricted latent class models with ordinal polytomous attributes \citep{davier2008general, ma2016seq, chen2018poly}, and multidimensional latent trait models \citep{embretson1991}. 
   First, an RLCM with ordinal polytomous attributes can be considered as an RLCM with binary attributes and a constrained relationship among the binary attributes. For instance, consider an ordinal attribute $\gamma$ that can take $C$ different values $\{0,1,\ldots,C-1\}$, then $\gamma$ can be equivalently viewed as a collection of $C-1$ binary random variables $\aaa^{\gamma}:=(\alpha_1,\ldots,\alpha_{C-1})$ with the following constraints. If $\alpha_i=0$ for some $i<C-1$, then $\alpha_j=0$ for all $j=i+1,\ldots,C-1$. In other words, any pattern $\aaa^{\gamma}$ with $\alpha_i=0$ and $\alpha_j=1$ for some $i<j$ is ``forbidden" and constrained to have proportion zero. The vector of the polytomous attributes can be augmented to a longer vector of binary attributes with constraints in this fashion. Then we can consider the restricted latent class model with the augmented proportion parameters, by constraining the proportions of those ``forbidden" binary attribute patterns to zero. In this scenario, it might be possible to extend the current theory and develop identifiability conditions for the case of polytomous attributes.

Second, if a multidimensional latent trait model includes both continuous and discrete latent traits, then the techniques of establishing identifiability for latent class models in this paper would also be useful when treating the discrete latent variables. For the continuous latent variables, the  techniques  developed in \cite{bai2012} for identifiability of the factor analysis model and those developed in traditional multivariate analysis \citep{anderson} would be helpful. 

\bigskip
In practice, the newly developed identifiability theory can serve as the foundation for designing statistically guaranteed estimation procedures. Specifically, consider the set of all $Q$-matrices that satisfy our identifiability conditions ($A$, $B$ and $C$ under the DINA model, or $D$ and $E$ under general RLCMs), and call it the ``identifiable $Q$-set". Then one can use likelihood-based approaches, such as that in \cite{xu2018jasa}, to jointly estimate $Q$ and model parameters by constraining $Q$ to the identifiable $Q$-set; or one can use Bayesian approaches to estimate $Q$ such as that in \cite{chen2018bayesian}.
Additionally, if under the DINA model   the $Q$-matrix does not contain a submatrix $I_K$, then according to \cite{partial}, certain attribute profiles would be equivalent  and the strongest possible identifiability argument therein is the so-called $\pp$-partial identifiability. In this scenario, it would be interesting to study the identifiability   of the incomplete $Q$-matrix under the notion of $\pp$-partial identifiability, and we leave this for future study.

\vskip 14pt
\noindent {\Large\bf Supplementary Materials}

The online supplementary material contains the proofs of Proposition 1, 4, Theorems 1--5, and additional simulation results.
\par
%%%%%%%%%%%%%%%%%%%%%%%%%%%%%%%%%%%%%%%%%%%%%%%%%%%%%%%%%%%%%%%%%%%%%%%%%%%%%%%%%%%%%%%%%%%%%%%%%%%%%%%%%%%%%%%%%%%%%%%%%%%%
\vskip 14pt
\noindent {\Large\bf Acknowledgements}

This work was supported in part by National Science Foundation grants SES-1659328 and DMS-1712717.
\par

\newpage
\centerline{\Large\bf Supplementary Material}

\vspace{3mm}
In this supplementary material, we present proofs of the theorems in the main text and also provide simulation results to support them. We first introduce some notations and a useful technical lemma. We then give proofs of Proposition \ref{prop-0vector}, Theorems 1--5, and Proposition \ref{prop-finite} in Sections S1-S7, respectively. 
{We perform various simulation studies to verify the proposed theoretical conditions in Section S8.}
\\
\setcounter{section}{0}
\setcounter{equation}{0}
\def\theequation{S\arabic{section}.\arabic{equation}}
\def\thesection{S\arabic{section}}

Before presenting the proofs of the theoretical results, we introduce a useful notation, the $T$-matrix $T(Q,\TT)$ of size $2^J\times 2^K$.
%Here $\TT$ denotes the item parameters for any restricted latent class model in general, and for the DINA model $\TT$ are just the slipping and guessing parameters, i.e.,  $\TT=(\cs,\cg)$.
The rows of $T(Q,\TT)$ are indexed by the $2^J$ different response patterns $\rr = (r_1,\ldots,r_J)^\top\in\{0,1\}^J$, and columns by attribute patterns $\aaa\in\{0,1\}^K$, while the $(\rr, \aaa)$th entry of $T(Q,\TT)$, denoted by $T_{\rr, \aaa}(Q,\TT)$, represents the marginal probability that subjects in latent class $\aaa$ provide positive responses to the set of items $\{j: r_j=1\}$, namely 
\[
T_{\rr,\aaa}(Q,\TT) = P(\boldsymbol R \succeq \rr\mid Q,\TT, \aaa) = \prod_{j=1}^J \theta_{j,\aaa}^{r_j}.
\]
We denote the $\aaa$th column vector and the $\rr$th row vector of the $T$-matrix by $T_{\Cdot,\aaa}(Q,\TT)$ and $T_{\rr,\Cdot}(Q,\TT)$ respectively. Let $\ee_j$ denote the $J$-dimensional unit vector with the $j$th element being one and all the other elements being zero, then any response pattern $\rr$ can be written as a sum of some $\ee$-vectors, namely $\rr=\sum_{j:r_j=1}\ee_j$.
%%%
The $\rr$th element of the $2^J$-dimensional vector $T(Q,\TT)\pp$ is
\[
T_{\rr,\Cdot}(Q,\TT)\pp = \sum_{\aaa\in\{0,1\}^K} T_{\rr,\aaa}(Q,\TT) p_{\aaa} = P(\boldsymbol R \succeq \rr \mid Q,\TT,\pp).
\]
Based on the $T$-matrix, there is an equivalent  definition of identifiability for ($Q, \TT, \pp$). The $T$-matrix also has a nice property that will be useful in proving the identifiability results. These are summarized in the following lemma, whose proof can be found in \cite{Xu2016}.

%equivalent to definition  \eqref{eq-orig} in Section 2.3 of the main text. The equivalence of the two definitions comes from that two sets of model parameters lead to the same marginal distribution of  responses $\{P(\boldsymbol R \succeq \rr \mid \TT), \forall \rr\in\{0,1\}^J\}$ if and only if they  lead to the same distribution of the responses $\{P(\boldsymbol R = \rr \mid \TT), \forall \rr\in\{0,1\}^J\}$.

\begin{lemma}\label{lem}
Under a restricted latent class model, $(Q,\TT, \pp)$ are identifiable if and only if for any $(Q,\TT, \pp)$ and  $(\bar Q,\bar\TT, \bar\pp)$, 
\begin{equation}\label{eq1}
T(Q,\TT)\pp = T(\bar Q, \bar\TT)\bar\pp%\quad \text{implies}\quad (\TT,\pp) = (\bar\TT,\bar\pp).
\end{equation}
implies $(Q,\TT,\pp) = (\bar Q,\bar\TT,\bar\pp)$.
For any $\ttt^*=(\theta_1,\ldots,\theta_J)^\top\in\mathbb R^J$, there exists an invertible matrix $D(\ttt^*)$ depending only on $\ttt^*$, such that %the diagonal elements of $D(\ttt^*)$ are all 1, %e.g., $diag\{D(\ttt^*)\}=\mathbf1$,  
\begin{equation}\label{eq-trans}
T(Q,\TT-\ttt^*\mathbf1^\top) = D(\ttt^*)T(Q,\TT).
\end{equation}
\end{lemma}

%Another useful property of the $T$-matrix  is given by the following lemma, whose proof is given in Section D.
%\begin{lemma}\label{lem-rk}
%Denote the $T$-matrix corresponding to a subset of items $S$ by $T(\TT_{S})$, where $\TT_S=(\theta_{j,\aaa},~j\in S,~\aaa\in\ma)$. 
%If for an item set $S$, the $\Gamma$-matrix $\Gamma^S$ of size $|S|\times m$ is separable, then the corresponding $T$-matrix $T(\TT_{S})$ of size $2^{|S|}\times m$  has full column rank $m$.
%\end{lemma}

\bigskip
We introduce some additional notations. For a submatrix $Q_1$ of $Q$ that has size $J_1\times K$,
we denote the item parameter matrix corresponding to these $J_1$ items by $\TT_{Q_1}$, then $\TT_{Q_1}$ is a $J_1\times K$ submatrix of $\TT$.
Denote $Q_1$'s corresponding $T$-matrix by $T(Q_1,\TT_{Q_1})$, then $T(Q_1,\TT_{Q_1})$ has size $2^{J_1}\times 2^K$.
For notational simplicity, in the following we denote $\cc\equiv\one-\cs$ under the DINA model, then $\TT = (\mathbf1-\cs,\cg)=(\cc,\cg)$ under DINA.

We   add some remarks on Lemma \ref{lem}. 
First,  Equation \eqref{eq1} can be written as that, for any response pattern $\rr\in\{0,1\}^J$,
%\begin{equation*}%\label{eq-tpr0}
$T_{\rr,\Cdot}(Q,\TT)\pp = T_{\rr,\Cdot}(\bar Q,\bar\TT)\bar\pp$.
Second, 
thanks to \eqref{eq-trans}, for any $\ttt^*=(\theta_1,\ldots,\theta_J)^\top\in\mathbb R^J$, equality \eqref{eq1} leads to
$$
T(Q,\TT-\ttt^*\mathbf1^\top)\pp = T(\bar Q,\bar\TT-\ttt^*\mathbf1^\top)\bar\pp,
$$
and further $T_{\rr,\Cdot}(Q,\TT-\ttt^*\mathbf1^\top)\pp = T_{\rr,\Cdot}(\bar Q,\bar\TT-\ttt^*\mathbf1^\top)\bar\pp$ for any $\rr\in\{0,1\}^J$.
Besides, If \eqref{eq1} holds, then for any submatrix $Q_1$ of $Q$, equality $T(Q_1,\TT_{Q_1})\pp \allowbreak = T(\bar Q_1,\bar\TT_{\bar Q_1})\bar\pp$ also holds.

\section{Proof of Proposition 1}
Consider a $Q$-matrix of size $J\times K$ in the form
\begin{equation*}%\label{eq-qq}
Q=\begin{pmatrix}
Q'\\
\mathbf0\\
\end{pmatrix},\end{equation*}
where $Q'$ is of size $J'\times K$ and contains those nonzero $\bq$-vectors of $Q$.
For any item   $j\in\{J'+1,\ldots,J\}$ which has   $\bq_j=\mathbf 0$, all the attribute profiles $\aaa$ satisfy $\aaa\succeq\qq_j$, so there is only one item parameter associated with $j$ under $Q$, and we denote it by $\theta_j$. 
%Firstly, if the $\bq$-vector of some item $j$ equals the zero vector in the $Q$-matrix, then all the attribute profiles have the same positive response probability $1-s_j$ (or $1-\bar s_j$) to this item, and the guessing parameter of this item is not needed in model specification. 
%\end{equation*}
Denote the first $J'$ rows of $\TT$ by $\TT'$.
% and similarly denote $\bar\TT'$. 
Denote the $2^{J'}\times 2^K$ $T$-matrix associated with matrix $Q'$ by $T'(Q',\TT')$.

First consider the case where $(Q',\TT',\pp)$ are strictly  (or generically) identifiable, and we will show $(Q,\TT,\pp)$ are also strictly  (or generically) identifiable.
Assume there is a $J\times K$ matrix $\bar Q$ and associated parameters $(\bar\TT,\bar\pp)$ such that \eqref{eq1} holds. Denote the submatrix of $\bar Q$ containing its first $J'$ rows by $\bar Q'$, and the submatrix of $\bar\TT$ containing its first $J'$ rows by $\bar\TT'$.
Then \eqref{eq1} implies $T(Q',\TT')\pp' = T(\bar Q',\bar \TT')\bar \pp'$, and the strict (or generic) joint identifiability of $(Q',\TT',\pp)$ gives that $\bar Q'\sim Q'$ and $(\bar\TT',\bar\pp)= (\TT',\pp)$.
% (or if $(Q',\TT')$ satisfy certain generically satisfied constraints). 
%We consider the DINA model and the general RCLMs separately.
For an arbitrary RLCM, the strict (or generic) identifiability of $(Q',\TT',\pp)$ implies that $T(Q',\TT')$ has full rank $2^K$ strictly (or generically). 
%This is because as mentioned in the main text, condition $A$, i.e., the completeness condition, is necessary for strict joint identifiability. Then without loss of generality assume the first $K$ rows of $Q'$ and $\bar Q'$ are $I_K$. Now that the first $K$ items are all single attribute items, each item $k\in\{1,\ldots,K\}$ only have two item parameters, which we denote by $\theta_{k,\mz}$ and $\theta_{k,\mo}$.
%Next consider any item $j\in\{J'+1,\ldots,J\}$, and define $\rr^*=\sum_{k=1}^K \ee_k$ and 
%$$
%\ttt^* = \sum_{1\leq k\leq K:\, \aaa\succeq\qq_k} \theta_{k,\mz}\ee_k + \sum_{1\leq k\leq K:\, \aaa\nsucceq\qq_k} \theta_{k,\mo}\ee_k.
%$$
%Then \eqref{eq-trans} implies that
%$$
%\theta_{j,\aaa} = 
%\frac{T_{\rr^*+\ee_j,\Cdot}(Q,\cc-\ttt^*, \cg-\ttt^*)\pp}{T_{\rr^*,\Cdot}(Q,\cc-\ttt^*, \cg -\ttt^*)\pp} = 
%\frac{T_{\rr^*+\ee_j,\Cdot}(\bar Q,\bar\cc-\ttt^*, \bar\cg-\ttt^*)\bar\pp}{T_{\rr^*,\Cdot}(\bar Q,\bar\cc-\ttt^*, \bar\cg -\ttt^*)\bar\pp}
%= \bar\theta_{j,\aaa}.
%$$
%%%%%%%%%
%so if $(\TT',\pp')$ satisfy certain generically true constraints, $\bar Q'\sim Q'$ and $(\bar\TT',\bar\pp)= (\TT',\pp)$. 
This is because if not so, then the proportion parameters $\pp$ can not be strictly (or generically) identifiable, in the sense that there exist multiple different $\pp$ such that $T(Q', \TT')\pp$ are all equal. This would contradict the assumption that $(Q',\TT',\pp)$ are strictly (or generically) identifiable. Therefore $T(Q', \TT')$  is strictly (or generically) full-rank.
%\textcolor{blue}{yet to be completed...
Then for each $\aaa\in\{0,1\}^K$ there must exist a $2^K$-dimensional vector $\vv_{\aaa}$ such that 
$$\vv_{\aaa}^\top\Cdot T(Q', \TT') = \vv_{\aaa}^\top\Cdot T(\bar Q', \bar \TT') = (\mathbf0, \underbrace{x_{\aaa}}_{\text{column }\aaa},\mathbf 0),\quad x_{\aaa}\neq 0,$$
and $\vv_{\aaa}^\top\Cdot T( Q',  \TT') \pp = \vv_{\aaa}^\top\Cdot T(\bar Q', \bar \TT') \bar\pp = x_{\aaa} p_{\aaa}\neq 0$.
Then again use the property \eqref{eq-trans} and we have the following equality for any $j\in\{J'+1,\ldots,J\}$,
\begin{align*}
\theta_{j,\aaa} = 
&\frac{\{T_{\ee_j,\Cdot}(Q,\TT) \odot [\vv_{\aaa}^\top\Cdot~ T(Q', \TT')]\,\}\pp }{\vv_{\aaa}^\top\Cdot~ T(Q', \TT')\pp}
\\
= &\frac{\{T_{\ee_j,\Cdot} (Q,\TT) \odot [\vv_{\aaa}^\top\Cdot~ T(\bar Q', \bar \TT')]\,\}\bar\pp }{\vv_{\aaa}^\top\Cdot~ T(\bar Q', \bar \TT')\bar\pp}
= \bar\theta_{j,\aaa},
\end{align*}
where ``$\odot$" represents the element-wise product of two vectors. This   proves $\TT = \bar\TT$ and $Q\sim \bar Q$. So $(Q,\TT,\pp)$ are   strictly (or generically) identifiable.

Next consider the case where $(Q',\TT',\pp)$ are \textbf{not} strictly (or generically) identifiable, so there exist $(\bar Q',\bar\TT',\bar\pp) \nsim (Q',\TT',\pp)$ such that $T'(\bar Q',\bar\TT')\bar\pp = T'(Q',\TT')\pp$. Now extend $\bar Q'$ to $\bar Q$ of size $J\times K$ by adding $J-J'$ all-zero $\qq$-vectors, i.e.,
$$
\bar Q=\begin{pmatrix}
\bar Q'\\
\mathbf 0
\end{pmatrix},
$$
and set $\bar\theta_j=\theta_j$ for $j\in\{J'+1,\ldots,J\}$.
Then for any $\rr=(r_1,\ldots,r_{J'},\allowbreak r_{J'+1},\ldots,r_J)\in\{0,1\}^J$ and the corresponding $\rr'=(r_1,\ldots,r_{J'})$,
\begin{align*}
T_{\rr,\Cdot}(Q,\TT)\pp &%\sum_{\aaa\in\{0,1\}^K}
%\Big\{T_{(\rr',\mathbf0),\Cdot}( Q,\TT)\pp\Big\} \prod_{j>J'}\theta_j^{r_j}
=\Big\{T'_{\rr',\Cdot}(Q',\TT')\pp\Big\} \prod_{j>J'}\theta_j^{r_j};\\
T_{\rr,\Cdot}(\bar Q,\bar\TT)\bar\pp &%\sum_{\aaa\in\{0,1\}^K}
%\Big\{T_{(\rr',\mathbf0),\Cdot}( \bar Q,\bar\TT)\bar\pp\Big\} \prod_{j>J'}\theta_j^{r_j}
=\Big\{T'_{\rr',\Cdot}( \bar Q',\bar\TT')\pp\Big\} \prod_{j>J'}\theta_j^{r_j}.
\end{align*}
Now that $T(Q,\TT)\pp = T(\bar Q,\bar\TT)\bar\pp$ but $(\bar Q,\bar\TT,\bar\pp) \nsim (Q,\TT,\pp)$, we obtain that $(Q,\TT,\pp)$ are not strictly (or generically) identifiable. The proof of the proposition is complete.

\section{Proof of Theorem 1}\label{sec-proof-dina1}

%\paragraph{Proof of Theorem \ref{thm-dina}}

%\begin{proof}[Proof of Theorem \ref{thm-dina}]
We first prove the sufficiency, and then show the necessity of the conditions. 
Under DINA, \eqref{eq1} can be equivalently written as that for any $\rr\in\{0,1\}^J$,
%\begin{equation}\label{eq-def}
%\mathbb P(\RR\mid Q,\cc,\cg,\pp) = \mathbb P(\RR\mid \bar Q,\bar \cc,\bar \cg,\bar \pp) 
%\end{equation}
\begin{equation}\label{eq-def-T}
T_{\rr,\Cdot}(Q, \cc,\cg) \pp = T_{\rr,\Cdot}(\bar Q,  \bar\cc,\bar\cg)\bar\pp.
\end{equation}
%\color{blue}
We first introduce some notations. In the following discussion, for an integer $M$, we denote $[M]=\{1,\ldots,M\}$. For an item set $S\subseteq[J]$, denote $\qq_S = \vee_{j\in S}\qq_j= (\max_{j\in S}q_{j,1},\allowbreak\max_{j\in S}q_{j,2},\ldots,\allowbreak\max_{j\in S}q_{j,K})$, then $\qq_S$ is also a $K$-dimensional binary vector, and we denote its $k$ element by $q_{S,k}$.
Recall  
\begin{equation*}
Q = \begin{pmatrix}
I_K\\
Q^\star
\end{pmatrix},%\quad\quad
%\bar Q = \begin{pmatrix}
%\bar Q_{1:K,\Cdot}\\
%\bar Q^\star
%\end{pmatrix}.
\end{equation*}
and we denote the submatrix of $\bar Q$ consisting of its first $K$ row vectors by $\bar Q_{1:K,\Cdot}$.
We next show in five steps that if \eqref{eq-def-T} holds, then  $\bar Q\sim Q$, and also $\cc=\bar\cc$, $\bar\cg=\cg$, $\bar\pp=\pp$. %And in the end, we show the necessity of the Conditions $A$, $B$ and $C$ for identifiability.

\noindent {\bf Step 1.} After some column rearrangement, $\bar Q_{1:K,\Cdot}$ is an upper-triangular matrix with all the diagonal elements being ones.

\noindent {\bf Step 2.} $\bar c_j=c_j$ for all $j\in\{K+1,\ldots,J\}$.

\noindent {\bf Step 3.} $\bar g_k=g_k$ for all $k\in\{1,\ldots,K\}$.

\noindent {\bf Step 4.}  $\bar Q_{1:K,\Cdot}\sim I_K$ %After some column rearrangement, $\bar Q_{1:K,\Cdot}$ is an $K\times K$ identity matrix.

\noindent {\bf Step 5.} $\bar Q\sim Q$, $\cc=\bar\cc$, $\bar\cg=\cg$, $\bar\pp=\pp$.
%$\bar c_k=c_k$ for all $j\in\{1,\ldots,K\}$, $\bar g_j=g_j$ for all $j\in\{K+1,\ldots,J\}$.

%%%
For any item set $S\subseteq\{1,\ldots,J\}$, denote $\cc_S=\sum_{j\in S}c_j\ee_j$, %$\cg_S = \sum_{j\in S}g_j\ee_j$,
 and denote $\cg_S$, $\bar\cc_S$, and $\bar\cg_S$ similarly. Consider the response pattern $\rr^\star= \sum_{j\in S}\ee_j$ and any $\ttt^\star = \sum_{j\in S}\theta^\star_j\ee_j$, then Equation \eqref{eq-def-T} together with Lemma \ref{lem} imply that
\begin{equation}\label{eq-tra}
T_{\rr^\star,\Cdot}( Q,\cc_S-\ttt^\star,\cg_S-\ttt^\star)\pp = T_{\rr^\star,\Cdot}(\bar Q,\bar\cc_S-\ttt^\star,\bar\cg_S-\ttt^\star)\bar\pp.
\end{equation} 
We will frequently use \eqref{eq-tra} in the following proof. And when the item set $S$ and response pattern $\rr^\star$ are clearly implied by the definition of $\ttt^\star$, we will omit the subscript $S$ in the above \eqref{eq-tra}.
We also frequently use the fact that when \eqref{eq-tra} holds, $c_j\neq\bar g_j$ and $g_j\neq \bar c_j$ for any item $j$. This is true because if $c_j=\bar g_j$, we would have
\begin{align*}
T_{\ee_j,\Cdot}(Q,\cc,\cg)\pp = 
&~c_j(\sum_{\aaa\succeq\qq_j} p_{\aaa}) + g_j(\sum_{\aaa\nsucceq\qq_j} p_{\aaa})< c_j = \bar g_j \\
 \leq & \bar c_j(\sum_{\aaa\succeq\qq_j} \bar p_{\aaa}) + \bar g_j(\sum_{\aaa\nsucceq\qq_j} \bar p_{\aaa})
 =T_{\ee_j,\Cdot}(\bar Q,\bar\cc,\bar\cg)\bar\pp,
\end{align*}
which contradicts \eqref{eq-def-T}. So $c_j\neq\bar g_j$ and similarly $g_j\neq \bar c_j$ for each $j$. 
As stated in the main text, we assume without loss of generality that there is no all-zero row vector in true $Q$-matrix. If, however, the $j$th row vector of $\bar Q$ equals $\mz$, then $\bar c_j$ would equal $\bar g_j$, and we denote this value by $\bar\theta_j$. Equation \eqref{eq-def-T} gives
$$
\bar\theta_j = c_j\Big(\sum_{\aaa:\,\aaa\succeq\qq_j} p_{\aaa}\Big) + g_j\Big(\sum_{\aaa:\,\aaa\nsucceq\qq_j} p_{\aaa}\Big),
$$
and hence $g_j<\bar\theta_j<c_j$ holds for this $j$.
\color{black}

%and the underlying true $Q$-matrix $Q$ takes the form of \eqref{eq-true-q}.

%\vspace{2mm}
%\color{blue}
\medskip
\noindent
\textbf{Step 1.} In this step we prove that $\bar Q_{1:K,\Cdot}$ must take the following form after some column rearrangement,
%form an identity matrix, up to some row permutation. 
\begin{equation}\label{eq-qbar-K}
\bar Q_{1:K,\Cdot} \sim
\begin{pmatrix}
    1 & * & \dots  & * \\
    0 & 1 & \dots  & * \\
    \vdots & \vdots & \ddots & \vdots \\
    0 & 0 & \dots  & 1
\end{pmatrix}.
\end{equation}
Namely, after properly rearranging the columns of $\bar Q_{1:K,\Cdot}$, we have $\bar Q_{k,k}=1$ and $\bar Q_{k,h}=0$ for any $k>h$.

We first introduce the following useful lemma.
%\begin{lemma}\label{lem-H1}
%Suppose the true $Q$ satisfies Condition $A$ that $Q_{1:K}=I_K$.
%If there exists an item set  $S\subseteq\{K+1,\ldots,J\}$ such that 
%$$\max_{m\in S}q_{m,h}=0,\quad \max_{m\in S} q_{m,j}=1$$ for some attributes $h,j\in[K]$ and $h\neq j$, then $\bar\qq_j\nsucceq\bar\qq_h$.
%\end{lemma}
%\color{purple}
\begin{lemma}\label{lem-H1}
Suppose the true $Q$ satisfies Condition $A$ that $Q_{1:K}=I_K$.
If there exists an item set  $S\subseteq\{K+1,\ldots,J\}$ such that 
$$\max_{m\in S}q_{m,h}=0,\quad \max_{m\in S} q_{m,j}=1~~\forall j\in\mathcal J$$ 
for some attributes $h\in[K]$ and a set of attributes $\mathcal J\subseteq[K]\setminus\{h\}$, then 
$$
\vee_{j\in\mathcal J}\,\bar\qq_j\nsucceq\bar\qq_h.
$$
\end{lemma}
%\color{blue}
{\bf Proof of Lemma \ref{lem-H1}.}
We prove by  contradiction. Assume there exist attribute $h\in[K]$ and a set of attributes $\mathcal J\subseteq[K]\setminus\{h\}$, such that $\vee_{j\in\mathcal J}\,\bar\qq_j\nsucceq\bar\qq_h$; and that there exists  $S\subseteq\{K+1,\ldots,J\}$ such that $\max_{m\in S}q_{m,h}=0$ and $\max_{m\in S} q_{m,j}=1$. %then if an attribute pattern $\aaa$ satisfies $\aaa\succeq \bar\qq_{j}$, it must also satisfy $\aaa\succeq\bar\qq_{h}$. 
Define
$$
\ttt^\star = \bar c_h\ee_{h} + \sum_{j\in\mt J}\bar g_j\ee_{j} 
+ \sum_{m=K+1}^J g_{m}\ee_m,\quad \rr^\star = \ee_{h} + \sum_{j\in\mt J}\ee_{j} + \sum_{m=K+1}^J \ee_m,
$$
and we claim that $T_{\rr^\star,\Cdot}(\bar Q,\bar\cc-\ttt^\star,\bar\cg-\ttt^\star)$ is an all-zero vector. This is because for any $\aaa\in\{0,1\}^K$, the corresponding element in $T_{\rr^\star,\aaa}(\bar Q,\bar\cc-\ttt^\star,\bar\cg-\ttt^\star)$ contains a factor $F_{\aaa} = (\bar \theta_{h,\aaa} - \bar c_{h})\prod_{j\in\mt J}(\bar \theta_{j,\aaa} - \bar g_{j})$. While this factor $F_{\aaa}\neq 0$ only if $\bar\theta_{h,\aaa}=\bar g_h$ and $\bar \theta_{j,\aaa} = \bar c_{j}$ for all $j\in\mt J$, which happens if and only if $\aaa\nsucceq\bar \qq_{h}$ and $\aaa\succeq\bar \qq_{j}$  for all $j\in\mt J$, which is impossible because $\vee_{j\in\mt J}\bar \qq_{j}\succeq\bar  \qq_{h}$ by our assumption. So the claim $T_{\rr^\star,\Cdot}(\bar Q,\bar\cc-\ttt^\star,\bar\cg-\ttt^\star)=\mz$ is proved, and further $T_{\rr^\star,\Cdot}(\bar Q,\bar\cc-\ttt^\star,\bar\cg-\ttt^\star)\bar\pp=0$. %Equation \eqref{eq-def-T} and Lemma \ref{lem} gives 
Equality \eqref{eq-tra} becomes
$
T_{\rr^\star,\Cdot}( Q,\cc-\ttt^\star,\cg-\ttt^\star)\bar\pp = T_{\rr^\star,\Cdot}(\bar Q,\bar\cc-\ttt^\star,\bar\cg-\ttt^\star)\bar\pp=0,$
which leads to 
\begin{align*}%\label{eq-t11}
0=T_{\rr^\star,\Cdot}( Q,\cc-\ttt^\star,\cg-\ttt^\star)\pp
=
%\Big(\sum_{\aaa:\alpha_{j}=0,\atop \alpha_{h}=1}p_{\aaa}\Big) 
p_{\mo}(c_{h} - \bar c_{h})\prod_{j\in\mt J} (c_{j} - \bar g_{j})\prod_{m>K}(c_m-g_m),
\end{align*}
which is because for any $\aaa\neq\mo$, we must have $\aaa\nsucceq\qq_m$ for some $m>K$ under Condition $C$, and hence the element $T_{\rr^\star,\aaa}( Q,\cc-\ttt^\star,\cg-\ttt^\star)$ contains a factor $(g_m - g_m)=0$. 
Since $c_m-g_m>0$ for $m>K$ and $c_{j} - \bar g_{j}\neq 0$, we obtain $c_{h} = \bar c_{h}$. 

We remark here that $c_{h} = \bar c_{h}$ also implies $\bar\qq_h\neq\mz$, because otherwise we would have $\bar\theta_h=\bar c_h=c_h$, which contradicts the $g_h<\bar\theta_h<c_h$ proved before the current Step 1. 
This indicates the $\bar Q_{1:K,\Cdot}$ can not contain any all-zero row vector, because otherwise $\bar\qq_j\succeq\bar\qq_h$ for the all-zero row vector $\bar\qq_h$, which we showed is impossible.

%We now define
%$$
%\ttt^\star = \bar g_h\ee_{h} + \bar c_j\ee_{j} 
%+ \sum_{m=K+1}^J g_{m}\ee_m,\quad \rr^\star = \ee_{h} + \ee_{j} + \sum_{m=K+1}^J \ee_m,
%$$
%then similar to the argument in the previous paragraph, we also have $T_{\rr^\star,\Cdot}( Q,\cc-\ttt^\star,\cg-\ttt^\star)=\mz$ and hence the right hand side of \eqref{eq-tra} is zero. And further,
%\begin{align*}
%0=T_{\rr^\star,\Cdot}( Q,\cc-\ttt^\star,\cg-\ttt^\star)\pp=
%p_{\mo}(c_{h} - \bar g_{h}) (c_{j} - \bar c_{j})\prod_{m>K}(c_m-g_m),
%\end{align*}
%since $c_h-\bar g_h\neq 0$, the above display implies $c_j=\bar c_j$.

%Under Condition $B$ that $Q^\star$ has distinct columns, we have that at least one of the two column vectors $Q^\star_{\Cdot,j}$ and $Q^\star_{\Cdot,h}$ does not equal the all-one vector, i.e., either $Q^\star_{\Cdot,j}\neq\mathbf1$ or $Q^\star_{\Cdot,h}\neq\mathbf1$. We next consider these two cases separately.

%\noindent
%\textbf{Case 1:}  
%Under the assumption $Q^\star_{\Cdot,j}\neq\mathbf1$, we assume without loss of generality $q_{l,j}=1$ for some $l>K$. 
%Define 
%$$
%\ttt^\star = \bar c_h\ee_{h} + \bar g_j\ee_{j} 
%%+ \sum_{m\in S}\theta_m\ee_m,
%$$
%and note that $\bar c_h = c_h$.
%Similar to the previous argument, the right hand side (RHS) of Equation \eqref{eq-tra} is zero, so we have its left hand side (LHS) being
%$$
%(g_{h} - \bar c_{h})\Big[(c_j-\bar g_j)\Big(\sum_{\aaa:\,\alpha_j=1, \alpha_h=0}p_{\aaa} \Big)
%+ (g_j-\bar g_j)\Big(\sum_{\aaa:\,\alpha_j=0, \alpha_h=0}p_{\aaa} \Big)\Big]=0.
%$$
Consider the item set $S$ in the lemma that satisfies $S\subseteq\{K+1,\ldots,J\}$ such that $\max_{m\in S}q_{m,h}=0$ and $\max_{m\in S}q_{m,j}=1$ for all $j\in\mt J$. Define
$$
\ttt^\star = \bar c_h\ee_{h} + \sum_{j\in\mt J}\bar g_j\ee_{j} +\sum_{m\in S} g_m\ee_{m}.
$$
Note that $c_h = \bar c_h$.
The RHS of \eqref{eq-tra} is zero, and so is the LHS of it. The row vector $T_{\rr^\star,\Cdot}(Q,\cc-\ttt^\star,\cg-\ttt^\star)$ has the following property
\begin{align*}
&T_{\rr^\star,\aaa}(Q,\cc-\ttt^\star,\cg-\ttt^\star)\\
=
&\begin{cases}
(g_{h} - \bar c_{h})\prod_{j\in\mt J}(c_j-\bar g_j)\prod_{m\in S}(c_m-g_m), 
& \aaa\nsucceq\qq_h,\,\aaa\succeq \qq_{\mt J},\, \aaa\succeq\qq_S;\\
0, & \text{otherwise}.
\end{cases}
\end{align*}
An important observation is that $\{\aaa\in\{0,1\}^K:\,\aaa\nsucceq\qq_h,\,\aaa\succeq \qq_{\mt J}, \,\aaa\succeq\qq_S\}=\mt A \neq\varnothing$. This is because %$\max_{m\in S}q_{m,h}=0$ and $\max_{m\in S}q_{m,j}=1$ 
$q_{S,h}=0$ and $q_{S,j}=1$ for all $j\in \mt J$ hold, and we can just choose $\aaa$ for which $\alpha_h=0$ and $\alpha_k=1$ for all $q_{S,k}=1$, then such $\aaa$ belongs to the set $\mathcal A$.
Therefore we have
\begin{align*}
&~T_{\rr^\star,\Cdot}(Q,\cc-\ttt^\star,\cg-\ttt^\star)\pp\\
=
&~(g_{h} - \bar c_{h})\prod_{j\in\mt J}(c_j-\bar g_j)\prod_{m\in S}(c_m-g_m)\Big(\sum_{\aaa\in\mt A}p_{\aaa} \Big)=0,
\end{align*}
which leads to a contradiction since $g_h- \bar c_h\neq 0$, $c_j-\bar g_j\neq 0$, $c_m-g_m\neq0$ and $\sum_{\aaa\in\mt A}p_{\aaa}>0$, i.e., every factor in the above product is nonzero. This completes the proof of Lemma \ref{lem-H1}.
\bigskip
%\end{proof}
%\noindent
%\textbf{Case 2:} $Q^\star_{\Cdot,h}\neq\mo$. Assume without loss of generality $q_{l,h}=0$ for some $h>K$.  Since $q_{l,h}=0$, if an attribute pattern $\aaa\nsucceq\qq_l$, then it must satisfy $\aaa\nsucceq\qq_h=(\mz,\underbrace{1}_{\text{column }h},\mz)$. Note that $c_j = \bar c_j$. Define 
%$\ttt^\star = \bar g_h\ee_{h} + \bar c_j\ee_{j} 
%+ c_l\ee_l$, and \eqref{eq-tra} gives
%$$
%(g_{h} - \bar g_{h}) (g_{j} - \bar c_{j})(g_l - c_l)\Big(\sum_{\aaa:\,\aaa\nsucceq\qq_l,\atop \alpha_h=0,\alpha_j=0} p_{\aaa}\Big)=0,
%$$
%which further implies $g_h = \bar g_h$. Now we have obtained $c_j=\bar c_j$ and $g_h = \bar g_h$, we further define $\ttt^\star = g_h\ee_h+c_j\ee_j$, then the RHS of \eqref{eq-tra}=0 and further the LHS equals
%\begin{align*}
%0=(c_h-g_h) (g_{j} - c_{j})\Big(\sum_{\aaa:\,\alpha_h=1,\,\alpha_j=0} p_{\aaa}\Big),
%\end{align*}
%which also clearly is a contradiction. This means under Case 2 that $Q^\star_{\Cdot,h}\neq\mo$, the assumption $\bar\qq_j\succeq\bar\qq_h$ also can not happen for any $j\in[K]\setminus\{h\}$.
%To summarize, by far we have shown that when Condition $B$ that \textit{$Q^\star$ has distinct columns} is satisfied, $\bar\qq_j\succeq\bar\qq_h$ can not happen for $j\neq h$, $j,h\in[K]$. This concludes the proof of Lemma \ref{lem-H1}.

%Lemma \ref{lem-H1} naturally leads to the following statement.
%\begin{state}\label{state-1}
%Suppose the true $Q$ satisfies Conditions A and B. Lemma \ref{lem-H1} implies $\bar Q_{1:K,\Cdot} \sim I_K$.
%\end{state}
We now proceed with the proof of Step 1 using an induction argument. %that $\bar Q_{1:K,\Cdot}\sim I_K$.
We first introduce the definition of \textit{lexicographic order} between two binary vectors of the same length. Specifically, for two binary vectors $\bo a=(a_1,\ldots,a_L)^\top$ and $\bo b=(b_1,\ldots,b_L)^\top$ both of length $L$, we say $\bo a$ is of smaller lexicographic order than $\bo b$ and denote $\bo a\prec_{\text{lex}}\bo b$, if either $a_1<b_1$, or there exists a integer $l\in\{2,\ldots,L\}$ such that $a_l<b_l$ and $a_m=b_m$ for all $m=1,\ldots,l-1$. It is not hard to see when Condition $B$ that $Q^\star$ contains $K$ distinct column vectors is satisfied, the $K$ columns of $Q^\star$ can be arranged in an increasing lexicographic order. Namely, under Condition $B$, there exists a permutation map $\sigma(\cdot):[K]\to[K]$ such that 
\begin{equation}\label{eq-qlex}
Q^\star_{,\sigma(1)}
\prec_{\text{lex}}
Q^\star_{,\sigma(2)}
\prec_{\text{lex}}
\cdots
\prec_{\text{lex}}
Q^\star_{,\sigma(K)}.
\end{equation}
Without loss of generality, next we consider the case where $\sigma(\cdot)$ is the identity map, i.e., $\sigma(k)=k$ for all $k\in[K]$.

%We first consider attribute $\sigma(1)$. Since $Q^\star_{,\sigma(1)}$ has the smallest lexicographic order among the columns of $Q^\star$, we have the conclusion that there must exist an item set $S\subseteq\{K+1,\ldots,J\}$ such that 
%$$
%q_{S,\sigma(1)}=0,\quad q_{S,\sigma(k)}=1~~\forall k=2,\ldots,K.
%$$
%Therefore Lemma \ref{lem-H1} implies $\bar\qq_{\sigma(k)}\nsucceq\qq_{\sigma(1)}$ for all $k=2,\ldots,K$, which  means 
%\begin{align*}
%&~(\max_{k\in[K]\setminus\{1\}}\,\bar  q_{\sigma(k),1},
%\max_{k\in[K]\setminus\{1\}}\,\bar  q_{\sigma(k),2},\ldots,
%\max_{k\in[K]\setminus\{1\}}\,\bar  q_{\sigma(k),K})\\
%\nsucceq
%&~(\bar q_{\sigma(1),1},~\ldots, ~\bar q_{\sigma(1),K}).
%\end{align*}
%This implies there must exist an attribute $m_1\in[K]$ such that %$\bar q_{\sigma(k),m_1} = 0$ for all $k\in\{2,\ldots,K\}$ and $\bar q_{\sigma(1),m_1} = 1$.
%\begin{equation}\label{eq-hh}
%%\begin{cases}
%\max_{k\in[K]\setminus\{1\}}\bar  q_{\sigma(k),m_1}=0,\quad
%\bar q_{\sigma(1),m_1}=1,
%%\end{cases}
%\end{equation} which exactly says the $m_1$-th column vector of $\bar Q_{1:K,\Cdot}$ must equal the basis vector $(\mz,\underbrace{1}_{\text{column }\sigma(1)},\mz)^\top=\ee_{\sigma(1)}$, i.e., we have $\bar Q_{1:K,m_1}=\ee_{\sigma(1)}$. 

We first consider attribute $1$. Since $Q^\star_{,1}$ has the smallest lexicographic order among the columns of $Q^\star$, we have the conclusion that there must exist an item set $S\subseteq\{K+1,\ldots,J\}$ such that 
$$
q_{S,1}=0,\quad q_{S,\ell}=1~~\forall \ell=2,\ldots,K.
$$
We apply Lemma \ref{lem-H1} to obtain {$\vee_{\ell\in\{2,\ldots,K\}}\bar\qq_\ell\nsucceq\bar\qq_1$}, which  means 
\begin{align*}
&~(\max_{m\in\{2,\ldots,K\}}\,\bar  q_{\ell,1},
\max_{m\in\{2,\ldots,K\}}\,\bar  q_{\ell,2},\ldots,
\max_{m\in\{2,\ldots,K\}}\,\bar  q_{\ell,K})\\
\nsucceq
&~(\bar q_{1,1},~\ldots, ~\bar q_{1,K}).
\end{align*}
This implies there must exist an attribute $m_1\in[K]$ such that %$\bar q_{\sigma(k),m_1} = 0$ for all $k\in\{2,\ldots,K\}$ and $\bar q_{\sigma(1),m_1} = 1$.
\begin{equation}\label{eq-hh}
%\begin{cases}
\max_{k\in[K]\setminus\{1\}}\bar  q_{k,m_1}=0,\quad
\bar q_{1,m_1}=1,
%\end{cases}
\end{equation} which exactly says the $m_1$-th column vector of $\bar Q_{1:K,\Cdot}$ must equal the basis vector $(\underbrace{1}_{\text{column }1},\mz)^\top=\ee_{1}$, i.e., we have $\bar Q_{1:K,m_1}=\ee_{1}$. 

Now we assume as the inductive hypothesis that for $h\in[K]$ and $h>1$, we have a distinct set of attributes $\{m_1,\ldots,m_{h-1}\}\subseteq[K]$ such that their corresponding column vectors in $\bar Q_{1:K,\Cdot}$ satisfy 
\begin{equation}\label{eq-induc-q1}
\forall i=1,\ldots,h-1,\quad \bar Q_{1:K,m_i}=(*,\ldots,*,\underbrace{1}_{\text{column }i}, 0,\ldots,0)^\top.
\end{equation}
 Now we focus on attribute $h$. By \eqref{eq-qlex}, the column vector $Q^\star_{,h}$ has the smallest lexicographic order among the $K-h-1$ columns in $\{Q^\star_{\Cdot,h},\,\allowbreak  Q^\star_{\Cdot,h+1},\,\allowbreak \ldots, Q^\star_{\Cdot,K}\}$, therefore similar to the argument in the previous paragraph, there must exist an item set $S\subseteq\{K+1,\ldots,J\}$ such that 
\begin{equation}\label{eq-sigma-hk}
q_{S,h}=0,\quad q_{S,\ell}=1~~\forall \ell=h+1,\ldots,K.
\end{equation}
Therefore Lemma \ref{lem-H1} implies {$\vee_{\ell\in\{h+1,\ldots,K\}}\bar\qq_\ell\nsucceq\bar\qq_1$},  and further leads to 
\begin{equation}\label{eq-mkh}
%\begin{cases}
\max_{\ell\in\{h+1,\ldots,K\}}\bar  q_{\ell,m_h}=0,\quad
\bar q_{h,m_h}=1.
%\end{cases}
\end{equation}
We point out that $m_h\not\in\{m_1,\ldots,m_{h-1}\}$, because by the induction hypothesis \eqref{eq-induc-q1} we have  $\bar q_{h,m_i}=0$ for $i=1,\ldots,h-1$.
%Without loss of generality, consider the case where $\sigma(\cdot)$ is the identity map and $k=k$ for all $k\in[K]$, 
So $\{m_1,\ldots,m_{h-1},m_h\}$ contains $h$ distinct attributes.
Furthermore, \eqref{eq-mkh} gives that 
$\bar Q_{\Cdot,m_h} = (*,\ldots,*,\allowbreak\underbrace{1}_{\text{column }h},\allowbreak 0,\ldots,0)^\top,$
which generalizes \eqref{eq-induc-q1} by extending $h-1$ there to $h$.
Therefore, we use the induction argument to obtain
$$\forall k\in\{1,\ldots,K-1\},\quad
\bar Q_{1:K,m_k} = (*,\ldots,*,\underbrace{1}_{\text{column }k}, 0,\ldots,0)^\top.$$
Furthermore, when considering the last attribute $K$, the $K$th item must have $\qq$-vector taking the form of $\bar \qq_K=(0,\ldots,0,\underbrace{*}_{\text{column }m_K},0,\ldots,0)$, where the ``$*$'' in $\bar\qq_K$ is the only element unspecified. Since previously we have shown in the proof of Lemma \ref{lem-H1} that $\bar\qq_j=0$ can not happen for any item $j$, there must be $\bar \qq_K=(0,\ldots,0,\underbrace{1}_{\text{column }m_K},0,\ldots,0)$.
Now we have  essentially obtained
\begin{equation}\label{eq-qmk}
\bar Q_{1:K,\,(m_1,\ldots,m_K)} =
\begin{pmatrix}
    1 & * & \dots  & * \\
    0 & 1 & \dots  & * \\
    \vdots & \vdots & \ddots & \vdots \\
    0 & 0 & \dots  & 1
\end{pmatrix},
\end{equation}
and the conclusion of Step 1 in \eqref{eq-qbar-K} is proved.

\color{black}

%\vspace{2mm}
\medskip
\noindent
\textbf{Step 2.} In this step we prove $c_j = \bar c_j$ for $j=K+1,\ldots,J$. 
%With the result of Step 1, suppose $\bar Q$ takes the following form
%\[
%\bar Q = \left(
%\begin{array}{c}
%I_K\\
%\bar Q^\star
%\end{array}
%\right).
%\]
%We first show for any item $h$, \eqref{eq-def-T} implies $c_h> \bar g_h$ and $g_h< \bar c_h$. Since $c_h>g_h$ and $p_{\aaa}>0$ for all $\aaa$,
%\begin{equation}\label{eq-cg}
%c_h = \sum_{\aaa\in\{0,1\}^K} c_h p_{\aaa}
%> \sum_{\aaa\in\{0,1\}^K} \theta_{h,\aaa}p_{\aaa} 
%= \sum_{\aaa\in\{0,1\}^K} \bar\theta_{h,\aaa}\bar p_{\aaa}
%>\sum_{\aaa\in\{0,1\}^K} \bar g_h \bar p_{\aaa} =\bar g_h,
%\end{equation}
%and similarly $g_h<\bar c_h$.
For an arbitrary item $j\in\{K+1,\ldots,J\}$, define a response vector $\rr^* = \sum_{h:\, h\neq j}\ee_j$ and
\[\ttt^* = \sum_{h=1}^K \bar g_h \ee_h + \sum_{h>K,\,h\neq j}g_h \ee_h.\]
{
We claim that $T_{\rr^*,\Cdot}(\bar Q,\bar\cc-\ttt^*,\bar\cg-\ttt^*)$ contains only one nonzero element corresponding to the all-one attribute pattern $\aaa=\mo$. The reasoning is as follows. Under the conclusion of Step 1, $\bar Q_{1:K,\Cdot}$ takes the form of \eqref{eq-qbar-K}, %up to some column permutation $\bar Q_{1:K,\Cdot}$ equals an upper triangular matrix with all the diagonal elements being one, 
which means each attribute is required by at least one item in $\{\bar\qq_1,\ldots,\bar\qq_K\}$. Then for any $\aaa\neq\mo$, there must exist some attribute $k\in[K]$ such that $\aaa\nsucceq\bar\qq_k$, which implies for this particular $\aaa$ the element $T_{\rr^*,\aaa}(\bar Q,\bar\cc-\ttt^*,\bar\cg-\ttt^*)$ contains a factor $(\bar g_h-\bar g_h)=0$. Therefore $T_{\rr^*,\aaa}(\bar Q,\bar\cc-\ttt^*,\bar\cg-\ttt^*)\neq 0$ only if $\aaa=\mo$. 
}% end of blue
Next consider $T_{\rr^*,\aaa}( Q,\cc-\ttt^*,\cg-\ttt^*)$. Under Condition $A$, in the true $Q$ each attribute is required by at least three items, so the row vector corresponding to response pattern $\rr^*$ in $T(Q,\cc-\ttt^*,\cg-\ttt^*)$ only contains one nonzero element, in column  $\aaa=\mo_K^\top$, representing the attribute profile mastering all the $K$ attributes. This is because for any other attribute profile $\aaa'$ that lacks at least one attribute $k$, there must be some item $h>K$, $h\neq j$ requiring attribute $k$ so that $\aaa'\nsucceq\qq_h$; and this results in $\theta_{\ee_h,\aaa'} = g_h$ and $T_{\rr^*,\aaa'}(Q,\cc-\ttt^*,\cg-\ttt^*)=0$. In summary,
\[\begin{aligned}
T_{\rr^*,\aaa}(Q,\cc-\ttt^*,\cg-\ttt^*) 
= \prod_{h=1}^K (\theta_{h,\aaa}-\bar g_h)\prod_{h>K:\atop h\neq j} (\theta_{h,\aaa}- g_h) \neq 0\quad\text{iff}\quad
\aaa &= \mo; 
\end{aligned}\]
\[\begin{aligned}
T_{\rr^*,\aaa}(\bar Q,\bar\cc-\ttt^*,\bar\cg-\ttt^*) 
= \prod_{h=1}^K (\bar\theta_{h,\aaa}-\bar g_h)\prod_{h>K:\atop h\neq j} (\bar\theta_{h,\aaa}- g_h)
\neq 0\quad\text{iff}\quad
\aaa &= \mo. \\
\end{aligned}\]
%where the inequalities $c_h> \bar g_h$, $g_h<\bar c_h$ ensure $T_{\rr^*, \mo_K^\top }(Q,\cc-\ttt^*,\cg-\ttt^*)\neq 0$ and $T_{\rr^*, \mo_K^\top }(\bar Q,\bar\cc-\ttt^*,\bar\cg-\ttt^*)\neq 0$.
Now further consider item $j$. Since $\mo^\top_K\succeq\qq_j$ and $\mo^\top_K\succeq\bar\qq_j$, one must have $\theta_{j,\mo^\top_K} = c_j$ and $\bar\theta_{j,\mo^\top_K} =\bar c_j$. Since we assume $p_{\aaa}>0$ for each $\aaa$, we have $T_{\rr^*,\Cdot}(Q,\cc-\ttt^*, \cg -\ttt^*)\pp = T_{\rr^*,\mo^\top_K}(Q,\cc-\ttt^*, \cg -\ttt^*) p_{\mo^\top_K} \neq 0$. So \eqref{eq-trans} in Lemma \ref{lem} implies that
\[
c_j = 
\frac{T_{\rr^*+\ee_j,\Cdot}(Q,\cc-\ttt^*, \cg-\ttt^*)\pp}{T_{\rr^*,\Cdot}(Q,\cc-\ttt^*, \cg -\ttt^*)\pp} = 
\frac{T_{\rr^*+\ee_j,\Cdot}(\bar Q,\bar\cc-\ttt^*, \bar\cg-\ttt^*)\bar\pp}{T_{\rr^*,\Cdot}(\bar Q,\bar\cc-\ttt^*, \bar\cg -\ttt^*)\bar\pp}
= \bar c_j.
\]
In the above argument $j$ is arbitrary, so $c_j = \bar c_j$ for any $j=K+1,\ldots,J$.

\medskip
\noindent
\textbf{Step 3.} In this step we prove $g_k = \bar g_k$ for $k=1,\ldots,K$. 
%We first introduce a definition of lexicographic order among binary vectors of the same length. For $\boldsymbol a=(a_1,\ldots,a_L)$, $\boldsymbol b=(b_1,\ldots,b_L)\in\{0,1\}^L$, we say $\boldsymbol a$ is of smaller lexicographic order then $\boldsymbol b$ if $a_1<b_1$ or there exists some $m>1$ such that $a_m<b_m$ and $a_l=b_l$ for $l=1,\ldots,m-1$. And we write $\boldsymbol a\prec_{\text{lex}} \boldsymbol b$ in this case.
%Without loss of generality, suppose in the true $Q$, the column vectors of the submatrix $Q^\star$ are ordered lexicographically (alphabetically), namely we assume that %for any $1\leq k_1<k_2\leq K$, 
Recall that in Step 1 we showed that \eqref{eq-qlex} about the lexicographic order holds and assumed $\sigma(k)=k$ for $k\in[K]$ without loss of generality.
%\begin{equation}\label{eq-lex}
%\text{for any }1\leq k_1<k_2\leq K,\quad
%Q^\star_{\Cdot,k_1} \prec_{\text{lex}} Q^\star_{\Cdot,k_2}.\end{equation} 
We now prove $g_1 = \bar g_1$. Define %$\rr^* = \sum_{h=1}^K\ee_h$ and 
\begin{equation}\label{eq-ttt}
\ttt^* = \sum_{h=1}^K \bar g_h \ee_h + \sum_{h>K:\atop q_{h,1}=0} g_h\ee_h + \sum_{h>K:\atop q_{h,1}=1} c_h\ee_h,
\end{equation}
then 
\[\begin{aligned}
T_{\sum_{h}\ee_h,\aaa}(Q,\cc-\ttt^*,\cg-\ttt^*) &= 
\prod_{h=1}^K (\theta_{h,\aaa}-\bar g_h)\prod_{h>K:\atop q_{h,1}=0}(\theta_{h,\aaa}-g_h)\prod_{h>K:\atop q_{h,1}=1}(\theta_{h,\aaa}-c_h);\\
T_{\sum_{h}\ee_h,\aaa}(\bar Q,\bar \cc-\ttt^*,\bar \cg-\ttt^*) &= 
\prod_{h=1}^K (\bar \theta_{h,\aaa}-\bar g_h)\prod_{h>K:\atop q_{h,1}=0}(\bar \theta_{h,\aaa}-g_h)\prod_{h>K:\atop q_{h,1}=1}(\bar \theta_{h,\aaa}-c_h).\\
\end{aligned}\]
First, the row vector $T_{\sum_{h=1}^J\ee_h,\Cdot}(\bar Q,\bar \cc-\ttt^*,\bar \cg-\ttt^*)$ equals the zero vector. This is because  $\bar Q_{1:K,\Cdot}$ takes the form in \eqref{eq-qbar-K} by Step 1, and  any attribute profile $\aaa\neq\mo_K^\top$ would have $\bar\theta_{h,\aaa}=\bar g_h$ for some $h\in\{1,\ldots,K\}$, which makes the corresponding element in the above row vector zero. Furthermore, $T_{\sum_{h=1}^J\ee_h,\mo_K^\top}(\bar Q,\bar \cc-\ttt^*,\bar \cg-\ttt^*)$ is also zero, because $\bar\theta_{h,\aaa} = \bar c_h=c_h$ for those $h>K$ such that $q_{h,1}=1$.
Since $Q^\star_{\Cdot,1}$ has the smallest lexicographic order among the columns of $Q^\star$, for any $k\in\{2,\ldots,K\}$, there must exist some item $h\in\{K+1,\ldots,J\}$ that requires attribute $1$, as a result 
\[\vee_{h>K:\,q_{h,1}=0}\, \qq_h = (0,1,\ldots,1).\] 
This ensures $T_{\sum_{h=1}^J\ee_h,\aaa}(Q,\cc-\ttt^*,\cg-\ttt^*)$ would equal zero if $\aaa$ lacks any attribute other than the first one. So the nonzero elements in the row vector $T_{\sum_{h=1}^J\ee_h,\Cdot}(Q,\cc-\ttt^*,\cg-\ttt^*)$ can only correspond to columns $\aaa^1=(0,1,\ldots,1)$ or $\aaa^2= \mo^\top_K$. Further, we claim $T_{\sum_{h=1}^J\ee_h,\aaa^2}(Q,\cc-\ttt^*,\cg-\ttt^*)=0$, this is because $\theta_{h,\aaa} = c_h$ for those $h$ such that $q_{h,1} = 1$. So the row vector $T_{\sum_{h=1}^J\ee_h,\Cdot}(Q,\cc-\ttt^*,\cg-\ttt^*)$ only contains one potentially nonzero element in column $\aaa_1 = (0,1,\ldots,1)$ as follows
\begin{equation}\label{eq-nzero}
T_{\sum_{h=1}^J\ee_h,\aaa_1}(Q,\cc-\ttt^*,\cg-\ttt^*) 
= (g_1 - \bar g_1) \prod_{h=2}^K(c_h - \bar g_h) \prod_{h>K:\atop q_{h,1}=0}(c_h-g_h)\prod_{h>K:\atop q_{h,1}=1}(g_h-c_h).
\end{equation}
Using the fact $T_{\sum_{h=1}^J\ee_h,\Cdot}(\bar Q,\bar \cc-\ttt^*,\bar \cg-\ttt^*)=\mz_{2^K}$, the equality 
\[T_{\sum_{h=1}^J\ee_h,\aaa^1}(Q,\cc-\ttt^*,\cg-\ttt^*)\pp = T_{\sum_{h=1}^J\ee_h,\aaa^1}(\bar Q,\bar\cc-\ttt^*,\bar\cg-\ttt^*)\bar\pp=0\] implies the element in \eqref{eq-nzero} must also be zero.
As shown earlier, $c_h-\bar g_h\neq 0$ for any $h$, so $g_1 = \bar g_1$ must hold.

Next we use an induction argument to prove that for $k=2,\ldots,K$, $g_k = \bar g_k$. In particular, suppose for any $1\leq m\leq k-1$, we already have $g_m = \bar g_m$. Define  
\begin{equation}\label{eq-ggc}
\ttt^* = \sum_{h=1}^K \bar g_h \ee_h + \sum_{h>K:\,q_{h,k}=0} g_h\ee_h + \sum_{h>K:\,q_{h,k}=1} c_h\ee_h.
\end{equation}
For the similar reason as stated before, $T_{\sum_{h=1}^J\ee_h,\Cdot}(\bar Q,\bar \cc-\ttt^*,\bar \cg-\ttt^*)$ equals the zero vector. We claim that the row vector $T_{\sum_{h=1}^J\ee_h,\Cdot}(Q,\cc-\ttt^*,\cg-\ttt^*)$ only contains one potentially nonzero element in column $\aaa' := (1,\ldots,1,\underbrace{0}_\text{column $k$},\allowbreak 1,\ldots,1)$. The reason is as follows. On the one hand, for any attribute profile $\aaa$ that lacks some attribute $l\in\{k+1,\ldots,K\}$, due to the assumption in \eqref{eq-qlex} that $Q_{\Cdot,k}^* \prec_{\text{lex}}Q_{\Cdot,l}^*$, there must exist some item $h>K$ such that $q_{h,k}=0$, $q_{h,l}=1$. So for this particular $\aaa$ we have $\aaa\nsucceq\qq_h$, $\theta_{h,\aaa}=g_h$, which makes $T_{\sum_{h=1}^J\ee_h,\aaa}(Q,\cc-\ttt^*,\cg-\ttt^*)=0$. On the other hand, for any attribute profile $\aaa'$ that lacks some attribute $m\in\{1,\ldots,k-1\}$, one has $\aaa'\nsucceq\qq_m=\ee_m$ and $\theta_{m,\aaa'}=g_m=\bar g_m$, where the last equality $g_m=\bar g_m$ comes from the induction assumption. This results in $T_{\sum_{h=1}^J\ee_h,\aaa'}(Q,\cc-\ttt^*,\cg-\ttt^*)=0$ for all such $\aaa'$. In conclusion, the nonzero elements in this transformed row vector can only be in columns $\aaa'$ or $\aaa_2 = \mo^\top_K$. For similar reason as in proving $g_1 = \bar g_1$, $T_{\sum_{h=1}^J\ee_h,\aaa_2}(Q,\cc-\ttt^*,\cg-\ttt^*)=0$. So the transformed row vector only contains one potentially nonzero entry corresponding to $\aaa'$:
\begin{align*}%\label{eq-nzero-2}
&~T_{\sum_{h}\ee_h,\aaa'}(Q,\cc-\ttt^*,\cg-\ttt^*) \\
= &~(g_k - \bar g_k) \prod_{1\leq h\leq K:\atop h\neq k}(c_h - \bar g_h) \prod_{h>K:\atop q_{h,k}=0}(c_h-g_h)\prod_{h>K:\atop q_{h,k}=1}(g_h-c_h).
\end{align*}
The same argument after \eqref{eq-nzero} gives $g_k = \bar g_k$. In conclusion, the induction method yields $g_k = \bar g_k$ for $k=1,\ldots,K$.

%\color{blue}
\medskip
\noindent
\textbf{Step 4.}  In this step we show that $\bar Q_{1:K,\Cdot}\sim I_K$. Recall that in Step 1 we already obtained  \eqref{eq-qmk}, and now we aim to show that the $\bar Q_{1:K,\,(m_1,\ldots,m_K)}$ in \eqref{eq-qmk} can be further written as
\begin{equation*}
\bar Q_{1:K,\,(m_1,\ldots,m_K)} =
\begin{pmatrix}
    1 & 0 & \dots  & 0 \\
    0 & 1 & \dots  & 0 \\
    \vdots & \vdots & \ddots & \vdots \\
    0 & 0 & \dots  & 1
\end{pmatrix}.
\end{equation*}
We now claim that $\bar\qq_j\nsucceq\bar\qq_h$ for any $1\leq j<h\leq K$. If this claim is true, then $\bar Q_{1:K,\,(m_1,\ldots,m_K)}=I_K$ must hold and the conclusion $\bar Q_{1:K,\Cdot}\sim I_K$ is reached.
We next prove that claim by contradiction. If there exist some $1\leq j<h\leq K$ such that $\bar\qq_j\succeq\bar\qq_h$, then define 
$$
\ttt^\star = \bar c_h\ee_h + \bar g_j\ee_j + \sum_{m=K+1}^J g_m\ee_m,
$$
we have 
\begin{align*}
0 = &~ T_{\rr^\star,\Cdot}(\bar Q,\bar \cc-\ttt^\star,\bar \cg-\ttt^\star)\bar\pp\\
= &~ T_{\rr^\star,\Cdot}(Q,\cc-\ttt^\star,\cg-\ttt^\star)\pp\\ =&~ p_{\mo}(c_h-\bar c_h)(c_j-\bar g_j)\prod_{m=K+1}^J (c_m-g_m),
\end{align*}
which implies $c_h = \bar c_h$. Note that we have obtained $g_j = \bar g_j$ in Step 3, and we next define $\ttt^\star = \bar c_h\ee_h + \bar g_j\ee_j$. The equality $T_{\rr^\star,\Cdot}(\bar Q,\bar \cc-\ttt^\star,\bar \cg-\ttt^\star)\bar\pp=0$ still holds and \eqref{eq-tra} gives
\begin{align}\label{eq-qcont}
0= &~T_{\rr^\star,\Cdot}(Q,\cc-\ttt^\star,\cg-\ttt^\star)\pp\\ \notag
= &~(g_h-\bar c_h)(c_j-\bar g_j)\Big(\sum_{\aaa:\,\aaa\nsucceq\qq_h,\aaa\succeq\qq_j}p_{\aaa}\Big)\\ \notag
=&~(g_h- c_h)(c_j- g_j)\Big(\sum_{\aaa:\,\aaa\nsucceq\qq_h,\aaa\succeq\qq_j}p_{\aaa}\Big).
\end{align}
Since $Q_{1:K,\Cdot}=I_K$, we have that $\qq_j$ and $\qq_h$ in the true $Q$ are distinct basis vectors, therefore $\Big(\sum_{\aaa:\,\aaa\nsucceq\qq_h,\aaa\succeq\qq_j}p_{\aaa}\Big)>0$. Therefore \eqref{eq-qcont} leads to a contradiction, and we have proved the claim that $\bar\qq_j\nsucceq\bar\qq_h$ for any $1\leq j<h\leq K$. As stated earlier, this claim naturally leads to the conclusion of Step 3 that $\bar Q_{1:K,\Cdot}\sim I_K$.

\color{black}

\medskip
\noindent
\textbf{Step 5.} In this step we prove that after reordering the columns in $\bar Q$ such that $\bar Q_{1:K} = I_K$, we must have $\qq_j = \bar{\qq}_j$ for $j = K+1,\ldots,J$.
%\color{blue}
In the following two parts, we first prove $\bar \qq_j\succeq\qq_j$ for all $j\in\{K+1,\ldots,J\}$ in part (a);  and then prove $\bar \qq_j=\qq_j$ for all $j\in\{K+1,\ldots,J\}$ in part (b).

\begin{enumerate}
\item[(a)] We next show $\bar \qq_j\succeq\qq_j$ for all $j\in\{K+1,\ldots,J\}$. We use proof by contradiction, and assume $\bar \qq_j\nsucceq\qq_j$ for some $j\in\{K+1,\ldots,J\}$. Then $\{\aaa:\,\aaa\succeq\bar\qq_j,\,\aaa\nsucceq\qq_j\}=\mathcal A\neq\varnothing$ and $\sum_{\aaa\in\mathcal A}p_{\aaa}\neq 0$.
Define
\begin{equation}\label{eq-def-exclude}
\ttt^* = \sum_{k\in[K]:\, \bar q_{j,k}=1}g_k\ee_k + c_j\ee_j,\end{equation}
then $T_{\rr^*,\Cdot}(\bar Q,\bar \cc-\ttt^*,\bar \cg-\ttt^*) = \mz$ and $T_{\rr^*,\Cdot}(\bar Q,\bar \cc-\ttt^*,\bar \cg-\ttt^*)\bar \pp=0$.  
However, for any $\aaa\in\mathcal A$, one has  $\theta_{j,\aaa}= g_j$ and $\theta_{k,\aaa}= c_k$ for any $k$ s.t. $\bar q_{j,k}=1$, so for any $\aaa\in\mathcal A$ we have
\begin{align*}
T_{\rr^*,\aaa}( Q,\cc-\ttt^*,\cg-\ttt^*) 
= & \prod_{1\leq k\leq K:\atop q_{j,k}=1} (\theta_{k,\aaa}-g_k)(\theta_{j,\aaa}-c_j) \\
= & \prod_{1\leq k\leq K:\atop q_{j,k}=1} ( c_k-g_k)( g_j-c_j) 
\neq 0,
\end{align*}
and hence
\[\begin{aligned}
T_{\rr^*,\Cdot}( Q,\cc-\ttt^*,\cg-\ttt^*)\pp
&=
\prod_{1\leq k\leq K:\atop q_{j,k}=1} ( c_k-g_k)( g_j-c_j) \sum_{\aaa\in\mathcal A}p_{\aaa} \\ %\Cdot \sum_{\aaa}I(\aaa\succeq \qq_j,~\aaa\nsucceq \bar\qq_j) \\
& \neq 0 = T_{\rr^*,\Cdot}(\bar Q,\bar \cc-\ttt^*,\bar \cg-\ttt^*)\bar \pp,
\end{aligned}\]
which contradicts \eqref{eq-def-T}. 

\item[(b)] Based on (a), we next show $\bar \qq_j=\qq_j$ for all $j\in\{K+1,\ldots,J\}$ using proof by contradiction. Since part (a) gives $\bar \qq_j\succeq\qq_j$, if $\bar \qq_j\neq\qq_j$, then there must exist some attribute $k\in[K]$ such that $\bar q_{j,k}=1$ and $q_{j,k}=0$. This implies $\bar\qq_j\succeq\bar\qq_k$. Define 
$$
\ttt^\star = \bar c_k\ee_k + \bar g_j\ee_j + \sum_{m>K:\,m\neq j}g_m\ee_m,
$$
then $T_{\rr^\star,\Cdot}(\bar Q,\bar\cc-\ttt^\star,\bar\cg-\ttt^\star)\bar\pp=0$. Since Condition $C$ holds, each attribute is required by at least one item in the set $\{m>K:\,m\neq j\}$, which implies $T_{\rr^\star,\aaa}( Q,\cc-\ttt^\star,\cg-\ttt^\star)\neq 0$ only if $\aaa=\mo$. Therefore \eqref{eq-tra} gives that
\begin{align*}
0=&~T_{\rr^\star,\Cdot}( Q,\cc-\ttt^\star,\cg-\ttt^\star)\pp\\
=&~(c_k-\bar c_k)(c_j-\bar g_j)\prod_{m>K:\,m\neq j}(c_m-g_m)p_{\mo},
\end{align*}
so $c_k=\bar c_k$. Now we further define 
$$
\ttt^\star = \bar c_k\ee_k + \bar g_j\ee_j + \sum_{h\in[K]\setminus\{k\}}g_m\ee_m,
$$
then $T_{\rr^\star,\Cdot}(\bar Q,\bar\cc-\ttt^\star,\bar\cg-\ttt^\star)\bar\pp=0.$ However,  $\qq_j\nsucceq\qq_k$ under the true $Q$, and \eqref{eq-tra} gives 
$$
T_{\rr^\star,\Cdot}( Q,\cc-\ttt^\star,\cg-\ttt^\star)\pp=(g_k-\bar c_k)\prod_{h\in[K]\setminus\{k\}}(c_h-g_h)(c_j-\bar g_j)p_{\aaa-\ee_k},
$$
where $\aaa-\ee_k=(\mo,\underbrace{0}_{\text{column }k},\mo)$, so the above display is nonzero. This contradicts \eqref{eq-tra}, and this means $\bar \qq_j\neq\qq_j$ can not happen. So we have $\bar\qq_j=\qq_j$ for $j\in\{K+1,\ldots,J\}$.
\end{enumerate}
\color{black}

Now we have proved $Q\sim\bar Q$.
Now that $Q\sim\bar Q$, Theorem 1 in \cite{dina} gives that Conditions A and B ensure the identifiability of the model parameters $(\cs:=\mo-\cc,\cg,\pp)$. This concludes the proof of the sufficiency of the conditions.

\medskip
In the end we show the necessity of the conditions.
By Theorem 1 in \cite{dina}, Conditions A and B are necessary for identifiability of the model parameters $(\cs,\cg,\pp)$ given a known $Q$, so they are also necessary for identifiability of $(Q, \cs,\cg,\pp)$.
%\end{proof}

\section{Proof of Theorem 2}
%\begin{proof}[Proof of Theorem \ref{thm-dina-gen}]
%\ \
\noindent
\textbf{Proof of the necessity of each attribute required by $\geq 2$ items.}
Suppose $Q$ takes the form of
$$
Q=\begin{pmatrix}
1 & \mathbf0^{\top} \\
\mathbf0 & Q^\star
\end{pmatrix},
$$
then for any valid $(\cc,\cg,\pp)$ associated with $Q$, we next construct $(\bar\cc,\bar\cg,\bar\pp)\neq(\cc,\cg,\pp)$ such that $T(Q,\cc,\cg)\pp=T(Q,\bar\cc,\bar\cg)\bar\pp$ holds. In particular, we arbitrarily choose $\bar c_1$ that is not equal to $c_1=1-s_1$ and set
$$
\bar p_{\aaa}=\begin{cases}
(c_1/\bar c_1) p_{\aaa},& \text{if}~\alpha_1=1,\\
p_{\aaa}+ (1-c_1/\bar c_1)p_{\aaa+\ee_1},& \text{if}~ \alpha_1=0.
\end{cases}
$$
Then set $\bar g_1=g_1$, and $\bar c_j=c_j$, $\bar g_j=g_j$ for $j=2,\ldots J$. Then it is not hard to check that $T(Q,\cc,\cg)\pp=T(Q,\bar\cc,\bar\cg)\bar\pp$. Since  $(\cc,\cg,\pp)$ are arbitrary, we have shown the non-identifiability set spans the entire parameter space and $(Q, \cc,\cg,\pp)$ are not generically identifiable.
Therefore, this proves that $(Q, \cc,\cg,\pp)$  are not generically identifiable if some attribute is required by only one item.

In the following we prove part (a), (b), and (c) when some attribute is required by only two items. 

\noindent
\textbf{Proof of Part (a).}
Under the assumption of part (a), $Q$ takes the form
$$
Q = \begin{pmatrix}
1 & \mathbf 0^\top \\
1 & \mathbf 1^\top \\
\mathbf 0 & Q^{\star}
\end{pmatrix}.
$$
Given arbitrary DINA model parameters $(\cc,\cg,\pp)$ under this $Q$,
we next construct another different set of DINA parameters $(\bar\cc,\bar\cg,\bar\pp)\neq (\cc,\cg,\pp)$ also associated with this $Q$, such that 
\begin{equation}\label{eq-tptp}
T(Q,\cc,\cg)\pp = T(Q,\bar\cc,\bar\cg)\bar\pp.
\end{equation}

In particular, we set $\bar c_j=c_j$ and $\bar g_j=g_j$ for all $j=3,\ldots,J$. Under this construction, \eqref{eq-tptp} simplifies to the following two sets of equations%\eqref{eq-anot1} and \eqref{eq-a1}
\begin{align}\label{eq-anot1}
\forall \aaa'\in\{0,1\}^{K-1},~\aaa'\neq\mathbf 1,\quad
&\begin{cases}
p_{(0,\aaa')}+p_{(1,\aaa')} =\bar p_{(0,\aaa')}+\bar p_{(1,\aaa')}, \\
g_1 p_{(0,\aaa')}+c_1 p_{(1,\aaa')} =\bar g_1\bar p_{(0,\aaa')}+\bar c_1 \bar p_{(1,\aaa')}, \\
g_2[ p_{(0,\aaa')}+ p_{(1,\aaa')} ] =\bar g_2[\bar p_{(0,\aaa')}+\bar p_{(1,\aaa')} ], \\
g_2[g_1 p_{(0,\aaa')}+c_2 p_{(1,\aaa')}] =\bar g_2[\bar g_1\bar p_{(0,\aaa')}+\bar c_1\bar p_{(1,\aaa')}];
\end{cases} 
\end{align}
and for $\aaa'=\mathbf 1$,
\begin{align} \label{eq-a1}
%\end{align}
%and
%\begin{align}\label{eq-a1}
%\text{for}~\aaa'=\mathbf 1,\quad
&\begin{cases}
p_{(0,\mathbf1)}+p_{(1,\mathbf1)} =\bar p_{(0,\mathbf1)}+\bar p_{(1,\mathbf1)}, \\
g_1 p_{(0,\mathbf1)}+c_1 p_{(1,\mathbf1)} =\bar g_1\bar p_{(0,\mathbf1)}+\bar c_1 \bar p_{(1,\mathbf1)}, \\
g_2 p_{(0,\mathbf1)}+c_2 p_{(1,\mathbf1)} =\bar g_2\bar p_{(0,\mathbf1)}+\bar c_2 \bar p_{(1,\mathbf1)}, \\
g_1 g_2 p_{(0,\mathbf1)}+c_1 c_2 p_{(1,\mathbf1)} =\bar g_1 \bar g_2\bar p_{(0,\mathbf1)}+\bar c_1 \bar c_2 \bar p_{(1,\mathbf1)}.\\
\end{cases}
\end{align}
The above \eqref{eq-anot1} obviously leads to $\bar g_2 = g_2$, and the last two equations of \eqref{eq-anot1} are automatically satisfied if the first two of \eqref{eq-anot1} are satisfied. Then the last two equations of \eqref{eq-a1} can be transformed to
\begin{align*}
\begin{cases}
(c_2-g_2)p_{(1,\mathbf1)} = (\bar c_2 - g_2)\bar p_{(1,\mathbf1)}, \\
c_1(c_2-g_2)p_{(1,\mathbf1)} = \bar c_1(\bar c_2 - g_2)\bar p_{(1,\mathbf1)};
\end{cases}
\end{align*}
which gives $\bar c_1 = c_1$. Additionally, when $\bar c_1 = c_1$, we also have that the last equality of \eqref{eq-a1} holds as long as the first three equalities of \eqref{eq-a1} hold. In summary, now there are $2^K+2$ parameters to be determined, which are
$\{\bar g_1,\bar c_2\}\cup \{\bar p_{\aaa}:\,\aaa\in\{0,1\}^K\},$
while they only have to satisfy the following $2\times(2^{K-1}-1)+3 = 2^K+1$ constraints,
\begin{align*}
\forall \aaa'\in\{0,1\}^{K-1},\text{  for}~\aaa'\neq\mathbf 1,\quad
&\begin{cases}
p_{(0,\aaa')}+p_{(1,\aaa')} =\bar p_{(0,\aaa')}+\bar p_{(1,\aaa')}, \\
g_1 p_{(0,\aaa')}+c_1 p_{(1,\aaa')} =\bar g_1\bar p_{(0,\aaa')}+ c_1 \bar p_{(1,\aaa')}; \\
\end{cases} \\
\\
\text{and for}~\aaa'=\mathbf 1,\quad
&\begin{cases}
p_{(0,\mathbf1)}+p_{(1,\mathbf1)} =\bar p_{(0,\mathbf1)}+\bar p_{(1,\mathbf1)}, \\
g_1 p_{(0,\mathbf1)}+c_1 p_{(1,\mathbf1)} =\bar g_1\bar p_{(0,\mathbf1)}+ c_1 \bar p_{(1,\mathbf1)}, \\
g_2 p_{(0,\mathbf1)}+c_2 p_{(1,\mathbf1)} = g_2\bar p_{(0,\mathbf1)}+\bar c_2 \bar p_{(1,\mathbf1)}. \\
\end{cases}
\end{align*}
Since the number of free variables $2^K+2$ is greater than the number of constraints $2^K+1$, there exist infinitely many different solutions to the above system of equations. This means that the $(Q,\cs,\cg,\pp)$ are not generically identifiable.
In particular, one can arbitrarily choose $\bar g_1$ close to but not equal to $g_1$, then solve for the remaining parameters $\{\bar p_{\aaa},~\aaa\in\{0,1\}^K\}$ and $\bar c_2$ as follows,
\begin{align*}
&\forall\aaa'\in\{0,1\}^{K-1},\quad
\begin{cases}
\bar p_{(0,\aaa')} = p_{(0,\aaa')}(g_1-c_1)/(\bar g_1 - c_1),\\
\bar p_{(1,\aaa')} =  p_{(0,\aaa')} +  p_{(1,\aaa')} -\bar p_{(0,\aaa')};\\
\end{cases}\\
%&\text{and}\\
& \bar c_2 = \frac{g_2[p_{(0,\mathbf1)} - \bar p_{(0,\mathbf1)}] + c_2 p_{(1,\mathbf1)}}{\bar p_{(1,\mathbf1)}}.
\end{align*}
This concludes the proof of part (a) of the theorem.

\bigskip
Next we first prove (b.2), i.e. when $Q^\star$ has two submatrices $\mathcal I_{K-1}$. In this case, the $Q$ contains a submatrix of the form $(I_K, I_K)^\top$.
%\begin{equation}\label{eq-2IK}
%\begin{pmatrix}
%I_K \\
%I_K \\
%\end{pmatrix}.
%\end{equation}
The proof of (b.1), i.e. when $Q^\star$ satisfies Conditions $A$, $B$ and $C$, is combined with the proof of part (c) later.

\bigskip
\noindent
\textbf{Proof of Part (b.2).}
We first give the proof when $Q$ only consists of two $I_K$'s, namely $Q=(I_K,I_K)^\top$.
In this case, we first prove that $\bar Q\sim Q$ must hold, using an argument similar to Step 1 of the proof of Theorem \ref{thm-dina}. Suppose $T(Q,\cc,\cg)\pp = T(\bar Q,\bar\cc,\bar\cg)\bar\pp$. Since $Q_{(K+1):(2K),\Cdot}=I_K$, we have that for each attribute $h\in[K]$, there is 
$$
\max_{m\in\{K+1,\ldots,2K\},\atop m\neq K+h}\,q_{m,h}=0,\quad
\max_{m\in\{K+1,\ldots,2K\}}\,q_{m,k}=1~\forall k\in[K]\setminus\{h\}.
$$ 
Therefore we can apply Lemma \ref{lem-H1} with $S=\{K+1,\ldots,2K\}\setminus\{K+h\}$ and $\mathcal J=[K]\setminus\{h\}$ to obtain 
$$
\max_{k\in\mathcal J}\,\bar\qq_k
\nsucceq\bar\qq_h.
$$
This essentially implies that for an arbitrary $h\in[K]$, there must be a $m_h\in[K]$ such that $\bar q_{h,m_h}=0$ and $\bar q_{k,m_h}=0$ for all $k\in[K]\setminus\{h\}$.
Moreover, the $K$ integers $m_1,m_2,\ldots,m_K$ must all be distinct, otherwise it is easy to see $\max_{k\in\mathcal J}\,\bar\qq_k\nsucceq\bar\qq_h$ would fail to hold for some $h\in[K]$.
So $(m_1,m_2,\ldots,m_K)$ is a permutation of $(1,2,\ldots,K)$.
Now we have obtained that $\bar Q_{1:K,(m_1,\ldots,m_K)}$ must be an identity matrix, i.e., $\bar Q_{1:K,\Cdot}\sim Q_{1:K,\Cdot}$.
Reasoning in exactly the same way gives $\bar Q_{(K+1):(2K),\Cdot}\sim Q_{(K+1):(2K),\Cdot}$, and we have $\bar Q\sim Q$.
Now for an arbitrary $\aaa=(\alpha_1,\alpha_2,\ldots,\alpha_K)\equiv(\alpha_1,\aaa')$, define 
\begin{align*}
\ttt^* = & \bar g_1\ee_1 + \bar c_{K+1}\ee_{K+1}
+ \sum_{k>1:\, \alpha_k=1} g_k\ee_k + \sum_{k>1:\, \alpha_k=0} c_k\ee_k\\
\equiv & \bar g_1\ee_1 + \bar c_{K+1}\ee_{K+1} + \ttt^{\aaa}
%\quad\rr^* = \ee_1 + \ee_{K+1}
\end{align*}
then $T_{\ee_1+\ee_{K+1}}(Q,\bar \cs-\ttt^*, \bar \cg-\ttt^*)=\mathbf 0$, so
\begin{align*}
0=&T_{\ee_1+\ee_{K+1}}(Q,\bar \cs-\ttt^*, \bar \cg-\ttt^*)\bar\pp
=T_{\ee_1+\ee_{K+1}}(Q,\cs-\ttt^*, \cg-\ttt^*)\pp\\
=& \prod_{k>1:\, \alpha_k=1}(c_k-g_k) \times \prod_{k>1:\, \alpha_k=0} (g_k-c_k)\times\\
&\Big[(g_1-\bar g_1)(g_{K+1} - \bar c_{K+1})p_{(0,\alpha_2,\ldots,\alpha_K)}+(c_1-\bar g_1)(c_{K+1} - \bar c_{K+1})p_{(1,\alpha_2,\ldots,\alpha_K)}\Big].
\end{align*}
This implies that for any $\aaa'=(\alpha_2,\ldots,\alpha_K)\in\{0,1\}^{K-1}$, we have 
$$
(g_1-\bar g_1)(g_{K+1} - \bar c_{K+1})p_{(0,\alpha_2,\ldots,\alpha_K)}+(c_1-\bar g_1)(c_{K+1} - \bar c_{K+1})p_{(1,\alpha_2,\ldots,\alpha_K)}=0.
$$
Since $g_{K+1}-\bar c_{K+1}\neq 0$, we have that 
$$
g_1-\bar g_1 = \frac{(c_1-\bar g_1)(c_{K+1} - \bar c_{K+1})p_{(1,\aaa')}}{(\bar c_{K+1} - g_{K+1})p_{(0,\aaa')}}, \quad \text{for any}~\aaa'\in\{0,1\}^{K-1}. 
$$
This equality indicates that if there exists $\aaa'_1\neq\aaa'_2$ such that
\begin{equation}\label{eq-ratio1}
\frac{p_{(1,\aaa'_1)}}{p_{(0,\aaa'_1)}}\neq \frac{p_{(1,\aaa'_2)}}{p_{(0,\aaa'_2)}},
\end{equation}
then one must have 
$$c_{K+1}-\bar c_{K=1}= 0,\quad g_1-\bar g_1=0.$$
Redefine $\ttt^* = \bar c_1\ee_1 + \bar g_{K+1}\ee_{K+1} +\ttt^{\aaa}$, then following the same procedure as above one gets that if $\pp$ satisfy \eqref{eq-ratio1}, then $g_{K+1}-\bar g_{K=1}= 0$ and  $c_1-\bar c_1=0.$

Similarly as the above procedure for $k=1$, we have that if for any attribute $k\in\{1,\ldots,K\}$, there exist  two attribute profiles $\aaa^{k,1},\aaa^{k,2}\in \{0,1\}^{k-1}\times\{0\}\times \{0,1\}^{K-k-1}$ such that 
\begin{equation}\label{eq-ratiok}
\frac{p_{\aaa^{k,1}+\ee_k}}{p_{\aaa^{k,1}}}\neq \frac{p_{\aaa^{k,2}+\ee_k}}{p_{\aaa^{k,2}}},
\end{equation}
then 
$$
\bar g_k = g_k,\quad \bar c_k = c_k~,\quad\bar g_{K+k} = g_{K+k},\quad \bar c_{K+k} = c_{K+k}\quad\text{for every}~k\in\{1,\ldots,K\}.
$$
Now that all the item parameters are identified under \eqref{eq-ratiok}, Equation \eqref{eq-tpr} gives $\bar \pp = \pp$. Therefore other than the measure zero set of the parameter space specified by constraints \eqref{eq-ratiok}, $(Q,\cs,\cg,\pp)$ are identifiable. This means $(Q,\cs,\cg,\pp)$ are generically identifiable.

In particular, if $Q$ takes form of the $Q_{2\times 4}$ in \eqref{eq-q24},
$$
Q_{2\times 4}=\begin{pmatrix}
1 & 0 \\
0 & 1 \\
1 & 0 \\
0 & 1 \\
\end{pmatrix},
$$
then constraints \eqref{eq-ratiok} just simplify to
$$
\frac{p_{(10)}}{p_{(00)}}\neq \frac{p_{(11)}}{p_{(01)}}
\quad\text{and}\quad
\frac{p_{(01)}}{p_{(00)}}\neq \frac{p_{(11)}}{p_{(10)}},
$$
which can be equivalently written as inequality \eqref{eq-q24-id} that $p_{(01)}p_{(10)}\neq p_{(00)}p_{(11)}$  in the main text.

Next we prove the conclusion when $Q$ contains other rows besides the two identity submatrices, namely $Q=(I_K,I_K,(Q^\star)^\top)^\top$.
Using exactly the same arguments as previously we have that generically, all the item parameters of the first $2K$ items as well as all the proportion parameters are satisfied. Now for any $J>2K$ and $\aaa\in\{0,1\}^K$ define $\rr^*=\sum_{k=1}^K \ee_j$ and 
$$
\ttt^* = \sum_{1\leq k\leq K:\, \aaa\succeq\qq_j} g_j\ee_j + \sum_{1\leq k\leq K:\, \aaa\nsucceq\qq_j} c_j\ee_j,
$$
then \eqref{eq-trans} implies that
$$
\theta_{j,\aaa} = 
\frac{T_{\rr^*+\ee_j}(Q,\cc-\ttt^*, \cg-\ttt^*)\pp}{T_{\rr^*}(Q,\cc-\ttt^*, \cg -\ttt^*)\pp} = 
\frac{T_{\rr^*+\ee_j}(\bar Q,\bar\cc-\ttt^*, \bar\cg-\ttt^*)\bar\pp}{T_{\rr^*}(\bar Q,\bar\cc-\ttt^*, \bar\cg -\ttt^*)\bar\pp}
= \bar\theta_{j,\aaa}.
$$
This proves that any slipping or guessing parameter associated with item $j>2K$ is identifiable under the generic constraints \eqref{eq-ratiok}, and this completes the proof of part (b.2)  of the theorem.

\bigskip

Next we   prove   (b.1) and   (c) in Theorem 2 in four steps. 

\noindent
\textbf{Proof of Part (b.1) and Part (c).}

\noindent\textbf{Step 1.}
In this step,  
we aim   to show that if 
\begin{equation}\label{eq-tpr}
T_{\rr,\Cdot}(Q,\cs,\cg)\pp = T_{\rr,\Cdot}(\bar Q,\bar\cs,\bar\cg)\bar\pp\quad\text{for every}~\rr\in\{0,1\}^J,
\end{equation}
then $\bar Q$ must take the following form up to column permutation
\begin{equation}\label{eq-qqbar}
%Q=\begin{pmatrix}
%1 & \mathbf0 \\
%1 & \vv \\
%\mathbf0 & Q^\star
%\end{pmatrix},\qquad\qquad
\bar Q=\begin{pmatrix}
1 & \mathbf0 \\
\bar u & \bar\vv \\
\mathbf0 & Q^\star
\end{pmatrix}.
\end{equation}
Here $(\bar u, \bar\vv)$ is a $K$ dimensional binary vector. The structure of $(\bar u, \bar\vv)$ will be studied in Steps 2 and 3.

Since the submatrix $Q^\star$ of $Q$ satisfies Conditions $A$, $B$ and $C$, the matrix $Q$ can be written as
$$
Q = \begin{pmatrix}
1 & \mathbf 0^\top \\
1 & \vv^\top \\
\mathbf 0 & \mathcal I_{K-1} \\
\mathbf 0 & Q^{\star\star}
\end{pmatrix},
$$  
then follow the same procedure as Step 1 in the proof of Theorem \ref{thm-dina} one has that, up to some column permutation, $\bar Q$ takes the form 
$$
\bar Q = \begin{pmatrix}
1 & \mathbf 0^\top \\
\bar u & \bar\vv^\top \\
\mathbf 0 & \mathcal I_{K-1} \\
\bar\bb & \bar Q^{\star\star}
\end{pmatrix}.
$$
For notational convenience and without loss of generality, in the following proof we rearrange the order of the row vectors of  $Q$ (and $\bar Q$) and rewrite them as follows %(put the second row to the $(K+1)$th row)
\begin{equation}\label{eq-permute}
Q = \begin{pmatrix}
1 & \mathbf 0^\top \\
\mathbf 0 & \mathcal I_{K-1} \\
1 & {\vv}^\top \\
\mathbf 0 & {Q^{\star\star}}
\end{pmatrix}
%=\begin{pmatrix}
%1 & \mathbf 0 \\
%\mathbf 0 & \mathcal I_{K_1} \\
%\vv & Q^{sub} \\ 
%\end{pmatrix}
,\qquad
\bar Q = \begin{pmatrix}
1 & \mathbf 0^\top \\
\mathbf 0 & \mathcal I_{K-1} \\
\bar u & \bar\vv^\top \\
\bar\bb & \bar Q^{\star\star}
\end{pmatrix}.
\end{equation}
Now that each column of $Q^{\star\star}$ contains at least two entries of ``$1$'' from the assumption of scenarios (b.1) and (c), following the same procedure as Step 2 in the proof of Theorem \ref{thm-dina} we can obtain
$$c_j = \bar c_j,\quad \text{for } j=K+2,\ldots,J.$$
Note that slightly different from Step 2 in the proof of Theorem \ref{thm-dina}, here we \textit{do not have} $c_{K+1}=\bar c_{K+1}$ due to the fact that the first attribute is  required by only two items.

%Next we follow a similar procedure as that of Step 3 in the proof of Theorem \ref{thm-dina},
Now denote the $(J-K)\times(K-1)$ bottom-right submatrix of $Q$ by $Q^{s}$ and the $(J-K)\times K$ bottom submatrix of $Q$ by $Q^l$, i.e.,
$$
Q^{s}=
\begin{pmatrix}
\vv^\top\\
Q^{\star\star}
\end{pmatrix},\quad
Q^l = \begin{pmatrix}
1 & {\vv}^\top \\
\mathbf 0 &  {Q^{\star\star}}
\end{pmatrix}
$$
and assume without loss of generality that {the $K-1$ column vectors of $Q^{s}$ are arranged in the lexicographic order.} Specifically, for any $1\leq k_1<k_2\leq K-1$,  assume $Q^s_{\Cdot,k_1} \prec_{\text{lex}} Q^s_{\Cdot,k_2}$.
%\begin{align*}
%\text{either }~ & Q^{s}_{1, \,k_1} < Q^{s}_{1,\, k_2};\\
%\text{or }~ &  Q^{s}_{j_0,\, k_1} < Q^{s}_{j_0,\, k_2} ~\text{and}~ 
%Q^{s}_{j,\, k_1} = Q^{s}_{j,\, k_2}~\text{for all}~j\leq j_0,
%\end{align*}
%and we denote this by $Q^s_{\Cdot,k_1}\prec_{\text{lex}} Q^s_{\Cdot,k_2}$.
This implies that the vector $\vv$ can be written as
$$
\vv = (0,\ldots,0, 1,\ldots,1)
$$
Note that in scenario (b.1), $\vv=\zero$ and $k_0=K-1$.
where its first $k_0$ elements are zero and the remain $K-1-k_0$ elements are one. So $\qq_2=(1,\vv)=(1,0,\ldots,0,1,\ldots,1)$. We now use an induction method to prove that
\begin{equation}\label{eq-gk}
g_k = \bar g_k,\quad \forall k=2,\, \ldots, \, 1+k_0.
\end{equation}
A key observation is that if considering the order of the columns of the larger submatrix $Q^l$ instead of $Q^s$, then the first column of $Q^l$, i.e. $Q^l_{\Cdot, 1}$ is of larger lexicographic order of $Q^l_{\Cdot, k}$ for any  $k=2,\ldots, 1+k_0$. This indicates that we can follow a similar induction argument as Step 3 in the proof of Theorem \ref{thm-dina} by defining $\ttt^*_k$
as (the same form as \eqref{eq-ggc})
\begin{equation}\label{eq-ggck}
\ttt^*_k = \sum_{h=1}^K\bar g_h\ee_h +\sum_{h>K:\,q_{h,k}=0} g_h\ee_h + \sum_{h>K:\,q_{h,k}=1} c_h\ee_h,
\end{equation}
for $k=2,\ldots, 1+k_0$ one after another, to obtain \eqref{eq-gk}.

\textit{We emphasize here that if $\vv=\mathbf0$, i.e. in scenario (b.1) of the theorem, %with $Q^\star$ satisfying Conditions $A$, $B$ and $C$, 
then $k_0=K-1$ and by far we have already obtained $\bar g_k=g_k$ for all $k=2,\ldots,K$.} So we can directly go to the next step, Step 2 of the proof, without  the local condition to appear in \eqref{eq-induc} later. That is why in scenario (b.1) of the theorem, we have  \textit{global} generic identifiability of $(Q,\cs,\cg,\pp)$.

%Now consider $\bar g_{1}$ and $g_{1}$ and define 
%$$
%\ttt^*_{k_0+2} = \sum_{h=1}^K\bar g_h\ee_h +\sum_{h=K+2}^J g_h\ee_h + \bar c_{K+1}\ee_{K+1},
%$$
%then $T_{\sum_{k}\ee_k,\Cdot}(\bar Q,\bar\cc-\ttt^*_{k_0+2},\bar\cg-\ttt^*_{k_0+2}) = \mathbf 0$ and $T_{\sum_{k}\ee_k,\Cdot}(\bar Q,\bar\cc-\ttt^*_{k_0+2},\bar\cg-\ttt^*_{k_0+2})\bar\pp=0$, so 
%\begin{equation}\label{eq-k02}
%T_{\sum_{k}\ee_k,\Cdot}( Q,\cc-\ttt^*_{k_0+2},\cg-\ttt^*_{k_0+2})\pp=0
%\end{equation}
% must hold. Since $Q^{\star\star}$ contains at least two entries of one in each column, the row vector $T_{\sum_{k}\ee_k,\Cdot}( Q,\cc,\cg)$ contains only two potentially nonzero elements corresponding to attribute profiles $\aaa_0:=(0,1,\ldots,1)$ and $\aaa_1:=(1,1,\ldots, 1)$. Then \eqref{eq-k02} can be written as
%$$
%\prod_{h=2}^K(c_h - \bar g_h)\prod_{h=K+2}^J(c_h - g_h)\Big[
%(g_1-\bar g_1)(g_{K+1} - \bar c_{K+1}) p_{01\cdots1} + (c_1-\bar g_1)(c_{K+1} - \bar c_{K+1}) p_{11\cdots1} \Big]=0,
%$$
%which further reduces to 
%\begin{equation}\label{eq-gc-tradeoff}
%(g_1-\bar g_1)(g_{K+1} - \bar c_{K+1}) p_{01\cdots1} + (c_1-\bar g_1)(c_{K+1} - \bar c_{K+1}) p_{11\cdots1}=0.
%\end{equation}
%This equality indicates that $g_1=\bar g_1$ if and only if $c_{K+1}=\bar c_{K+1}$, and generally,
%$$
%\bar c_{K+1} = \frac{g_{K+1}(g_1-\bar g_1)p_{01\cdots1} + c_{K+1}(c_1-\bar g_1)p_{11\cdots1}}{(g_1-\bar g_1)p_{01\cdots1} +(c_1-\bar g_1)p_{11\cdots1}}.
%$$

Next we consider the case $\vv\neq\zero$, i.e. in scenario (c) of the theorem, then $k_0<K-1$. We will use another induction argument to show $\bar g_{k}=g_k$ for $k=k_0+2, \ldots,K$, under an additional local condition.
First we consider $\bar g_k$ and $g_{k}$ for $k=k_0+2$. Note that $Q_{\Cdot,k} \succ_{\text{lex}} Q_{\Cdot, 1}$, and $Q_{\Cdot,k} \prec_{\text{lex}} Q_{\Cdot, m}$ for any $m=k+1,\ldots, K$. Define $\ttt^*_k$ the same as in \eqref{eq-ggck}, then $T_{\rr^*,\Cdot}(\bar Q,\bar\cc-\ttt^*_k,\bar\cg-\ttt^*_k) = \mathbf0$ and $T_{\rr^*,\Cdot}(\bar Q,\bar\cc-\ttt^*_k,\bar\cg-\ttt^*_k) \bar\pp=0$, so
$
T_{\rr^*,\Cdot}( Q,\cc-\ttt^*_k,\cg-\ttt^*_k)\pp=0.
$
We claim that in the the vector $T_{\rr^*,\Cdot}( Q,\cc-\ttt^*_k,\cg-\ttt^*_k)$, denoted by $T_{\rr^*,\Cdot}$ afterwards for notational simplicity, only contains two potentially nonzero elements corresponding to attribute profiles $\aaa_{1k} = \sum_{m=1}^K\ee_m - \ee_k=(1,\ldots,1,\alpha_k=0,1,\ldots,1)$ and $\aaa_{0k} = \aaa_1 - \ee_1=(\alpha_1=0,1,\ldots,1,\alpha_k=0,1,\ldots,1)$. This is because on the one hand, for any attribute profile $\aaa$ that lacks some attribute $m\in\{k+1,\ldots,K\}$, $\theta_{h,\aaa} = g_h$ for some item $h>K$ with $q_{h,k}=0$, which makes $T_{\rr^*,\aaa}=0$; and on the other hand, for any attribute profile that lacks some attribute $m\in\{2,\ldots, k-1\}$, since we already have \eqref{eq-gk}, $\theta_{h,m}=g_h=\bar g_h$ for some $h\in\{2,\ldots, K\}$, which makes $T_{\rr^*,\aaa}=0$. Now $T_{\rr^*,\aaa}\neq 0$ would only happen if 
$
\aaa = (\alpha_1,1,\ldots,1,\alpha_{k},1,\ldots,1).
$
However, if $\alpha_k=1$ and $\aaa = (\alpha_1,1,\ldots, 1)$, then $\theta_{h,\aaa}=c_h$ for some item $h>K$ with $q_{h,k}=1$, which also makes $T_{\rr^*,\aaa}=0$. Now we have proven the claim that $T_{\rr^*,\Cdot}$ has only two potentially nonzero elements corresponding to $\aaa_{1k}$ and $\aaa_{0k}$. Therefore we have for $k=k_0+2$,
\begin{align*}
0=&T_{\rr^*,\Cdot}( Q,\cc-\ttt^*_k,\cg-\ttt^*_k) \pp\\
=&\prod_{h=2}^K(c_h-\bar g_h)\prod_{h>K:\atop q_{h,k}=0}(c_h-g_h)\prod_{h>K:\atop q_{h,k}=1}(g_h-c_h)\\
&\times\Big[ (g_1-\bar g_1)p_{\aaa_{0k}} + (c_1-\bar g_1) p_{\aaa_{1k}} \Big]
(g_k-\bar g_k),
\end{align*}
which further gives
\begin{equation}\label{eq-a0k}
\Big[(g_1-\bar g_1)p_{\aaa_{0k}} + (c_1-\bar g_1) p_{\aaa_{1k}} \Big]
(g_k-\bar g_k) = 0~\text{for}~k=k_0+2.
\end{equation}
Note that if $\bar g_1=g_1$, then the part in the bracket in the above display becomes $(c_1-g_1)p_{\aaa_1}$, which is nonzero. Therefore, when $\bar g_1$ is sufficiently close to  the true parameter $g_1$, the part in the bracket in \eqref{eq-a0k} would be nonzero. We formally write it as
\begin{align}\label{eq-local}
\text{for}~k=k_0+2,\quad&\forall \bar g_1\in\mathcal{N}_k,\quad (g_1-\bar g_1)p_{\aaa_{0k}} + (c_1-\bar g_1)p_{\aaa_{1k}} \neq 0, \\ \notag
&\text{where}~\mathcal{N}_k = \{x: 0<x<\frac{g_1 p_{\aaa_{0k}} + c_1 p_{\aaa_{1k}}}{p_{\aaa_{0k}} + p_{\aaa_{1k}}}\}.
\end{align} 
This indicates that in the neighborhood $\mathcal N_k$ of $g_1$, \eqref{eq-a0k} leads to $g_k = \bar g_k$ for $k=k_0+2$.

Then we use induction to prove $g_k = \bar g_k$ for all $k=k_0+3,\ldots, K$. As the induction assumption, assume that when 
%\begin{equation*}
$\bar g_1\in\bigcap_{m=k_0+2}^{k-1}\mathcal N_m$
%\end{equation*} 
holds, we have $g_m = \bar g_m$ for all $m=2,\ldots, k-1$. Then define $\ttt^*$ the same as in \eqref{eq-ggck}, and deduce in the same way as in proving $g_{k_0+2}=\bar g_{k_0+2}$, we have 
\begin{equation*}
\Big[(g_1-\bar g_1)p_{\aaa_{0k}} + (c_1-\bar g_1) p_{\aaa_{1k}} \Big]
(g_k-\bar g_k) = 0,
\end{equation*}
and further for any $\bar g_1\in\mathcal N_k$ (more accurately any $\bar g_1\in\big[\cap_{m=k_0+2}^{k-1}\mathcal N_m \big]\cap \mathcal N_k $), we must have $\bar g_k=g_k$. Here $\mathcal N_k$ takes the same form as that in \eqref{eq-local}.
Now by induction, we have that if
\begin{equation}\label{eq-induc}
\bar g_1\in\bigcap_{m=k_0+2}^{K}\mathcal N_m,
\quad 
\end{equation}
then $g_k = \bar g_k$ for $k=k_0+2, \ldots, K$. Combined with the previous results shown in \eqref{eq-gk}, now we have proven that in scenario (c) of the theorem, if the local condition \eqref{eq-induc} is satisfied, then $\bar g_k = g_k$ for $k=2, \ldots, K$.

In summary, we have shown $\bar g_k=g_k$ for $k=2, \ldots, K$ (under \eqref{eq-induc} if in scenario (c)) and $\bar c_j = c_j$ for $j=K+2,\ldots,J$. Based on these, following  similar procedures as in Step 5 of the proof of Theorem \ref{thm-dina}, we obtain that 
$$
\bar \qq_j = \qq_j,\quad\forall j=K+2,\ldots,J.
$$
%Finally, we have shown \eqref{eq-qqbar} holds. %Then it is easy to check that \eqref{eq-q1q2} also holds.
%and $\bar g_k = g_k$ for $k=2,\ldots,K$ and $\bar c_j = c_j$ for $j=K+2,\ldots,J$. 

%\begin{lemma}\label{lem-zero}
%Consider two $Q$-matrices $Q$ and $\bar Q$ shown below,
%$$
%Q = \begin{pmatrix}
%Q_{11} & Q_{12} \\
%\zero & Q_{22} \\
%\end{pmatrix}
%,\qquad
%\bar Q = \begin{pmatrix}
%\bar Q_{11} & \bar Q_{12} \\
%\zero & Q_{22} \\
%\end{pmatrix},
%$$
%where the bottom-left matrix $\zero$ denotes a zero matrix of certain size. Then if 
%\end{lemma}

\medskip
\noindent\textbf{Step 2.}
In this step we show $\bar u=1$ in \eqref{eq-qqbar}.
If $\bar u=0$, set
$$\ttt^* = c_1\ee_1 + \bar c_2\ee_2+ \sum_{j=3}^{K+3}g_k\ee_k,\quad \rr^*=\sum_{j=1}^{K+3}\ee_j,$$
then 
\begin{align*}
T_{\rr^*,\Cdot}(\bar Q,\bar \cc-\ttt^*, \bar\cg - \ttt^*)\bar \pp &= \mathbf 0^\top \cdot \bar \pp=0, \\
T_{\rr^*,\Cdot}( Q, \cc-\ttt^*, \cg - \ttt^*) \pp &= (g_1-c_1)(g_2-\bar c_2)\prod_{j=3}^{K+3}(c_j-g_j)p_{(0,1,\ldots,1)}\neq 0,
\end{align*}
which contradicts Equation \eqref{eq-def-T}. So $\bar u=1$. 
Now we have obtained
\begin{equation}\label{eq-step2}
Q = \begin{pmatrix}
1 & \mathbf 0^\top \\
\mathbf 0 & \mathcal I_{K-1} \\
1 & {\vv}^\top \\
\zero & {Q^{\star\star}}
\end{pmatrix}
,\qquad
\bar Q = \begin{pmatrix}
1 & \mathbf 0^\top \\
\mathbf 0 & \mathcal I_{K-1} \\
1 & \bar\vv^\top \\
\zero & Q^{\star\star}
\end{pmatrix}.
\end{equation}

\medskip
\noindent\textbf{Step 3.} In this step we show $\bar\vv=\vv$. For notational simplicity in the following proof, we rearrange the order of the row vectors in $Q$ and $\bar Q$ in \eqref{eq-step2} again to the following forms
\begin{equation}\label{eq-step3}
Q = \begin{pmatrix}
1 & \mathbf 0^\top \\
1 & {\vv}^\top \\
\mathbf 0 & \mathcal I_{K-1} \\
\zero & {Q^{\star\star}}
\end{pmatrix}
,\qquad
\bar Q = \begin{pmatrix}
1 & \mathbf 0^\top \\
\bar u & \bar\vv^\top \\
\mathbf 0 & \mathcal I_{K-1} \\
\zero & Q^{\star\star}
\end{pmatrix},
\end{equation}
and our conclusions proved so far are $\bar g_k=g_k$ for $k=3,\ldots,K+1$  and $\bar c_j=c_j$ for $j=K+2,\ldots,J$ (under the local condition \eqref{eq-induc} if in scenario (b.1)).
%Next it suffices to show that when $\vv\neq \mathbf 1$, parameters $c_1,c_2,g_1,g_2$ as well as $\pp$ are generically identifiable. %given a fixed $Q$ in the form of \eqref{eq-prop1}. 
Given that the last $J-2$ rows of $Q$ and $\bar Q$ are equal,  %following the same argument as the proof of Theorem 2 of \cite{partial} we have the following equalities,
%Now consider an arbitrary response pattern $\rr=(r_1,r_2,r_3,\ldots,\allowbreak r_J)$. 
we claim that \eqref{eq-tpr} for response pattern $\rr$ can be equivalently written as
\begin{align}\label{eq-qprime-n0}
%\begin{aligned}
&\sum_{\aaa'\in\atop\{0,1\}^{K-1}}\prod_{j>2\atop\,r_j=1}\theta_{j,\,(0,\aaa')}\cdot \mathbb P(R_1\geq r_1,\,R_2\geq r_2,\, \ba_{2:K}=\aaa'\mid  Q,\TT,\pp)\\
=
&\sum_{\aaa'\in\atop\{0,1\}^{K-1}}\prod_{j>2\atop\,r_j=1}\bar\theta_{j,\,(0,\aaa')}\cdot {\mathbb P}(R_1\geq r_1,\,R_2\geq r_2,\,\ba_{2:K}=\aaa'\mid \bar Q,\bar\TT,\bar\pp).\notag
%\end{aligned}
\end{align}
%$f_{(r_1,r_2),\,\aaa^*} = P(R_1=r_1,R_2=r_2\mid \aaa_{2:K}=\aaa^*)$
%where $f_{(r_1,r_2),\,\aaa^*}$ and $\bar f_{(r_1,r_2),\,\aaa^*}$ are some quantities depending on $\aaa^*$ and $(r_1,r_2)$, the first two entries of $\rr$.
%where $\mathbb P(R_1\geq r_1,\,R_2\geq r_2,\, \ba_{2:K}=\aaa')$ represents the probability of $\{R_1\geq r_1,\,R_2\geq r_2\}$ and the random attribute profile $\ba$ has its last $K-1$ entries being $\aaa'$ under $(Q,\cs,\cg,\pp)$, while $\overline {\mathbb P}(R_1\geq r_1,\,R_2\geq r_2,\, \ba_{2:K}=\aaa')$ represents that under $(\bar Q,\bar \cs,\bar \cg,\bar \pp)$.
Here $\ba=(A_1,\ldots,A_K)$ denotes a random attribute profile following a categorical distribution with proportion parameters $\pp$, and $\ba_{2:K}$ denotes the vector consisting of the last $K-1$ elements of $\ba$.
The reason for the equivalence of \eqref{eq-qprime-n0} and \eqref{eq-tpr} is stated as follows. Since all items other than the first two do not require the first attribute, we have that for any $\aaa'\in\{0,1\}^{K-1}$, the two attribute profiles $(0,\aaa')$ and $(1,\aaa')$ always have the same response probability $\theta_{j,(0,\aaa')}$ to any item $j>2$. This indicates that the left hand side of \eqref{eq-tpr} can be written as 
$$T_{\rr,\Cdot}(Q,\cs,\cg)\pp=\sum_{\aaa'\in\atop\{0,1\}^{K-1}}\prod_{j>2\atop\,r_j=1}\theta_{j,\,(0,\aaa')}\cdot \mathbb P(R_1\geq r_1,\,R_2\geq r_2,\, \ba_{2:K}=\aaa'\mid  Q,\TT,\pp),$$
and this further leads to the equivalence between \eqref{eq-tpr} and \eqref{eq-qprime-n0}.
In particular, when $(r_1,r_2)=(0,0)$, we have $\mathbb P(R_1\geq r_1,\,R_2\geq r_2,\, \ba_{2:K}=\aaa'\mid  Q,\TT,\pp) = p_{(0,\aaa')}+p_{(1,\aaa')}$.
Now for any $J$-dimensional response pattern $\rr$ with $(r_1,r_2)=(0,0)$, then the constraint $T_{\rr,\Cdot}(Q,\cc,\cg)\pp = T_{\rr,\Cdot}(\bar Q,\bar\cc,\bar\cg)\bar\pp$ simply becomes
$$\sum_{\aaa'\in\atop\{0,1\}^{K-1}}\prod_{j>2\atop\,r_j=1}\theta_{j,\,(0,\aaa')}\cdot (p_{(0,\aaa')}+p_{(1,\aaa')})
=\sum_{\aaa'\in\atop\{0,1\}^{K-1}}\prod_{j>2\atop\,r_j=1}\bar\theta_{j,\,(0,\aaa')}\cdot (\bar p_{(0,\aaa')}+\bar p_{(1,\aaa')}).$$
Since the above equality holds for any $(r_3,r_4,\ldots,r_J)\in\{0,1\}^{J-2}$, we claim that, parameters $\theta_{j,(0,\aaa')}$ and $\bar\theta_{j,(0,\aaa')}$ for $j=3,\ldots,J$ can be equivalently viewed as all the item parameters (slipping or guessing) associated with the submatrix $Q^\star$, while grouped proportion parameters $p_{(0,\aaa')}+p_{(1,\aaa')}$ and $\bar p_{(0,\aaa')}+\bar p_{(1,\aaa')}$ can be viewed as all the ``proportion parameters" associated with $Q^\star$. Since $Q^\star$ satisfy the sufficient conditions $A$, $B$, $C$ in Theorem \ref{thm-dina} for identifiability, by Theorem \ref{thm-dina} we conclude that 
$\theta_{j,(0,\aaa')}=\bar\theta_{j,(0,\aaa')}$ for any $j\in\{3,\ldots,J\}$ and any $\aaa'\in\{0,1\}^{K-1}$. This indicates $\bar c_k=c_k$ for $k=3,\ldots,K+1$ and $\bar g_j=g_j$ for $j=K+2,\ldots,J$.
%The reason \eqref{eq-tprr} can be equivalently written as \eqref{eq-qprime-n0} is that, given any $\aaa_{2:K}\in\mathcal R^{Q'}$ and any item $j\in\{3,\ldots,J\}$, the positive response probability of $[\alpha_1,\aaa_{2:K}]$ to item $j$ only depends on $\aaa^*$ part, regardless of the value of $\alpha_1$,  as shown in \eqref{eq-pp}. 
%Therefore the terms in $T(\TT)_{\rr,\Cdot}\nnu$ can be grouped in such a way that it becomes the summation over all the $\aaa_{2:K}\in\mathcal R^{Q'}$, exactly as presented in Equation \eqref{eq-qprime-n0}.

Then an important observation is that, fix any particular pair of $(r_1, r_2)\in\{0,1\}^2$,  quantities in \eqref{eq-qprime-n0} can be viewed parameters associated with the $(J-2)\times (K-1)$ matrix $Q^\star$, just similar to the argument in the previous paragraph. Specifically, $\theta_{j,(0,\aaa')}$ and $\bar\theta_{j,(0,\aaa')}$ for $j=3,\ldots,J$ are item parameters (slipping or guessing) associated with the  $Q^\star$, and  $\mathbb  P(R_1\geq r_1,\,R_2\geq r_2,\, \ba_{2:K}=\aaa'\mid  Q, \TT, \pp)$ and ${ \mathbb P}(R_1\geq r_1,\,R_2\geq r_2,\, \ba_{2:K}=\aaa'\mid \bar Q,\bar\TT,\bar\pp)$ for each $\aaa'\in\{0,1\}^{K-1}$ can be viewed as the ``proportion parameters" associated with $Q^\star$.
%$T_{\rr',\Cdot}(Q^\star,\TT_{Q^\star})\pp^{Q^\star}$
%taking $(r_1,r_2)$ to be $(0,0)$, $(0,1)$, $(1,0)$, $(1,1)$ in \eqref{eq-tpr} respectively, 
Now because the submatrix $Q^\star$ satisfy the identifiability conditions $A$, $B$, $C$; and $\bar Q_{3:J,\Cdot}= Q_{3:J,\Cdot}=Q^\star$ and $\bar c_j=c_j$, $\bar g_j=g_j$ for $j=3,\ldots,J$, we must have
\begin{align}\label{eq-q1q2}
%&\forall (r_1,r_2)\in\{0,1\}^2,\quad 
\forall \aaa'\in\{0,1\}^{K-1},\qquad
&\mathbb  P(R_1\geq r_1,\,R_2\geq r_2,\, \ba_{2:K}=\aaa'\mid  Q, \TT, \pp)\\ \notag
=~& { \mathbb P}(R_1\geq r_1,\,R_2\geq r_2,\, \ba_{2:K}=\aaa'\mid \bar Q,\bar\TT,\bar\pp).
\end{align}
Now take $(r_1,r_2)$ to be $(0,0)$, $(0,1)$, $(1,0)$, $(1,1)$ in the above \eqref{eq-q1q2} respectively, we obtain
\begin{align}\label{eq-allv-dina}
\begin{cases}
 p_{(0,\aaa')}+  p_{(1,\aaa')} = \bar p_{(0,\aaa')}+\bar  p_{(1,\aaa')};\\
\theta_{1,\,(0,\aaa')}\cdot p_{(0,\aaa')} + \theta_{1,\,(1,\aaa')}\cdot p_{(1,\aaa')} = \bar\theta_{1,\,(0,\aaa')}\cdot\bar p_{(0,\aaa')} + \bar \theta_{1,\,(1,\aaa')}\cdot\bar p_{(1,\aaa')};\\
\theta_{2,\,(0,\aaa')}\cdot p_{(0,\aaa')} + \theta_{2,\,(1,\aaa')}\cdot p_{(1,\aaa')} = \bar\theta_{2,\,(0,\aaa')}\cdot\bar p_{(0,\aaa')} + \bar \theta_{2,\,(1,\aaa')}\cdot\bar p_{(1,\aaa')};\\
\theta_{1,\,(0,\aaa')}\theta_{2,\,(0,\aaa')}\cdot p_{(0,\aaa')} + \theta_{1,\,(1,\aaa')}\theta_{2,\,(1,\aaa')}\cdot p_{(1,\aaa')} 
\\
\qquad\qquad\qquad = \bar\theta_{1,\,(0,\aaa')}\bar\theta_{2,\,(0,\aaa')}\cdot\bar p_{(0,\aaa')} + \bar \theta_{1,\,(1,\aaa')}\bar\theta_{2,\,(1,\aaa')}\cdot\bar p_{(1,\aaa')}.
\end{cases}
\end{align}
Next we show $\vv=\bar\vv$.
\eqref{eq-allv-dina} implies that,
\begin{align*}
&\forall \aaa'\geq \vv,~\aaa'\ngeq\bar\vv,\quad
\begin{cases}
p_{(0,\aaa')}+p_{(1,\aaa')} =\bar p_{(0,\aaa')}+\bar p_{(1,\aaa')} \\
g_1 p_{(0,\aaa')}+c_1 p_{(1,\aaa')} =\bar g_1\bar p_{(0,\aaa')}+ \bar c_1 \bar p_{(1,\aaa')} \\
g_2 p_{(0,\aaa')}+c_2 p_{(1,\aaa')} =\bar g_2[\bar p_{(0,\aaa')}+\bar p_{(1,\aaa')} ]\\
g_1 g_2 p_{(0,\aaa')}+c_1 c_2 p_{(1,\aaa')} =\bar g_2[\bar g_1\bar p_{(0,\aaa')}+\bar c_1\bar p_{(1,\aaa')}] \\
\end{cases}
\end{align*}
If $\bar\vv\ngeq\vv$, then taking $\aaa'=\vv$ in the above equation and doing some transformation gives 
\begin{align*}
\begin{cases}
(g_2-\bar g_2)p_{(0,\aaa')}+(c_2-\bar g_2)p_{(1,\aaa')}=0,\\
(g_1-c_1)(g_2-\bar g_2)p_{(0,\aaa')}=0.
\end{cases}
\end{align*}
Since $g_1\neq c_1$, we have $g_2-\bar g_2=0$, which further gives $c_2-\bar g_2=0$. This contradicts $c_h>\bar g_h$ for any item $h$, so $\bar\vv\ngeq\vv$ can not happen. Similarly $\bar\vv\nleq\vv$ also can not happen, so $\bar\vv=\vv$.

\medskip
\noindent
\textbf{Step 4.} In the final step we show $c_1,c_2,g_1,g_2$ and $\pp$ are generically identifiable if $\vv\neq \mathbf1$.
First we show that if %the true model parameters $(\TT,\pp)$ satisfy the constraint that 
there exist $\aaa'_1$, $\aaa'_2\in\{0,1\}^{K-1}$, $\aaa'_1\neq\aaa'_2$ such that  
\begin{equation}\label{eq-s1s2}
%\frac{p_{(1,\aaa'_1)}}{p_{(0,\aaa'_1)}} \neq\frac{p_{(1,\aaa'_2)}}{p_{(0,\aaa'_2)}},
p_{(1,\aaa'_1)} p_{(0,\aaa'_2)} \neq p_{(1,\aaa'_2)}  p_{(0,\aaa'_1)},
\end{equation}
then one must have
\begin{equation}\label{eq-cg12}
%\theta_{j,(\alpha_1,\aaa'_1)} = \bar\theta_{j,(\alpha_1,\aaa'_1)},\quad\forall j=1,2,~~\forall\alpha_1 = 0,1.
c_i=\bar c_i, ~ g_i = \bar g_i,\quad i=1,2.
\end{equation}
After some transformations, the system of equations \eqref{eq-allv} yields
\[\begin{cases}
( g_1 -  c_1)\cdot( g_2 - \bar c_2)\cdot p_{(0,\aaa')}
=
(\bar  g_1 -  c_1)\cdot(\bar  g_2 - \bar c_2)\cdot \bar p_{(0,\aaa')},\\
( g_2 - \bar c_2)\cdot  p_{(0,\aaa')}  + ( c_2 - \bar c_2)\cdot \bar p_{(1,\aaa')}
= (\bar g_2 - \bar c_2)\cdot \bar p_{(0,\aaa')}.
\end{cases}\]
Since we have $\bar g_1\neq c_1$, the left hand side of the first equation above is nonzero. 
And obviously the right hand side of the second equation above is nonzero.
Taking the ratio of the above two equations gives
\[\begin{aligned}
% g_2 - \bar c_2 
%=&\frac{( g_1 -  c_1)\cdot(\bar  g_2 - \bar c_2)\cdot \bar p_{(0,\aaa')}}{(\bar g_1 -  c_1)\cdot \bar p_{(0,\aaa')} + (\bar c_1 -  c_1)\cdot \bar p_{(1,\aaa')}}\\
%%%%%%%%
& \frac{( g_1 -  c_1)\cdot( g_2 - \bar c_2)}{( g_2 - \bar c_2)  + ( c_2 - \bar c_2)\cdot  p_{(1,\aaa')} /  p_{(0,\aaa')}}
=(\bar  g_1 -  c_1)\equiv f(\aaa').
\end{aligned}\]
The right hand side of the above display does not involve any proportion parameter $\pp$ or $\bar\pp$. So for $\aaa'_1$, $\aaa'_2$ satisfying \eqref{eq-s1s2}, $f(\aaa'_1)=f(\aaa'_2)$. Note that the left hand side of the above equation involves a ratio $p_{(1,\aaa')} /  p_{(0,\aaa')}$ depending on $\aaa'$. %If $p_{(1,\aaa'_1)} /  p_{(0,\aaa'_2)} \neq p_{(1,\aaa'_1)} /  p_{(0,\aaa'_2)}$, 
Equality $f(\aaa'_1)=f(\aaa'_2)$ along with \eqref{eq-s1s2} imply
\[
(c_2 - \bar c_2)\cdot  \frac{p_{(1,\aaa'_1)}}{p_{(0,\aaa'_1)}}
=(c_2 - \bar c_2)\cdot  \frac{p_{(1,\aaa'_2)}}{p_{(0,\aaa'_2)}}\]
%&= (\theta_{2,(1,\aaa'_1)} - \bar\theta_{2,(1,\aaa'_1)})\cdot  \frac{p_{(1,\aaa'_2)}}{p_{(0,\aaa'_2)}} \\
\[
(c_2 - \bar c_2) \cdot \left(\frac{p_{(1,\aaa'_1)}}{p_{(0,\aaa'_1)}} - \frac{p_{(1,\aaa'_2)}}{p_{(0,\aaa'_2)}}\right) = 0
\]
then since $p_{(1,\aaa'_1)} p_{(0,\aaa'_2)} \neq p_{(1,\aaa'_2)}  p_{(0,\aaa'_1)}$ by assumption \eqref{eq-s1s2}, one must have $c_2 = \bar c_2$.
By symmetry of the four item parameters $ g_1$, $ c_1$, $ g_2$ and $ c_2$ in \eqref{eq-allv}, equalities \eqref{eq-cg12} hold as claimed following similar arguments. Now that all the item parameters are identified, $\pp = \bar \pp$. This completes the proof of part (b.1) and part (c) of the theorem.

\noindent
The proof of Theorem \ref{thm-dina-gen} is now complete.

\section{Proof of Theorem 3}\label{sec-proof-gen1}
When Condition $C$ fails and some attribute is required by less than three items, there are two possible scenarios: some attribute is required by only one item, or only two items. We consider them separately, and in both cases prove that $(Q,\TT,\pp)$ are not generically identifiable.
\begin{itemize}
\item[(a)] If some attribute is required by only one item. Then $Q$ must take the following form in \eqref{eq-Q1v} up to column and row permutations, where $\vv_1$ is a binary  vector  of length $K-1$.  %the same conclusion as that in part (a) holds and $(Q,\TT,\pp)$ are not generically identifiable.
\begin{equation}\label{eq-Q1v}
Q = \begin{pmatrix}
1 & \vv_1^\top \\
\mathbf0 & Q^\star
\end{pmatrix};\quad\quad
\bar Q = \begin{pmatrix}
1 & \mathbf 1^\top \\
\mathbf0 & Q^\star
\end{pmatrix}.
\end{equation}

Now for arbitrary model parameters $(\TT,\pp)$ associated with $Q$, 
we also construct $(\bar\TT,\bar\pp)$ associated with the $\bar Q$ in \eqref{eq-Q1v}, such that \eqref{eq1} holds.
Firstly, for any item $j\geq 2$, set $\bar\theta_{j,\aaa} = \theta_{j,\aaa}$ for all $\aaa\in\{0,1\}^K$, then following a similar argument as in Step 3 of the proof of Theorem \ref{thm-dina-gen} (b.1) and (c), we have that \eqref{eq1} hold as long as the following constraints are satisfied:
for any $\aaa'\in\{0,1\}^{K-1}$,
\begin{align}\label{eq-allv1}
\begin{cases}
 p_{(0,\aaa')}+  p_{(1,\aaa')} = \bar p_{(0,\aaa')}+\bar  p_{(1,\aaa')};\\
\theta_{1,\,(0,\aaa')}\cdot p_{(0,\aaa')} + \theta_{1,\,(1,\aaa')}\cdot p_{(1,\aaa')} = \bar\theta_{1,\,(0,\aaa')}\cdot\bar p_{(0,\aaa')} + \bar \theta_{1,\,(1,\aaa')}\cdot\bar p_{(1,\aaa')}.\\
\end{cases}
\end{align}
For each $\aaa'\in\{0,1\}^{K-1}$, we now still arbitrarily set the value of $\bar\theta_{1,\,(0,\aaa')}$ and $\bar\theta_{1,\,(1,\aaa')}$, and set the proportions parameters to be
\begin{align*}
\bar p_{(1,\aaa')}&=
\frac{(\theta_{1,\,(0,\aaa')} - \bar\theta_{1,\,(0,\aaa')})p_{(0,\aaa')} + (\theta_{1,\,(1,\aaa')} - \bar\theta_{1,\,(0,\aaa')})p_{(1,\aaa')} }{\bar\theta_{1,\,(1,\aaa')} - \bar\theta_{1,\,(0,\aaa')}} \\
\bar p_{(0,\aaa')} &= p_{(0,\aaa')} + p_{(1,\aaa')} - \bar p_{(1,\aaa')},
\end{align*}
for each $\aaa'\in\{0,1\}^{K-1}$. Then \eqref{eq-allv1} holds and further \eqref{eq1} holds.  Since the choice of the $2^{K}$ item parameters $\{\theta_{1,\,(0,\aaa')},\,\theta_{1,\,(1,\aaa')}:\allowbreak\aaa'\in\{0,1\}^{K-1}\}$ are arbitrary,
the original $Q$ and associated parameters are not generically identifiable.

\item[(b)] If some attribute   is required by only two items, then $Q$ takes the form in \eqref{eq-Q2v} up to column/row permutations, where $\vv_1$ and $\vv_2$ are  vectors of length $K-1$ and  $Q^\star$ is a submatrix of size $(J-2)\times(K-1)$.
\begin{equation}\label{eq-Q2v}
Q = \begin{pmatrix}
1 & \vv_1^\top \\
1 & \vv_2^\top \\
\mathbf0 & Q^\star
\end{pmatrix};\quad\quad
\bar Q = \begin{pmatrix}
1 & \mathbf 1^\top \\
1 & \mathbf 1^\top \\
\mathbf0 & Q^\star
\end{pmatrix},
\end{equation}
%Then  for any   model parameters $(\TT,\pp)$ associated with $Q$, there exist infinitely many   sets of   parameters $(\bar \TT,\bar \pp)\neq (\TT,\pp)$ associated with  the $\bar Q$ in \eqref{eq-Q} such that \eqref{eq-orig} holds, making $(Q,\TT,\pp)$ not generically identifiable.

%Consider the true $Q$ and another potentially different $\bar Q$ in the following forms,
%\begin{equation}\label{eq-Q}
%Q = \begin{pmatrix}
%1 & \vv_1^\top \\
%1 & \vv_2^\top \\
%\mathbf0 & Q^\star
%\end{pmatrix},\qquad
%\bar Q = \begin{pmatrix}
%1 & \mathbf1^\top \\
%1 & \mathbf1^\top \\
%\mathbf0 & Q^\star
%\end{pmatrix}.
%\end{equation}
Then for arbitrary model parameters $(\TT,\pp)$ associated with $Q$, 
%$$
%\bar Q = \begin{pmatrix}
%1 & \mathbf1^\top \\
%1 & \mathbf1^\top \\
%\mathbf0 & Q^\star
%\end{pmatrix}
%$$
we next carefully construct $(\bar\TT,\bar\pp)$ associated with the $\bar Q$ in \eqref{eq-Q2v}, such that \eqref{eq1} holds. This would prove the conclusion that joint generic identifiability fails.
%\begin{equation}\label{eq-tp}
%T( Q,\TT)\pp = T(\bar Q,\bar\TT)\bar\pp.
%\end{equation}
Firstly, for any item $j\geq 3$, set $\bar\theta_{j,\aaa} = \theta_{j,\aaa}$ for all $\aaa\in\{0,1\}^K$, then following the same argument as in Step 3 of the proof of Theorem \ref{thm-dina-gen} (b.1) and (c), we have that \eqref{eq1} hold as long as the following constraints are satisfied
%\begin{eqnarray}\label{eq-nf}
%&&\forall (r_1,r_2)\in\{0,1\}^2,~\forall \aaa'\in\{0,1\}^{K-1},\\
%&&\mathbb  P(R_1\geq r_1,\,R_2\geq r_2,\, \ba_{2:K}=\aaa'\mid  Q, \TT, \pp)= { \mathbb P}(R_1\geq r_1,\,R_2\geq r_2,\, \ba_{2:K}=\aaa'\mid \bar Q,\bar\TT,\bar\pp),\notag
%\end{eqnarray}
%and this yields that 
for every $\aaa'\in\{0,1\}^{K-1}$,
%\begin{equation}
%\begin{aligned}
%&\theta_{1,(0,\aaa')}^{r_1}\, \theta_{2,(0,\aaa')}^{r_2} \,p_{(0,\aaa')} + \theta_{1,(1,\aaa')}^{r_1}\, \theta_{2,(1,\aaa')}^{r_2}\, p_{(1,\aaa')}\\
%=&
%\bar\theta_{1,(0,\aaa')}^{r_1}\, \bar\theta_{2,(0,\aaa')}^{r_2} \,\bar p_{(0,\aaa')} + \bar\theta_{1,(1,\aaa')}^{r_1} \,\bar\theta_{2,(1,\aaa')}^{r_2}\,\bar p_{(1,\aaa')}
%\end{aligned}
%\end{equation}
\begin{align}\label{eq-allv}
\begin{cases}
 p_{(0,\aaa')}+  p_{(1,\aaa')} = \bar p_{(0,\aaa')}+\bar  p_{(1,\aaa')};\\
\theta_{1,\,(0,\aaa')}\cdot p_{(0,\aaa')} + \theta_{1,\,(1,\aaa')}\cdot p_{(1,\aaa')} = \bar\theta_{1,\,(0,\aaa')}\cdot\bar p_{(0,\aaa')} + \bar \theta_{1,\,(1,\aaa')}\cdot\bar p_{(1,\aaa')};\\
\theta_{2,\,(0,\aaa')}\cdot p_{(0,\aaa')} + \theta_{2,\,(1,\aaa')}\cdot p_{(1,\aaa')} = \bar\theta_{2,\,(0,\aaa')}\cdot\bar p_{(0,\aaa')} + \bar \theta_{2,\,(1,\aaa')}\cdot\bar p_{(1,\aaa')};\\
\theta_{1,\,(0,\aaa')}\theta_{2,\,(0,\aaa')}\cdot p_{(0,\aaa')} + \theta_{1,\,(1,\aaa')}\theta_{2,\,(1,\aaa')}\cdot p_{(1,\aaa')} 
\\
\qquad\qquad\qquad = \bar\theta_{1,\,(0,\aaa')}\bar\theta_{2,\,(0,\aaa')}\cdot\bar p_{(0,\aaa')} + \bar \theta_{1,\,(1,\aaa')}\bar\theta_{2,\,(1,\aaa')}\cdot\bar p_{(1,\aaa')}.
\end{cases}
\end{align}
%Since the structure of $\bar Q$ allows the $2\times 2^{K-1}$ item parameters $\{\bar \theta_{(0,\aaa')},~\bar \theta_{(1,\aaa')}:~\aaa'\in\{0,1\}^{K-1}\}$ to vary freely, the above system of equations have infinitely many sets of different solutions and the conclusion follows. 
For each $\aaa'\in\{0,1\}^{K-1}$, arbitrarily choose $\bar\theta_{1,\,(0,\aaa')}$ and $\bar\theta_{2,\,(0,\aaa')}$ from the neighborhood of the true parameter values $\theta_{1,(0,\aaa')}$ and $\theta_{2,(1,\aaa')}$ respectively. Then set
\begin{align}\label{eq-sol}
\begin{cases}
\bar \theta_{1,\,(1,\aaa')} = \theta_{1,(0,\aaa')} + \frac{([\theta_{1,\,(1,\aaa')} - \theta_{1,\,(0,\aaa')}][\theta_{2,\,(1,\aaa')} - \bar\theta_{2,\,(0,\aaa')}]p_{(1,\aaa')}}{[\theta_{2,\,(0,\aaa')} - \bar\theta_{2,\,(0,\aaa')}]p_{(0,\aaa')} + [\theta_{2,\,(1,\aaa')} - \bar\theta_{2,\,(0,\aaa')}] p_{(1,\aaa')}},\\ 
%%%
\bar \theta_{2,\,(1,\aaa')} = \theta_{2,(0,\aaa')} + \frac{[\theta_{2,\,(1,\aaa')} - \theta_{2,\,(0,\aaa')}][\theta_{1,\,(1,\aaa')} - \bar\theta_{1,\,(0,\aaa')}]p_{(1,\aaa')}}{[\theta_{1,\,(0,\aaa')} - \bar\theta_{1,\,(0,\aaa')}]p_{(0,\aaa')} + [\theta_{1,\,(1,\aaa')} - \bar\theta_{1,\,(0,\aaa')}] p_{(1,\aaa')}},\\ 
%%%
\bar p_{(1,\aaa')} = \frac{[\theta_{2,\,(0,\aaa')} - \bar\theta_{2,\,(0,\aaa')}]p_{(0,\aaa')} + [\theta_{2,\,(1,\aaa')} - \bar\theta_{2,\,(0,\aaa')}]p_{(1,\aaa')}}{\bar\theta_{2,\,(1,\aaa')} - \bar\theta_{2,\,(0,\aaa')}}, \\ 
%%%
\bar p_{(0,\aaa')} = p_{(0,\aaa')} + p_{(1,\aaa')} - \bar p_{(1,\aaa')}.
\end{cases}
\end{align}
Then  one can check that \eqref{eq-allv} holds and further \eqref{eq1} holds. Since in the above construction the choice of the $2^{K}$ item parameters $\{\theta_{1,\,(0,\aaa')},\, \allowbreak \theta_{2,\,(0,\aaa')}:\allowbreak\aaa'\in\{0,1\}^{K-1}\}$ are arbitrary, we have proved that the $Q$ and associated model parameters are not generically identifiable.
%\end{proof}

\end{itemize}

\section{Proof of Theorem 4}
We prove this theorem following a similar argument as the proof of Theorem 7 in \cite{partial}.
Assume $Q$ takes the form $Q=(Q_1^\top,Q_2^\top,(Q^\star)^\top)^\top$,
%\[
%Q = 
%\left(\begin{array}{c} 
%Q_1 \\
%Q_2 \\
%Q^\star
%\end{array} \right),
%\]
where $Q_1$ and $Q_2$ have all diagonal elements being 1. %To ensure generic identifiability for the big family of models, here we consider the saturated model where all the main effect and interaction effect terms are included in modeling the item parameters, namely 
Assume
\[\begin{aligned}
\theta_{j,\aaa} 
=& f\Big(\beta_{j,0} + \sum_{k=1}^K\beta_{j,k}q_{j,k}\alpha_k  + \sum_{k'=k+1}^K\sum_{k=1}^{K-1}\beta_{j,kk'}(q_{j,k}\alpha_k)(q_{j,k'}\alpha_{k'}) + \cdots + \beta_{j,12\cdots K}\prod_k (q_{j,k}\alpha_k)\Big),
\end{aligned}\]
where $f(\Cdot)$ is some link function and when $f(\Cdot)$ is the identify function, the model is the GDINA model.
We first show that under Condition $D$, the $2^K\times 2^K$ matrices $T(Q_1,\TT_{Q_1})$ and $T(Q_2,\TT_{Q_2})$ both have full rank $2^K$ generically. %To show the generic identifiability, 
It suffices to find some valid $\TT$ (i.e., $\TT_Q$) that gives 
\begin{equation}\label{eq-det}
\det(T(Q_1,\TT_{Q_1}))\neq 0,\quad \det(T(Q_2,\TT_{Q_2}))\neq 0.
\end{equation}
The reason is as follows. \eqref{eq-det} would imply the polynomials defining the two matrix determinants are not zero polynomials in the $Q$-restricted parameter space. Therefore for almost all parameters, $T(Q_1,\TT_{Q_1})$ and $T(Q_2,\TT_{Q_2})$ would have full rank.
Next we only focus on $T(Q_1,\TT_{Q_1})$. For every item $k=1,\ldots,K$, we set $\beta_{k,k}=1$, $\beta_{k,k'}=0$ for any $k'\neq k$, and set all the interaction effects to zero. Then  $T(Q_1,\TT_{Q_1})$ becomes identical to $T({I_K},\widehat \TT_{I_K})$ under a $Q$-matrix $I_K$ with associated item parameters $\widehat \TT_{I_K}$ defined as follows:
$\hat\theta_{\ee_k,\mz} = \beta_{k,0}$, and  $\hat\theta_{\ee_k,\ee_k} = \hat\theta_{\ee_k,\mo} = \beta_{k,0}+\beta_{k,k}$ for $k\in\{1,\ldots,K\}$.
It is not hard to see that $T({I_K},\widehat \TT_{I_K})$ can be viewed as a $T$-matrix under the DINA model with the $Q$-matrix equal to $I_K$, and guessing parameters $\beta_{k,0}$, slipping parameters $1-\beta_{k,0}-\beta_{k,k}$ for $k=1,\ldots, K$. 
Therefore $T({I_K},\widehat \TT_{ I_K})$ has full rank as argued in Step 1 of the proof of Theorem \ref{thm-dina}.
So $T(Q_1,\TT_{Q_1})$ has full rank generically.
%%%%%%%%%
%Applying the result that $T_{I_K}$ is of full rank $2^K$, we have $T_{Q_1}(\TT)$ is generically of full rank.
%
%We next show that if Condition B additionally holds, then any two different columns indexed by attribute profiles $\aaa$ and $\aaa'$ of $T_{Q_1}(\TT_{(2K+1):J})$ are generically distinct. For distinct $\aaa$, $\aaa'\in\{0,1\}^K$, they at least differ in one attribute $k$. Without loss of generality, assume $\alpha_k=1>0=\alpha'_k$. Condition B ensures that there exists some item $j>2K$ such that $q_{j,k}=1$.  Under the saturated model considered here, this implies $\theta_{j,\aaa}\neq\theta_{j,\aaa'}$.  
%%%%%%%%%%%%%%%%%%%%%%%%%%
%This further implies that $T(Q_{(2K+1):J}$, $\TT_{(2K+1):J})$ is generically of Kruskal rank 2. Then 
%\[\begin{aligned}
%&\text{rank}_{K} \{T(Q_1, \TT_{Q_1})\} + \text{rank}_{K} \{T(Q_2, \TT_{Q_2})\} + \text{rank}_K \{T(Q_{(2K+1):J}, \TT_{(2K+1):J})\} \\
%=& 2 \times 2^K + 2.
%\end{aligned}\]
%Then apply Corollary 2 of \cite{rhodes2010concise} to this $2^K$-class latent class model to get $T(Q, \TT) = T(Q, \bar\TT)$ and $\pp = \bar\pp$. This proves the generic identifiability of all the parameters in the model.

We next prove that if Condition $E$ holds in addition, then any two different columns of $T(Q^\star,\TT_{Q^\star})$ are  distinct generically. For $\aaa$, $\aaa'\in\{0,1\}^K$ and $\aaa\neq\aaa'$, they at least differ in one element.  Assume without loss of generality that $\alpha_k=1>0=\alpha'_k$. Then Condition $E$ ensures that there is some item $j>2K$ with $q_{j,k}=1$.  Under the general RLCM, this implies $\theta_{j,\aaa}\neq\theta_{j,\aaa'}$ generically.  
%%%%%%%%%%%%%%%%%%%%%%%%%
By \cite{kruskal1977three}, a matrix's Kruskal rank is the largest number $I$ such that every set of $I$ columns of the matrix are independent. When a matrix has full rank,  its Kruskal rank equals its rank.
By this definition, $T(Q^\star$, $\TT_{Q^\star})$ has Kruskal rank at least 2 generically, and  $T(Q_1, \TT_{Q_1})$, $T(Q_2, \TT_{Q_2})$ have Kruskal rank $2^K$ generically. Then for generic $\TT_Q$, we have
\begin{align}\label{eq-krank}
\text{rank}_{K} \{T(Q_1, \TT_{Q_1})\} + \text{rank}_{K} \{T(Q_2, \TT_{Q_2})\} + \text{rank}_K \{T(Q^\star, \TT_{Q^\star})\}\geq 2 \times 2^K + 2.
\end{align}
Applying Corollary 2 of \cite{rhodes2010concise} to this $2^K$-class latent class model, we get $T(Q, \TT) = T(Q, \bar\TT)$ and $\pp = \bar\pp$ up to column permutation. This proves  generic identifiability of $(Q,\TT,\pp)$ in the model.
Moreover, we also have the following form of the identifiable set
\begin{align*}
&\boldsymbol\vartheta_{Q}\setminus \boldsymbol\vartheta_{non}
= \{
 (\TT_Q,\pp):   \det (T(Q_1,\TT_{Q_1}))\neq 0, \det (T(Q_2,\TT_{Q_2}))\neq 0,\\
 &  T(Q^\star,\TT_{Q^\star})\cdot \text{Diag}(\pp)~\text{has column vectors different from each other}\}.
\end{align*}
This is because when $(\TT_Q,\pp)\in \boldsymbol\vartheta_{Q}\setminus \boldsymbol\vartheta_{non}$, the rank condition \eqref{eq-krank} is satisfied and joint identifiability of $(Q,\TT_Q,\pp)$ follows.  
 
\noindent
\section{Proof of Theorem 5.}
%We will show that if $Q$ is not generically complete, then $(Q,\TT,\pp)$ are not be generically identifiable.
%\begin{itemize}
%\item[$H$.] after any row and/or column permutation, $Q$ can not be transformed to contain a $K\times K$ submatrix whose diagonal entries are all ``1"s.
%\end{itemize} 
We prove the theorem in two steps.
In the first step, we show that if $Q$ is not generically complete, than it must take the following form (up to  column/row permutations) for some $k>m$,

\begin{align}\label{eq-bigq}
Q = \left(\begin{array}{ccc|ccc}
q_{1,1}  & \cdots & q_{1,k}      & *      & \cdots & *     \\
\vdots & \vdots & \vdots & \vdots & \vdots & \vdots\\
q_{m,1}  & \cdots & q_{m,k}      & *      & \cdots & *     \\
\hline
0      & \cdots & 0      & *      & \cdots & *     \\
\vdots & \vdots & \vdots & \vdots & \vdots & \vdots\\
0      & \cdots & 0      & *      & \cdots & *     \\
\end{array}\right)
=\left(\begin{array}{c|c}
Q_{11} & Q_{12} \\
\hline
Q_{21} & Q_{22} \\
\end{array}\right)
=\left(\begin{array}{c}
Q_1\\
Q_2
\end{array}\right).
\end{align}
The bottom-left submatrix $Q_{21}=\mathbf0_{(J-m)\times k}$.  Any entry not in $Q_{21}$ can be either 0 or 1. 
We introduce some definitions first. Given a $Q$-matrix $Q$, define a family $S_Q$ of $K$ finite sets $S_Q=\{\ma_1,\ma_1,\ldots,\ma_K\}$, where $\ma_k = \{1\leq j\leq J:~ q_{j,k}=1\}$ for each $k$. Then $\ma_k$ denotes the set of items that require attribute $k$. For the family $S_Q$, a \textit{transversal} is a system of distinct representatives from each of its elements $\ma_1,\ldots,\ma_K$. For example, for 
$$
Q = \begin{pmatrix}
1 & 1 & 0\\
0 & 1 & 1\\
1 & 0 & 1\\
\end{pmatrix},
$$
we have $S_Q=\{\ma_1=\{1,3\}$, $\ma_2=\{1,2\}$, $\ma_3=\{2,3\}\}$. Then $(1,2,3)$ is a valid transversal of $S_Q$, and so as $(3,1,2)$; but $(1,1,2)$ is not a transversal.
Now it is not hard to see that, the assumption that $Q$ is not generically complete is equivalent to the following statement $H^\star$,
\begin{itemize}
\item[$H^\star$.] Given $Q$, the family $S_Q$ does not have a valid transversal. 
\end{itemize}
Then by Hall's Marriage Theorem \citep{hall}, the nonexistence of a transversal indicates the failure of the marriage condition. So there must exist a subfamily $W\subseteq S_Q$ such that
$|W| > |\bigcup_{\ma\in W} \ma|.$ 
More specifically, this means there exist some $l_1,l_2,\ldots,l_k\in\{1,\ldots,K\}$ and $W=\{\ma_{l_1},\ldots,\ma_{l_k}\}$ such that 
$$|W|=k > |\ma_{l_1}\cup\cdots\cup\ma_{l_k}|\stackrel{\text{def}}{=}m.$$
In other words, we have shown that there exist some attributes, the number of which (e.g., $k$) exceeds the number of items that require any of these attributes (e.g., $m$). This is exactly saying that $Q$ has to take the form of \eqref{eq-bigq} with $k>m$ after some column/row permutation.

In the second step, we show that %if condition $H$ is true, then 
if $Q$ takes the form of \eqref{eq-bigq}  with $k>m$, then $(Q,\TT,\pp)$ under general RLCMs are not generically identifiable. 
Now we define another potentially different $\bar Q$ as 
$$
\bar Q =\left(\begin{array}{c|c}
Q_{11} & \bar Q_{12} \\
\hline
Q_{21} & Q_{22} \\
\end{array}\right)
=\left(\begin{array}{c}
\bar Q_1\\
Q_2
\end{array}\right), 
\quad\text{where}~ \bar Q_{12} = \mathbf 1_{m\times(K-k)}.
$$
Then given arbitrary $(\TT,\pp)$ associated with $Q$, we set $\bar \theta_{j,\aaa}=\theta_{j,\aaa}$ for every $j=m+1,\ldots,J$ and  every $\aaa\in\{0,1\}^K$. Because $Q_{21}$ is a $(J-m)\times k$ zero matrix, we claim that under the current construction, the original $2^J$ constraints in \eqref{eq1} are satisfied as long as the following constraints are satisfied %for any $\aaa'=(\alpha_{k+1},\ldots,\alpha_{K})\in\{0,1\}^{K-k}$, and any $\rr'=(r_1,\ldots,r_{m})\in \{0,1\}^{m}$,
\begin{align*}
&\qquad\forall \aaa'=(\alpha_{k+1},\ldots,\alpha_{K})\in\{0,1\}^{K-k},\quad 
\forall \rr'=(r_1,\ldots,r_{m})\in \{0,1\}^{m},\\
&\sum_{\aaa^\star\in\{0,1\}^{k}} T_{\rr',\,(\aaa^\star,\aaa')}(Q_1,\TT_{Q_1}) \cdot p_{(\aaa^\star,\aaa')} = 
\sum_{\aaa^\star\in\{0,1\}^{k}} T_{\rr',\,(\aaa^\star,\aaa')}(\bar Q_1,\bar \TT_{\bar Q_1})\cdot \bar p_{(\aaa^\star,\aaa')}.
\end{align*}
This claim can be shown following a similar argument as that in Step 3 of the proof of Theorem \ref{thm-dina-gen} (b.1) and (c).
Then the above system of equations contain $2^{K-k}\times 2^m$ constraints, while under the general RLCMs the number of free variables in $(\bar \TT, \bar\pp)$ involved is 
\begin{align*}
&~\left\lvert \{\bar p_{\aaa}:\aaa\in\{0,1\}^K\}\bigcup
\{ \bar \theta_{j,\aaa}: j\in\{1,\ldots,m\}, \aaa\in\{0,1\}^K \} \right\rvert\\
=&~2^K + 2^{K-k}\times \Big(\sum_{j=1}^m 2^{ q_{j,1} + \cdots + q_{j,k} } \Big) 
\geq 2^K + 2^{K-k}\times m.
\end{align*}
Under the assumption $m<k$, we have that the number of constraints $2^{K-k}\times 2^m$ is smaller than the number of variables to solve (which is lower bounded by $2^{K-k}\times(2^k+m)$), because $2^m<2^k+m$. So there exist infinitely many different sets of solutions of $(\bar\TT,\bar\pp)$ associated with $\bar Q$ such that $T(Q,\TT)\pp = T(\bar Q,\bar\TT)\bar\pp$. Therefore $(Q,\TT,\pp)$ are not generically identifiable and the proof of the theorem is complete.

\bigskip
\noindent 
\section{Proof of Proposition 4}
%\textbf{Proof of Proposition \ref{prop-finite}.}
We show the conclusion following a similar argument as the proof of Proposition 1 in \cite{xu2018jasa}. To establish the bound \eqref{eq-bound} in the proposition, we check the technical conditions in Theorem 1 in \cite{shen2012likelihood}. We first define some notations. For a family of probability mass functions $\mathcal F$, define $H(\Cdot,\mathcal F)$ to be the bracketing Hellinger metric entropy of $\mathcal F$. We call a finite set of function pairs $S(\epsilon,n)=\{(f^l_1,f^u_1),\ldots, (f^l_n,f^u_n)\}$ a \textit{Hellinger $\epsilon$-bracketing} of $\mathcal F$ if the $L_2$ norm $\norm{\sqrt{f^l_i}-\sqrt{f^u_i}}\leq \epsilon$ for all $i=1,\ldots,n$; and further fur any $f\in\mathcal F$, there is an $i$ such that $f^l_i\leq f\leq f^u_i$. The \textit{bracketing Hellinger metric entropy} is defined to be the logarithm of the cardinality of the $\epsilon$-bracketing with the smallest size, namely $H(\Cdot,\mathcal F)=\log\min\{n:\,S(\epsilon,n)\}$.
We next argue that the size of the parameter space of $(\TT,\pp)$ is well controlled under the Hellinger metric. Recall $S$ is defined in the main text before Proposition \ref{prop-finite}, and we define $\mathcal B_S=\mathcal F_S\cap \{h(\eeta,\eeta^0)\leq 2\epsilon\}$ as the local parameter space with $\eeta=(\bo B,\pp)$ denoting general model parameters and $\eeta^0=(\bo B^0,\pp^0)$ denoting the true model parameters.
 According to the argument in the proof of Proposition 1 in \cite{xu2018jasa}, in the considered scenario with fixed $J$ and $K$, for any $\epsilon<1$ and any $t\in(\epsilon/2^4,\epsilon)$, there is $H(t,\mathcal B_S)\leq c\log(J 2^K) |S|\log(2\epsilon/t)$; indeed, there is $H(t,\mathcal B_S) = O(\log(2\epsilon/t))$ uniformly for any $S$, $\epsilon$ and $t$. 
 
 With this upper bound on the Hellinger bracketing entropy, we can apply Theorem 1 in \cite{shen2012likelihood} to obtain 
 $$
 \mathbb P(
	%\ell(\hat\TT,\hat\pp) > \ell(\TT^0,\pp^0))
	\widehat Q \neq Q^0)
	%\leq \mathbb P(\widehat {\bo B}\neq \bo B^0)
	\leq \mathbb P(\widehat {\bo \eta}\neq \widehat\eeta^0)
	\leq c_2\exp\{-c_1 N C_{\min}(\TT^0,\pp^0) 
	\},
$$
where $C_{\min}(\TT^0,\pp^0):=\inf_{\eeta:\,|S|\leq m,S\neq S_0} h^2(\eeta,\eeta^0)$. The above display is the desired \eqref{eq-bound} in the proposition.
%which is \eqref{eq-bound} in the main text.
	
Next we show that when the proposed sufficient conditions for joint strict identifiability hold, the $C_{\min}(\TT^0,\pp^0)$ in \eqref{eq-bound} is bounded away from zero by some positive constant.
When the proposed conditions for joint strict identifiability (such as Conditions $A$, $B$ and $C$ under DINA model are satisfied), the $(\bo B^0,\pp^0)$ here are strictly identifiable. The consequence is that there exists a constant $\delta>0$ such that $h^2(\eeta,\eeta^0)\geq \delta$, where the $m$ denotes the number of free parameters under the $Q^0$ and the RLCM specification. Therefore,
\begin{align*}
C_{\min}(\TT^0,\pp^0)
\geq &~ 
\inf_{\eeta:\,|S|\leq m,S\neq S_0} \frac{h^2(\eeta,\eeta^0)}{2m}\geq \frac{\delta}{2m} >0,
%\\
%\geq &~d_0 \frac{\log(J 2^K)}{N},
\end{align*}
so $C_{\min}(\TT^0,\pp^0) \geq c_0$ for some positive constant $c_0$ holds.
This proves the conclusion that under the proposed strict identifiability conditions, the finite sample error bound $\mathbb P(\widehat Q \neq Q^0)$ has an exponential rate. This completes the proof of the proposition.
%\qed

%\newpage
\section{Simulation Studies}
In this section, we provide more simulation results to verify the developed identifiability theory. In Section \ref{sec-simu-dina}, we perform simulation studies to verify Theorems 1 and 2 for the DINA model. In Section \ref{sec-simu-gdina}, we perform simulation studies to verify Theorems 3 and 4 for the GDINA model. 
The Matlab code for performing the simulation studies are available at \verb|https://github.com/yuqigu/Identify_Q|.

To better illustrate the identifiability or non-identifiability phenomena of $Q$-matrix, in some of the following simulation studies, we conduct exhaustive search of all possible $Q$-matrices of a certain size $5\times 2$. Specifically, consider the set of all the $5\times 2$ binary $Q$-matrices other than those containing some all-zero row vectors. If treating two $Q$-matrices that are identical up to permuting the two columns as equivalent (because they are indeed equivalent in terms of model identifiability), then there are in total $121$ types of $Q$-matrices. 
We denote such a set of $Q$-matrices by $\text{Exhaus}(Q_{5\times 2})$, and denote its elements by $Q^{1},Q^{2},\ldots,Q^{121}$. 
For example, the first three and the last three $Q$-matrices in $\text{Exhaus}(Q_{5\times 2})$ are

\begin{equation*}
	% \label{eq-exhaus}
	\quad\quad\quad
	Q^1 = 
	\begin{pmatrix}
	 0  &   1 \\ 
     0  &   1 \\ 
     0  &   1 \\ 
     0  &   1 \\ 
     0  &   1 \\		
	\end{pmatrix};
     \quad
     Q^2 = 
	\begin{pmatrix}
	 0  &   1 \\ 
     1  &   0 \\ 
     0  &   1 \\ 
     0  &   1 \\ 
     0  &   1 \\		
	\end{pmatrix};
     \quad
     Q^3 = 
	\begin{pmatrix}
	 0  &   1 \\ 
     1  &   1 \\ 
     0  &   1 \\ 
     0  &   1 \\ 
     0  &   1 \\		
	\end{pmatrix}; \quad \cdots\cdots
\end{equation*}
\begin{equation*}
	% \label{eq-exhaus}
	Q^{119} = 
	\begin{pmatrix}
	 1  &   1 \\ 
     1  &   1 \\ 
     1  &   1 \\ 
     0  &   1 \\ 
     1  &   0 \\		
	\end{pmatrix};
     \quad
     Q^{120} = 
	\begin{pmatrix}
	 1  &   1 \\ 
     1  &   1 \\ 
     1  &   1 \\ 
     0  &   1 \\ 
     1  &   1 \\	
	\end{pmatrix};
     \quad
     Q^{121} = 
	\begin{pmatrix}
	 1  &   1 \\ 
     1  &   1 \\ 
     1  &   1 \\ 
     1  &   1 \\ 
     1  &   1 \\	
	\end{pmatrix}.
\end{equation*}
The complete list of the 121 $Q$-matrices in the set $\text{Exhaus}(Q_{5\times 2})$ is available in the Matlab file \verb|Q_aa.mat| at \verb|https://github.com/yuqigu/Identify_Q|.

In the exhaustive-search scenario, to illustrate the identifiability/non-identifiability phenomenon, we will generate data using some particular $Q$-matrix, and fit the dataset using all the 121 candidate $Q$-matrices in $\text{Exhaus}(Q_{5\times 2})$ and plot the log-likelihood values corresponding to all these 121 $Q$-matrices. Investigating whether the true data-generating $Q$-matrix achieves the maximum of the likelihood would help gain insight into whether this true $Q$-matrix is identifiable in the considered practical setting. We will see from these simulations how the developed identifiability theory matches the practice.

\subsection{Two-Parameter RLCM: DINA Model}\label{sec-simu-dina}
In this section, we carry out four simulation studies.

%\newpage
\noindent
\textbf{Study I: When $Q$-matrix satisfies the necessary and sufficient conditions $A$, $B$ and $C$ for strict identifiability.}

In this simulation study, we choose those $Q$-matrices from $\text{Exhaus}(Q_{5\times 2})$ that satisfies the proposed necessary and sufficient identifiability conditions $A$, $B$ and $C$ in Theorem 1 of the main text.
In particular, after rearranging rows, there are exactly two forms the $5\times 2$ $Q$-matrix that satisfies $A$, $B$ and $C$. Their representatives are $Q^{18}$ and $Q^{15}$ as follows,
\begin{equation*}
	Q^{18}=\begin{pmatrix}
		0 & 1 \\ 
		1 & 1 \\ 
		1 & 1 \\ 
		1 & 0 \\ 
		0 & 1 \\
	\end{pmatrix};\qquad\qquad
		Q^{15}=\begin{pmatrix}
		0 & 1 \\ 
		1 & 1 \\ 
		1 & 0 \\ 
		1 & 0 \\ 
		0 & 1 \\
	\end{pmatrix}.
\end{equation*}
Note that $Q^{18}$ contains only on identity submatrix $I_2$, while $Q^{15}$ contains two copies of submatrix $I_2$. As introduced prior to this section \ref{sec-simu-dina}, we generate datasets with sample size $N=10^5$ with true $Q$-matrix being  $Q^{18}$ and $Q^{15}$, respectively; and for each case, we run EM algorithm with several random initializations to fit the dataset with all the 121 $Q$-matrixes in $\text{Exhaus}(Q_{5\times 2})$ and obtain their log-likelihood values.

Figure \ref{fig-q18} and \ref{fig-q15} present the log-likelihood plots, with $x$-axis denoting the indices of the 121 candidate $Q$-matrices in $\text{Exhaus}(Q_{5\times 2})$, and $y$-axis denoting the log-likelihood values. Each blue triangle denotes a candidate $Q$-matrix; the red star denotes the true data-generating $Q$-matrix, and the purple square denotes the $Q$-matrix that achieves the largest likelihood.

We can see from these two plots in Figure \ref{fig-dina-abc} that when the true data-generating $Q$-matrix ($Q^{15}$ and $Q^{18}$) satisfies our proposed  conditions $A$, $B$ and $C$, it indeed achieves the largest likelihood compared to all other possible candidate $Q$-matrices. 
Therefore for any algorithm seeking the maximum likelihood estimator of $(Q,\cc,\cg,\pp)$, the true $Q$-matrix can be identified and any other $Q$-matrix will not be confused with the true $Q$.
Another observation from Figure \ref{fig-q18} and \ref{fig-q15} is that, for $Q^{15}$ that contains one more identity submatrix $I_2$ than $Q^{18}$, the true $Q$ can be relatively better distinguished from the other $Q$'s due to the larger gap in the likelihood values. This phenomenon might imply that the more identity submatrices the true data-generating $Q$-matrix contain, the easier the estimation for the true structure would be.

\begin{figure}[H]
\caption{DINA: exhaustive search in the set of $5\times 2$ $Q$-matrices with a true $Q$-matrix satisfying Conditions $A$, $B$ and $C$ in Theorem 1.}
\label{fig-dina-abc}

\centering
\begin{subfigure}{0.7\textwidth}
\includegraphics[width=\linewidth]{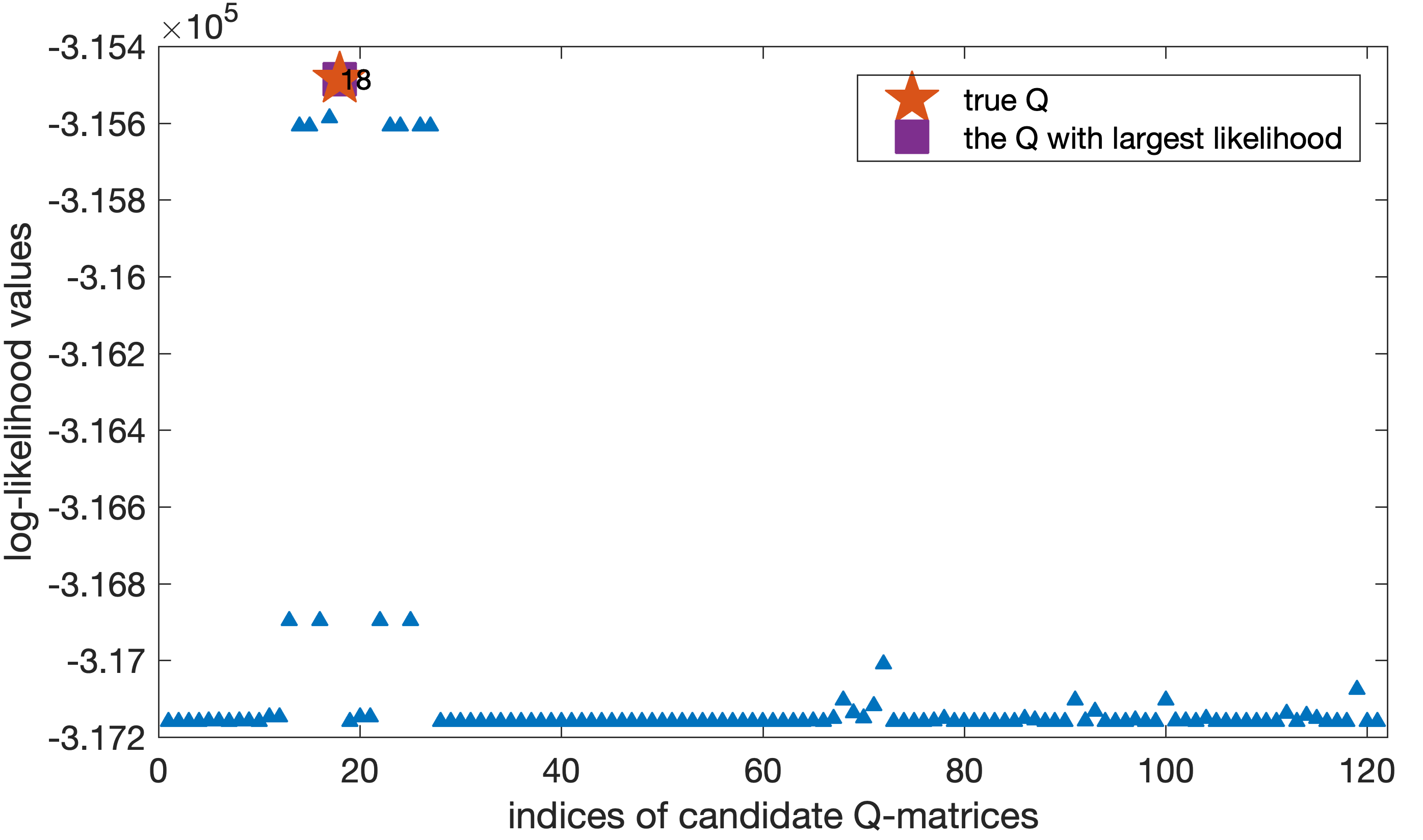}
\caption{true $Q$ containing one $I_2$:
$Q^{18}=\begin{pmatrix}
0  &   1  &   1  &   1  &   0 \\ 
1   &  1   &  1  &   0  &   1
\end{pmatrix}^\top$}
\label{fig-q18}
\end{subfigure}
\end{figure}
\begin{figure}[H]\ContinuedFloat
\centering
\begin{subfigure}{0.7\textwidth}
\includegraphics[width=\linewidth]{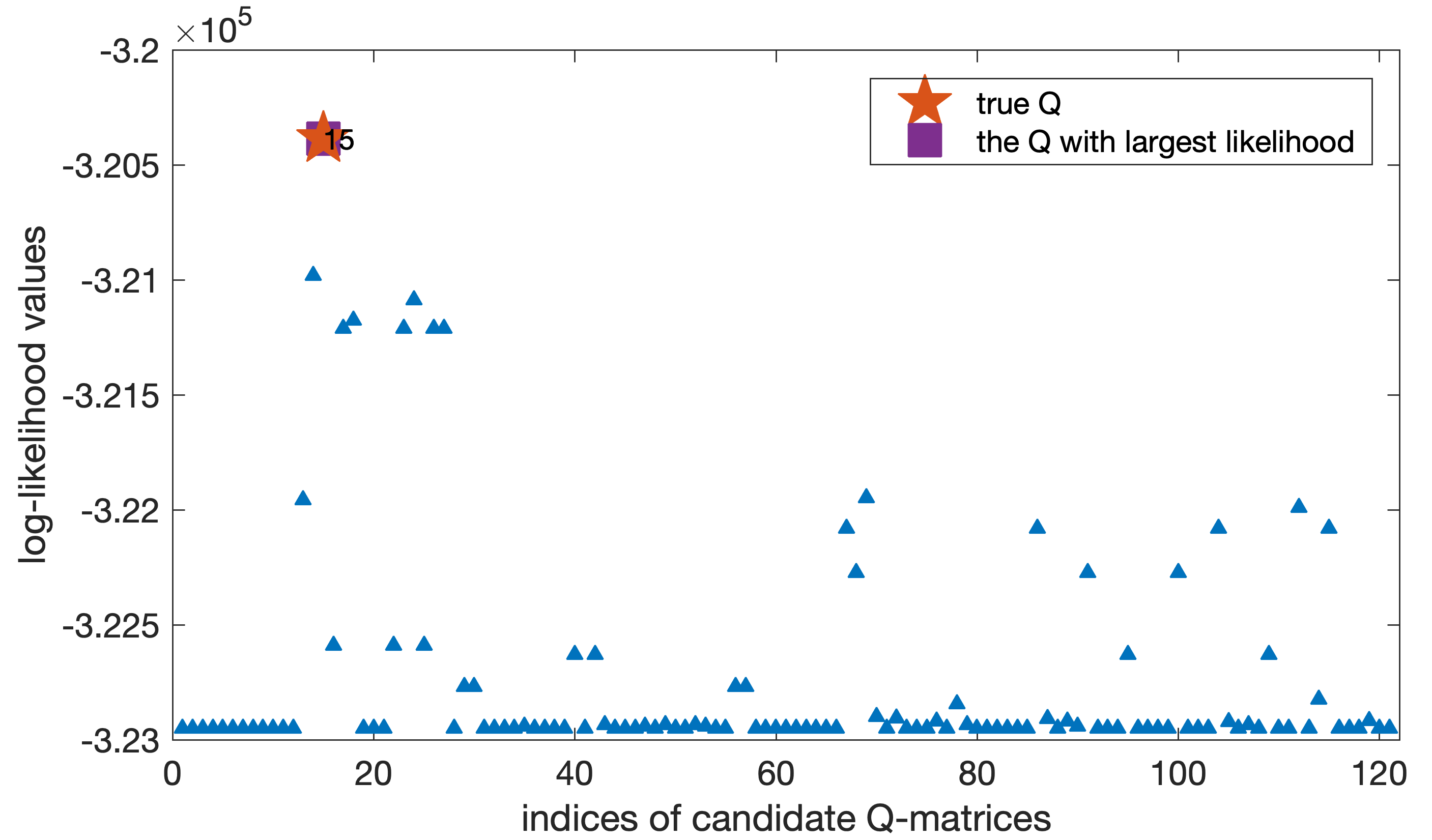}
\caption{true $Q$ containing two $I_2$'s:
$Q^{15}=\begin{pmatrix}
     0  &   1   &  1   &  1  &   0 \\ 
     1  &   1   &  0   &  0  &   1
\end{pmatrix}^\top$} 
\label{fig-q15}
\end{subfigure}

\end{figure}

\noindent
\textbf{Study II: When $Q$-matrix does not satisfy all of Conditions $A$, $B$, $C$ but satisfies conditions in Theorem 2 for generic identifiability.}

\begin{figure}[H]
\centering
\includegraphics[width=0.7\linewidth]
{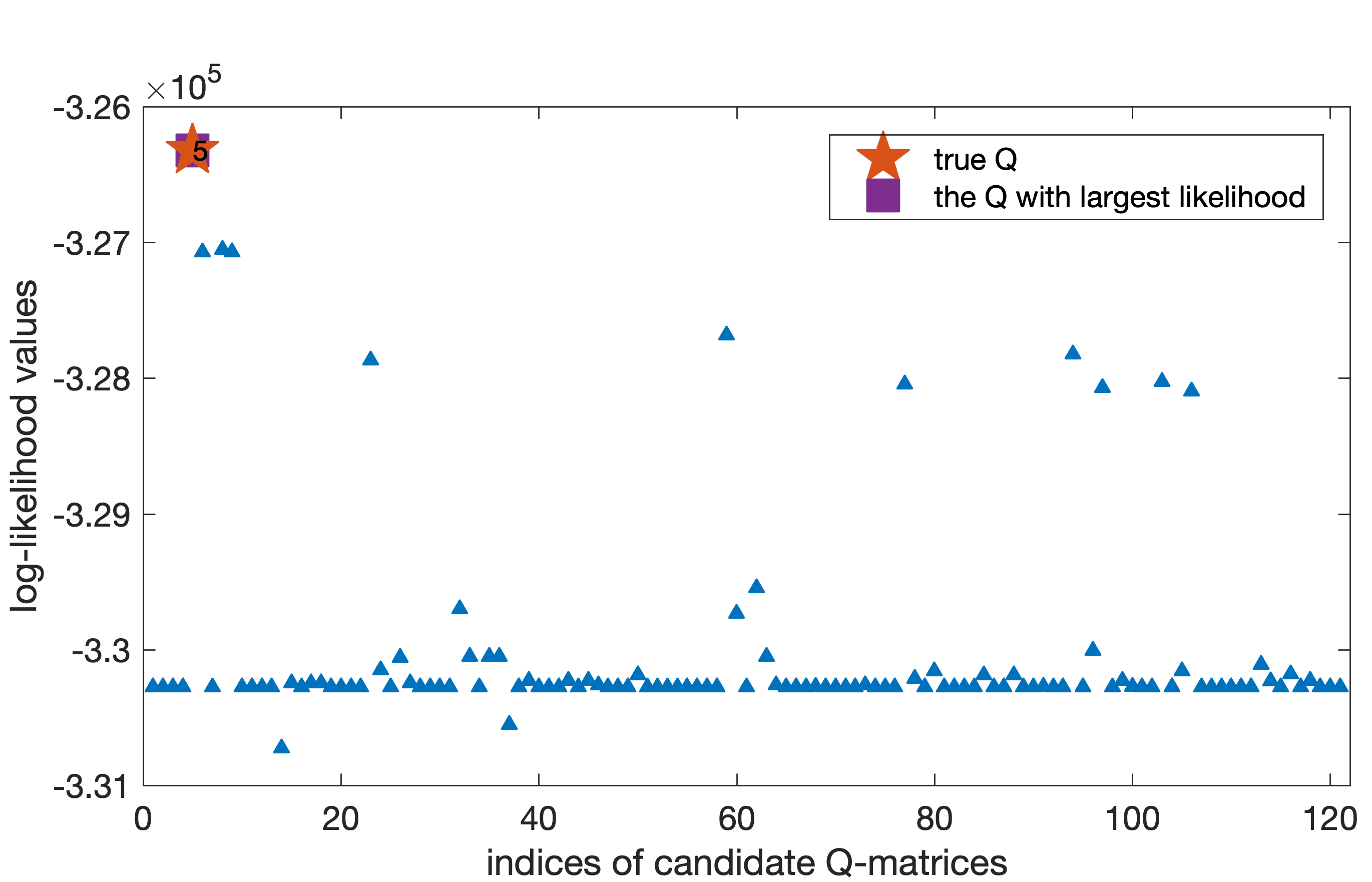}
\caption{DINA: exhaustive search in $5\times 2$ $Q$-matrices with a true $Q$-matrix $Q^5$ in \eqref{eq-q5} generically identifiable, corresponding to scenario (b.2) in Theorem 2.}
\label{fig-dina-gid}
\end{figure}

In this simulation study, we take the data-generating $Q$-matrix from $\text{Exhaus}(Q_{5\times 2})$ that do NOT satisfy some of Conditions $A$, $B$ and $C$, but satisfy the conditions in Theorem 2 for joint generic identifiability of $(Q,\cc,\cg,\pp)$.
In particular, for the considered case of $K=2$, the only possibility for (global) generic identifiability is scenario (b.2) described in Theorem 2, where Condition $C$ is violated and some column of $Q$ contains only two entries of ``1''.
After rearranging the rows of $Q$, it is not hard to see that there is only one possible case of the form of $Q$ leading to generic identifiability, and the following $Q^5$ is a representative,
\begin{equation}\label{eq-q5}
	Q^5=
%	\begin{pmatrix}
%     0  &   1  &   1   &  0  &   0 \\
%     1  &   0  &   0  &   1   &  1
%\end{pmatrix}^\top.
 	\begin{pmatrix}
		0 & 1 \\ 
		1 & 0 \\ 
		1 & 0 \\ 
		0 & 1 \\ 
		0 & 1
	\end{pmatrix}.
\end{equation}
The log-likelihood value plot is presented in Figure \ref{fig-dina-gid}. One can see in this generically identifiable scenario, with randomly generated true parameters, the true $Q$-matrix $Q^5$ achieves the largest likelihood and hence can be identified from data. We point out that although only the result of one simulated dataset is presented here,  the generically identifiable $Q$-matrix (as the true $Q$-matrix) generally can achieve the largest likelihood among all the candidate $Q$-matrices, based on our experience in various simulations.

%\newpage
\noindent
\textbf{Study III: When $Q$-matrix does not even lead to  local identifiability.}

In this simulation study, we take the data-generating $Q$-matrix from $\text{Exhaus}(Q_{5\times 2})$ that \textit{do not even lead to local identifiability}. That is, under such true $Q$-matrix, even in a small neighborhood of the true parameters, there exist infinitely many different alternative parameters that are not distinguishable from the true one.

Consider the following three different forms of $Q$-matrices from the set $\text{Exhaus}(Q_{5\times 2})$,
\begin{equation*}
	\quad\quad\quad
	Q^{10} = 
	\begin{pmatrix}
	 0  &   1 \\ 
     0  &   1 \\ 
     0  &   1 \\ 
     1  &   0 \\ 
     0  &   1 \\	
	\end{pmatrix};
     \quad
     Q^{21} = 
	\begin{pmatrix}
	 0  &   1 \\ 
     1  &   1 \\ 
     0  &   1 \\ 
     1  &   1 \\ 
     0  &   1 \\		
	\end{pmatrix};
     \quad
     Q^{55} = 
	\begin{pmatrix}
	 0  &   1 \\ 
     0  &   1 \\ 
     0  &   1 \\ 
     0  &   1 \\ 
     0  &   1 \\		
	\end{pmatrix},
\end{equation*}
where $Q^{10}$ contains only one entry of ``1'' in one column, $Q^{21}$ is incomplete (i.e., lacks $I_2$), and $Q^{55}$ contains an all-zero column.
Their corresponding log-likelihood plots in the exhaustive-search scenario are presented in Figure \ref{fig-q10}, \ref{fig-q21} and \ref{fig-q55}.
One can see from these plots that in these no even locally identifiable settings, the true data-generating $Q$-matrix does not achieve the maximum of the likelihood. Instead, many other alternative $Q$-matrices would have larger likelihood, and a wrong $Q$-matrix will be selected as the maximum likelihood estimator.

\begin{figure}[H]
\caption{DINA: exhaustive search in $5\times 2$ $Q$-matrices with a true $Q$-matrix \textit{not even locally} identifiable, corresponding to scenario (b.1) in Theorem 2.}
\label{fig-exhaus-dina-nolid}

\centering
\begin{subfigure}[b]{0.7\textwidth}
\includegraphics[width=\linewidth]
{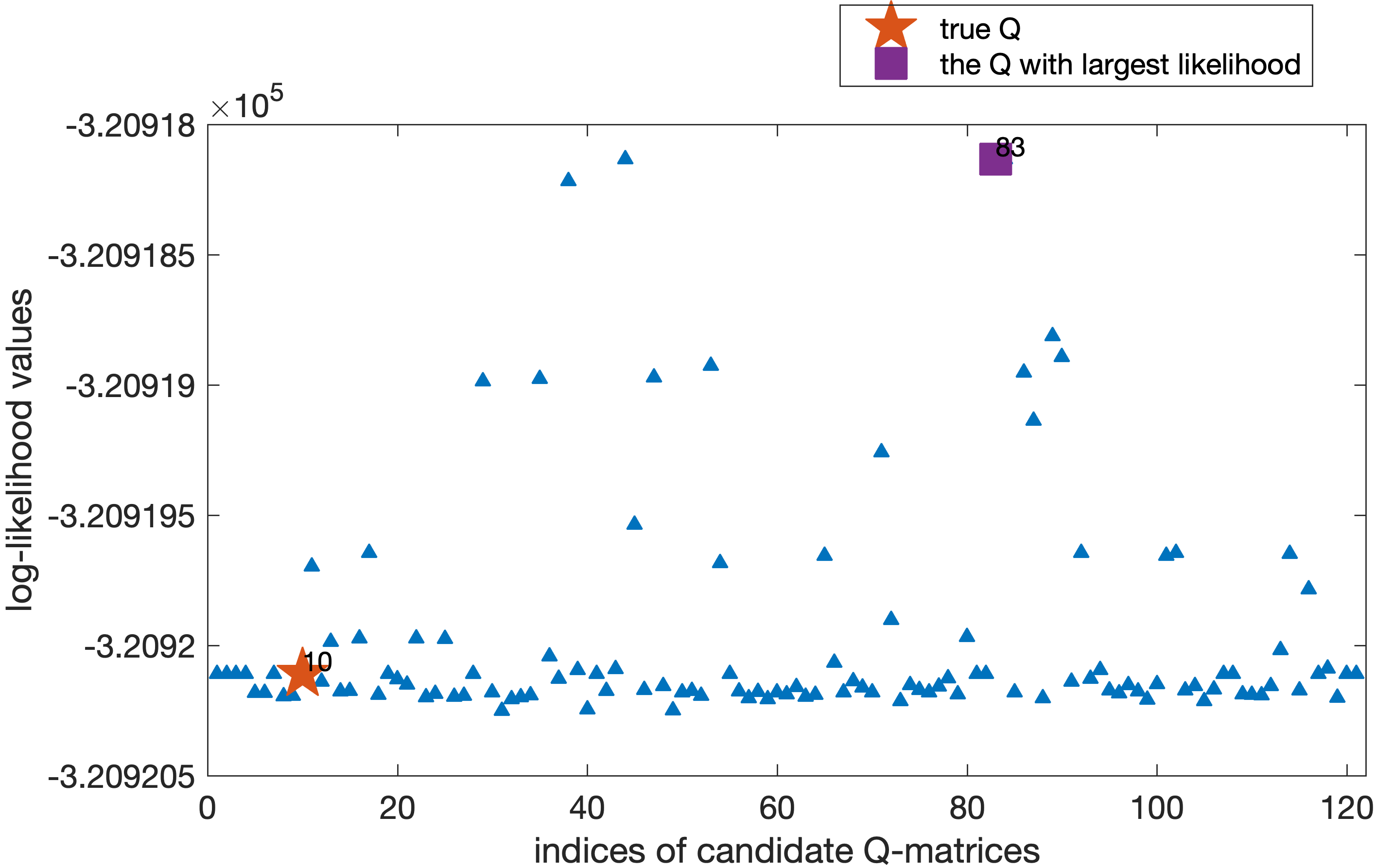}
\caption{true $Q$ not even locally identifiable:
$
Q^{10} = \begin{pmatrix}
     0   &  0   &  0  &   1   &  0 \\ 
     1   &  1   &  1  &   0   &  1
\end{pmatrix}^\top
$} 
\label{fig-q10}
\end{subfigure}
\end{figure}

\begin{figure}[H]\ContinuedFloat
\centering
\begin{subfigure}[b]{0.7\textwidth}
\includegraphics[width=\linewidth]
{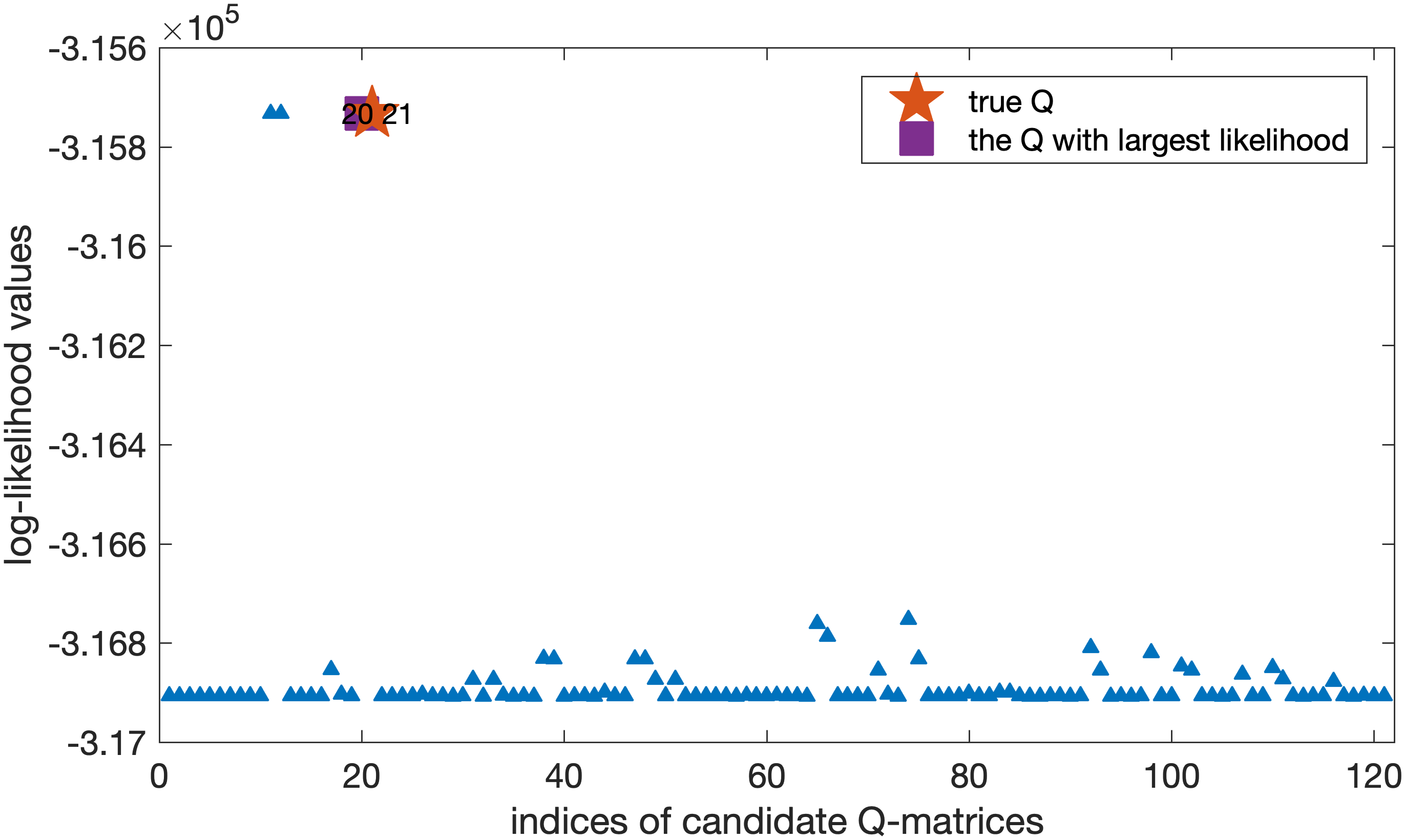}
\caption{true $Q$ not even locally identifiable:
$
Q^{21} = \begin{pmatrix}
     0  &   1  &   0  &   1  &   0 \\ 
     1   &  1   &  1   &  1   &  1
\end{pmatrix}^\top
$}
\label{fig-q21}
\end{subfigure}
\end{figure}

\begin{figure}[H]\ContinuedFloat
\centering
\begin{subfigure}[b]{0.7\textwidth}
\includegraphics[width=\linewidth]
{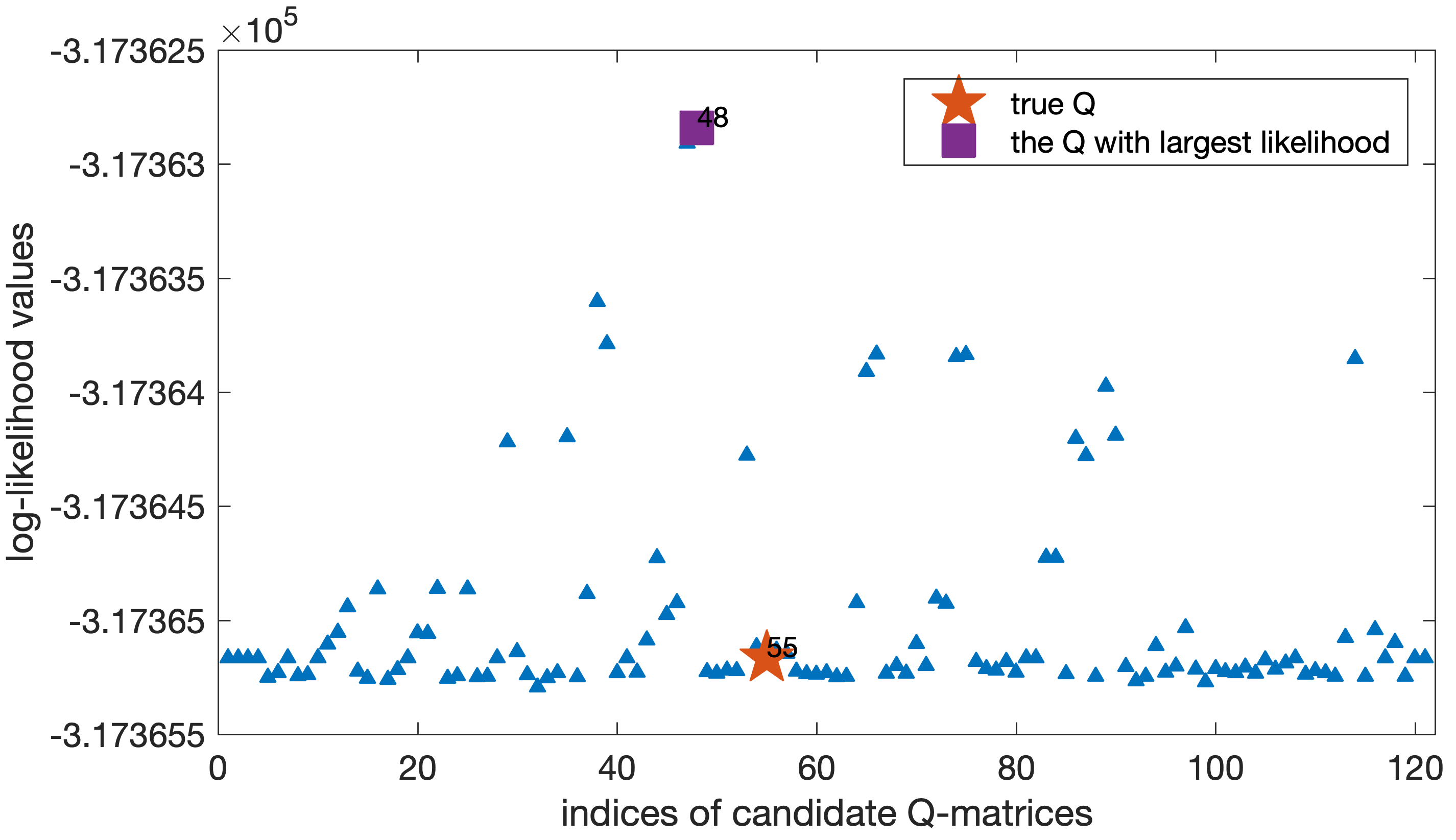}
\caption{true $Q$ not even locally identifiable:
$
Q^{55} = \begin{pmatrix}
     0 &  0  &   0  &   0  &   1 \\ 
     1  &   1   &  1   &  1   &  1
\end{pmatrix}^\top
$}
\label{fig-q55}
\end{subfigure}

\end{figure}

\newpage
\noindent
\textbf{Study IV: Verifying necessity of Condition $A$ ``completeness''.}

We verify the necessity of Condition $A$ ``\textit{completeness}" of the $Q$-matrix for identifiability. Consider two settings of incomplete $Q$-matrices, $Q_1$ with $(K, J)=(3,20)$ and $Q_2$ with $(K,J)=(5,20)$.
For $i=1,2$, for the matrix $Q=Q_i$ and arbitrary DINA model parameters $(\cc,\cg,\pp)$, we follow our theoretical derivations  to construct two alternative $Q$-matrices $Q'=Q'_i$ and $Q''=Q''_i$ and corresponding parameters $(\cc^{'},\cg^{'},\pp^{'})$ and $(\cc^{''},\cg^{''},\pp^{''})$.
Then we compute the marginal probabilities for all the possible $2^{20}\approx 10^6$ response patterns under each of the $Q$, $Q'$ and $Q''$, which characterize the distribution of the $20$-dimensional binary vector $\RR$. We give visualization plots to show how these different $Q$-matrices and different model parameters lead to exactly the same distribution of the observed responses $\RR$.

%\newpage
\begin{equation}\label{eq-inc-k3}
\small
Q_1 = \begin{pmatrix}
     1  &  0   &  0 \\ 
     0  &   1  &   0 \\ 
     \mathbf 1  & \mathbf  1  &  \mathbf 1 \\ 
     1  &   0  &   0 \\ 
     1  &   1  &   0 \\ 
    \mathbf 1  & \mathbf  1  &  \mathbf 1 \\ 
     1  &   0  &   0 \\ 
     1  &   1  &   0\\ 
     \mathbf 1  & \mathbf  1  &  \mathbf 1 \\ 
     1  &   0  &   0\\ 
     1  &   1  &   0\\ 
     \mathbf 1  & \mathbf  1  &  \mathbf 1 \\ 
     1  &   0  &   0\\ 
     1  &   1  &   0\\ 
     \mathbf 1  & \mathbf  1  &  \mathbf 1 \\ 
     1  &   0  &   0\\ 
     1  &   1  &   0\\ 
    \mathbf 1  & \mathbf  1  &  \mathbf 1 \\ 
     1  &   1  &   1\\ 
     1  &   1  &   1
\end{pmatrix}_{20\times 3}\quad\quad
Q'_1 = \begin{pmatrix}
     1  &  0   &  0 \\ 
     0  &   1  &   0 \\ 
     \mathbf 0  &  \mathbf 1  &  \mathbf 1 \\ 
     1  &   0  &   0 \\ 
     1  &   1  &   0 \\ 
     \mathbf 0  &  \mathbf 1  &  \mathbf 1\\ 
     1  &   0  &   0 \\ 
     1  &   1  &   0\\ 
     \mathbf 0  &  \mathbf 1  &  \mathbf 1\\ 
     1  &   0  &   0\\ 
     1  &   1  &   0\\ 
     \mathbf 0  &  \mathbf 1  &  \mathbf 1\\ 
     1  &   0  &   0\\ 
     1  &   1  &   0\\ 
     \mathbf 0  &  \mathbf 1  &  \mathbf 1\\ 
     1  &   0  &   0\\ 
     1  &   1  &   0\\ 
     \mathbf 0  &  \mathbf 1  &  \mathbf 1\\ 
     1  &   1  &   1\\ 
     1  &   1  &   1
\end{pmatrix}_{20\times 3}
\quad\quad
Q''_1 = \begin{pmatrix}
     1  &  0   &  0 \\ 
     0  &   1  &   0 \\ 
     \mathbf0  &   \mathbf 0  &   \mathbf 1 \\ 
     1  &   0  &   0 \\ 
     1  &   1  &   0 \\ 
     \mathbf0  &   \mathbf 0  &   \mathbf 1 \\ 
     1  &   0  &   0 \\ 
     1  &   1  &   0\\ 
     \mathbf0  &   \mathbf 0  &   \mathbf 1 \\ 
     1  &   0  &   0\\ 
     1  &   1  &   0\\ 
     \mathbf0  &   \mathbf 0  &   \mathbf 1 \\ 
     1  &   0  &   0\\ 
     1  &   1  &   0\\ 
     \mathbf0  &   \mathbf 0  &   \mathbf 1 \\ 
     1  &   0  &   0\\ 
     1  &   1  &   0\\ 
     \mathbf0  &   \mathbf 0  &   \mathbf 1 \\ 
     1  &   1  &   1\\ 
     1  &   1  &   1
\end{pmatrix}_{20\times 3}
\end{equation}
\begin{equation}\label{eq-inc-k5}
\small
%%%%%%%%
Q_2 = \begin{pmatrix}
     1  &   0  &   0  &   0  &   0\\ 
     0  &   1  &   0  &   0  &   0\\ 
     0  &   0  &   1  &   0  &   0\\ 
     0  &   0  &   0  &   1  &   0\\ 
     \mathbf1  &   \mathbf 1  &   \mathbf1  &  \mathbf 1  &  \mathbf 1\\ 
     1  &   0  &   0  &   0  &   0\\ 
     1  &   1  &   0  &   0  &   0\\ 
     1  &   1  &   1  &   0  &   0\\ 
     1  &   1  &   1  &   1  &   0\\ 
     \mathbf1  &   \mathbf 1  &   \mathbf1  &  \mathbf 1  &  \mathbf 1\\ 
     1  &   0  &   0  &   0  &   0\\ 
     1  &   1  &   0  &   0  &   0\\ 
     1  &   1  &   1  &   0  &   0\\ 
     1  &   1  &   1  &   1  &   0\\ 
     1  &   1  &   1  &  1  &   1\\ 
     1  &   0  &   0  &   0  &   0\\ 
     1  &   1  &   0  &   0  &   0\\ 
     1  &   1  &   1  &   0  &   0\\ 
     1  &   1  &   1  &   1  &   0\\ 
     1  &    1  &   1  &   1  &   1
\end{pmatrix}_{20\times 5}
~
Q'_2 = \begin{pmatrix}
     1  &   0  &   0  &   0  &   0\\ 
     0  &   1  &   0  &   0  &   0\\ 
     0  &   0  &   1  &   0  &   0\\ 
     0  &   0  &   0  &   1  &   0\\ 
     \mathbf0  &   \mathbf 0  &   \mathbf1  &  \mathbf 1  &  \mathbf 1\\ 
     1  &   0  &   0  &   0  &   0\\ 
     1  &   1  &   0  &   0  &   0\\ 
     1  &   1  &   1  &   0  &   0\\ 
     1  &   1  &   1  &   1  &   0\\ 
      \mathbf0  &   \mathbf 0  &   \mathbf1  &  \mathbf 1  &  \mathbf 1\\ 
     1  &   0  &   0  &   0  &   0\\ 
     1  &   1  &   0  &   0  &   0\\ 
     1  &   1  &   1  &   0  &   0\\ 
     1  &   1  &   1  &   1  &   0\\ 
     1  &    1  &   1  &   1  &   1\\ 
     1  &   0  &   0  &   0  &   0\\ 
     1  &   1  &   0  &   0  &   0\\ 
     1  &   1  &   1  &   0  &   0\\ 
     1  &   1  &   1  &   1  &   0\\ 
      1  &    1  &   1  &   1  &   1
\end{pmatrix}_{20\times 5}
~
Q''_2 = \begin{pmatrix}
     1  &   0  &   0  &   0  &   0\\ 
     0  &   1  &   0  &   0  &   0\\ 
     0  &   0  &   1  &   0  &   0\\ 
     0  &   0  &   0  &   1  &   0\\ 
     \mathbf0  &   \mathbf 0  &   \mathbf0  &  \mathbf0  &  \mathbf 1\\ 
     1  &   0  &   0  &   0  &   0\\ 
     1  &   1  &   0  &   0  &   0\\ 
     1  &   1  &   1  &   0  &   0\\ 
     1  &   1  &   1  &   1  &   0\\ 
    \mathbf0  &   \mathbf 0  &   \mathbf0  &  \mathbf0  &  \mathbf 1\\ 
     1  &   0  &   0  &   0  &   0\\ 
     1  &   1  &   0  &   0  &   0\\ 
     1  &   1  &   1  &   0  &   0\\ 
     1  &   1  &   1  &   1  &   0\\ 
    1  &    1  &   1  &   1  &   1\\ 
     1  &   0  &   0  &   0  &   0\\ 
     1  &   1  &   0  &   0  &   0\\ 
     1  &   1  &   1  &   0  &   0\\ 
     1  &   1  &   1  &   1  &   0\\ 
      1  &    1  &   1  &   1  &   1
      \end{pmatrix}_{20\times 5}
\end{equation}

\begin{figure}[H]
\caption{DINA: true $Q$-matrix of size $20\times 5$ is not complete and hence not identifiable.}
\label{fig-basis-k5}

\centering
\begin{subfigure}{0.58\textwidth}
\includegraphics[width=\linewidth]
{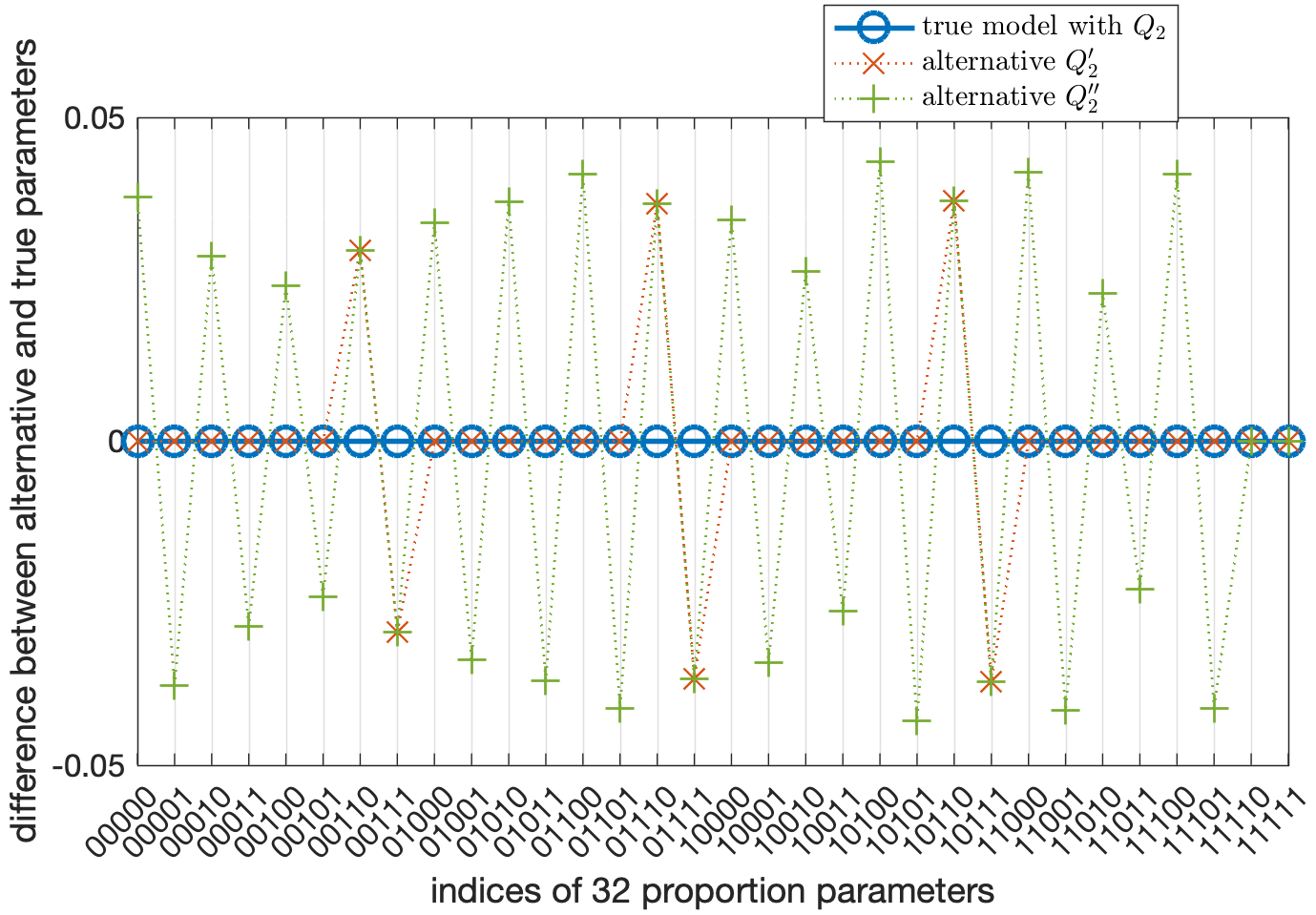}
\caption{$K=5$, $J=20$, three sets of parameters} 
\end{subfigure}

\begin{subfigure}{0.58\textwidth}
\includegraphics[width=\linewidth]
{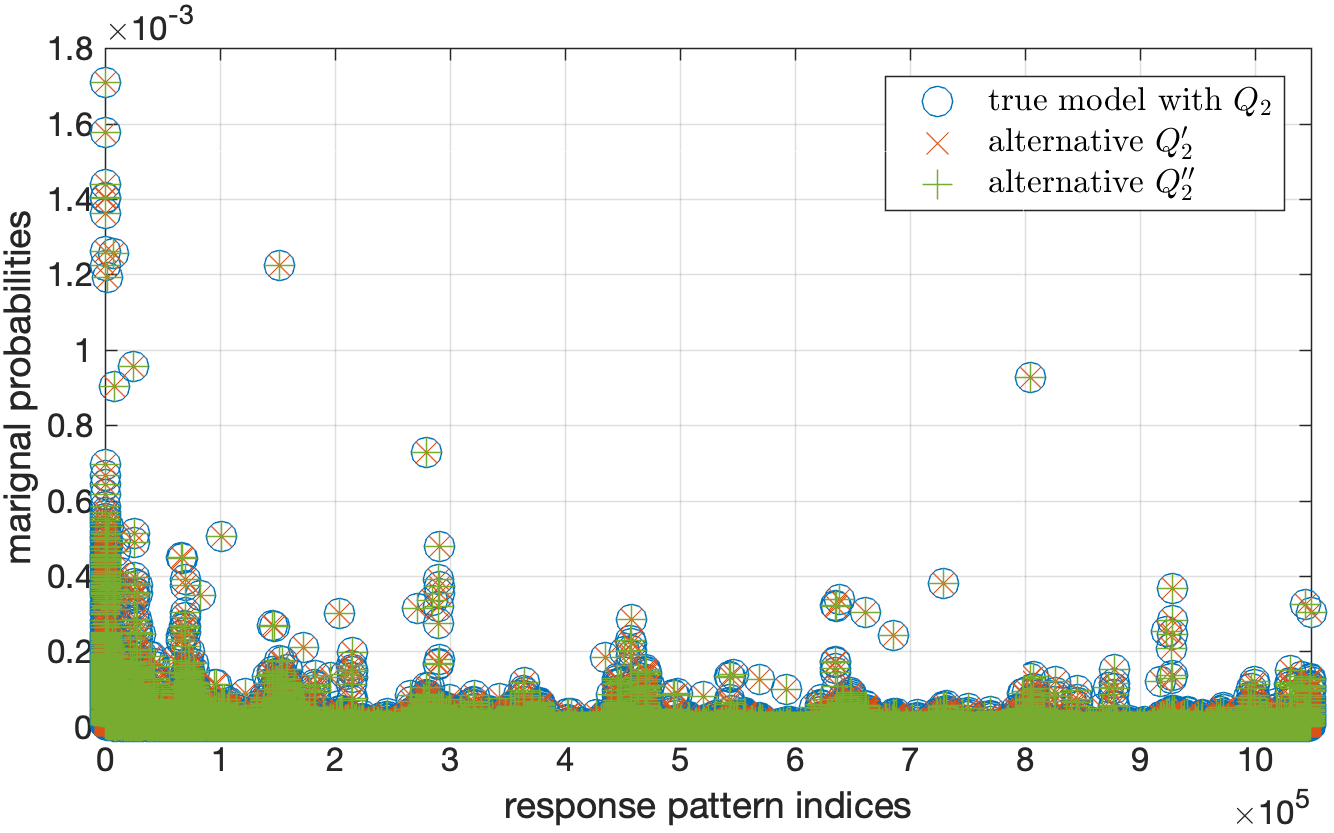}
\caption{$K=5$, $J=20$, $|\{0,1\}^{20}|=2^{20}$ response probabilities}
\end{subfigure}

\begin{subfigure}{0.58\textwidth}
\includegraphics[width=\linewidth]
{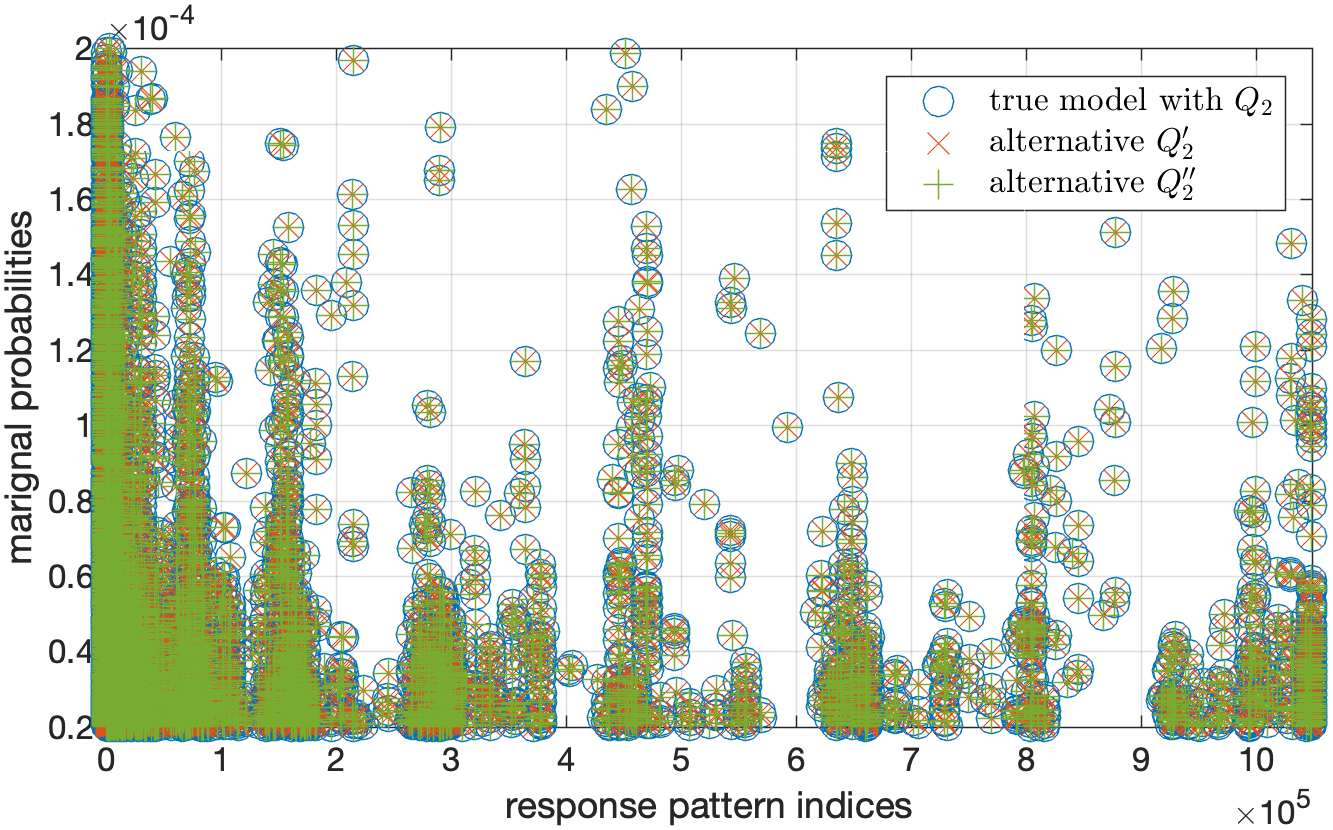}
\caption{$K=5$, $J=20$, response probabilities zoomed in}
\end{subfigure}
\end{figure}

%In particular, we set the item parameters to be the same as those under the true $Q^i$, i.e., $\cc^{''}=\cc^{'}=\cc$ and $\cg^{''}=\cg^{'}=\cg$, and only change the proportion parameters under the alternative $Q$-matrices $Q'_i$ and $Q''_i$. For instance, for the case of $K=3$, we set $\pp^{'}$ and $\pp^{''}$ in the following way,
%\newpage
First, consider the following $Q_1$ with $(K,J)=(3,20)$ in \eqref{eq-inc-k3}, which is incomplete because its row vectors does not contain the unit vector $(0, 0, 1)$. For arbitrarily generated parameters $(\cc,\cg,\pp)$, we set $\cc^{''}=\cc^{'}=\cc$ and $\cg^{''}=\cg^{'}=\cg$ and set the proportion parameters as follows,
\begin{align}
\label{eq-ind-k3-prop}
\begin{cases}
 p'_{(011)}=0,\\
 p'_{(010)} = p_{(010)}+p_{(011)},\\
 p'_{\aaa} = p_{\aaa},~\forall\aaa\neq(011),(010);
%\bar p^1_{(011)}=0,\\
% \bar p^1_{(010)} = p_{(010)}+p_{(011)}; 
\end{cases} \quad
\begin{cases}
 p''_{(001)}= p^2_{(011)}= p^2_{(111)}=0,\\
 p''_{(000)} = p_{(000)}+p_{(001)},\\
%\bar p^2_{(011)}=0,\\
 p''_{(010)} = p_{(010)}+p_{(011)}, \\
%\bar p^2_{(111)}=0,\\
 p''_{(110)} = p_{(110)}+p_{(111)}, \\
 p''_{\aaa} = p_{\aaa},~\forall\aaa=(100),(101).
\end{cases}
\end{align}
We define a notation $\Gamma(Q)$ to briefly explain the rationale behind the above constructions. The $\Gamma(Q)$ is a $J\times 2^K$ binary matrix defined based on $Q$. The columns and rows of $\Gamma(Q)$ are indexed by the $J$ items and the $2^K$ possible attribute patterns, respectively; and the $(j,\aaa)$th entry of it is defined to be $\Gamma_{j,\aaa}(Q)=I(\aaa\succeq\qq_j)$.
An important observation is that, due to the forms of $Q$, $Q'$ and $Q''$, the unique column vectors in $\Gamma(Q)$ form a subset of those of $\Gamma(Q')$; and further the unique column vectors of $\Gamma(Q')$ form a subset of those of $\Gamma(Q'')$.
Therefore, to construct $\pp'$ such that $(Q,\cc,\cg,\pp)$ and $(Q',\cc,\cg,\pp')$ that are non-distinguishable, we only need to set $p'_{\aaa}=0$ for those $\aaa$ whose corresponding column vector in $\Gamma(Q')$ does not appear as the column vector of $\Gamma(Q)$; and let the proportions (in vector $\pp'$) of other attribute patterns to absorb the proportions of these $\aaa$'s in the vector $\pp'$. The proportions $\pp''$ under $Q''$ are constructed similarly.
This is exactly how Equation \eqref{eq-ind-k3-prop} are derived.
For the $Q_2$, $Q'_2$ and $Q''_2$ defined in \eqref{eq-inc-k5}, we construct the proportion parameters $\pp'$ under $Q'_2$ and $\pp''$ under $Q''_2$ following the same rationale; the details of defining them are omitted but their values are later revealed in Figure \ref{fig-basis-k5}(c).
%$T(Q,\cc,\cg)\pp = T(Q',\cc,\cg)\pp''=T(Q'',\cc,\cg)\pp''$ under the DINA model.

\begin{figure}[H]
\caption{DINA: true $Q$-matrix of size $20\times 3$ is not complete and hence not identifiable.}
\label{fig-basis-k3}

\centering
\begin{subfigure}{0.58\textwidth}
\includegraphics[width=\linewidth]
{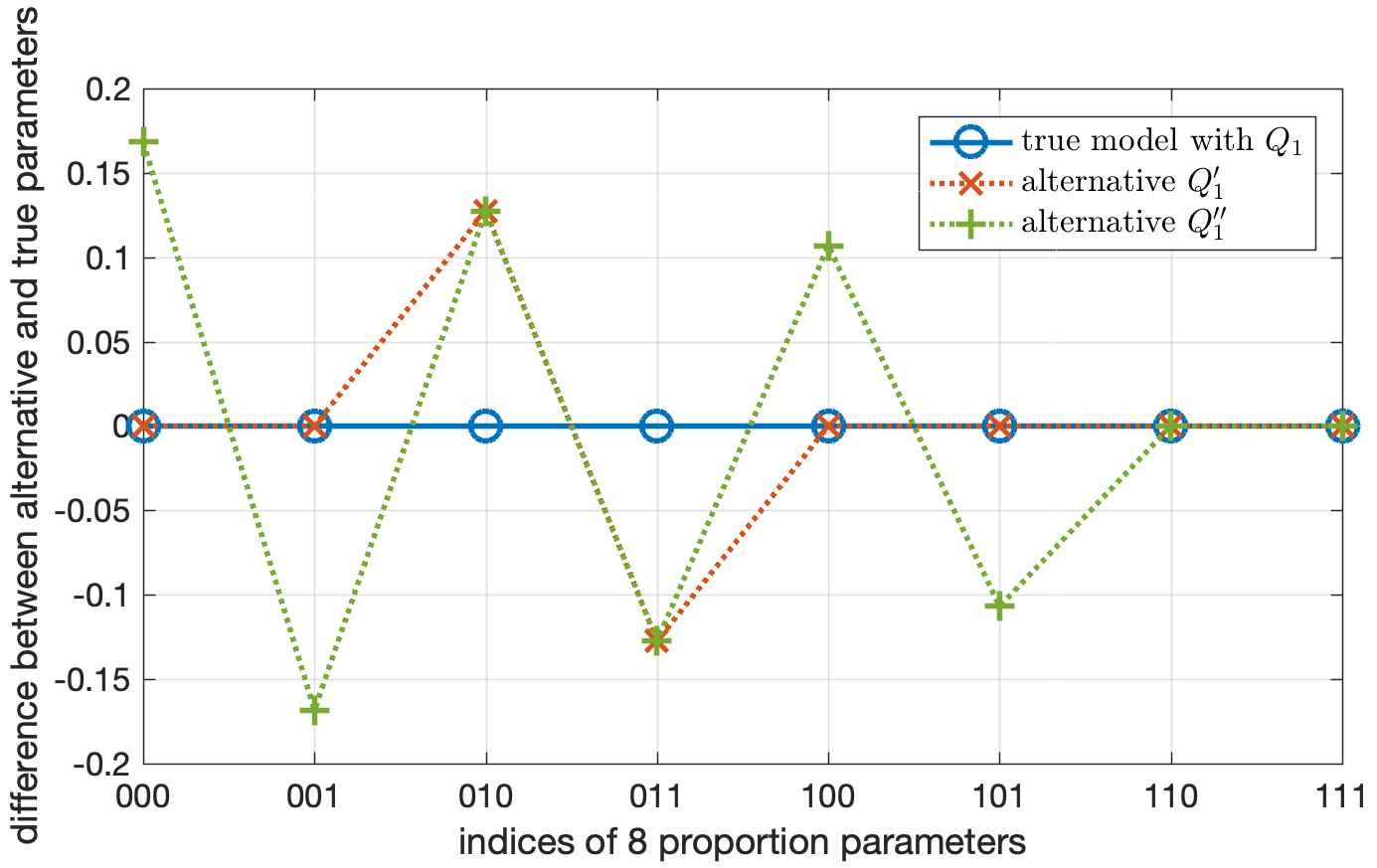}
\caption{$K=3$, $J=20$, three sets of parameters} 
\end{subfigure}

\begin{subfigure}{0.58\textwidth}
\includegraphics[width=\linewidth]
{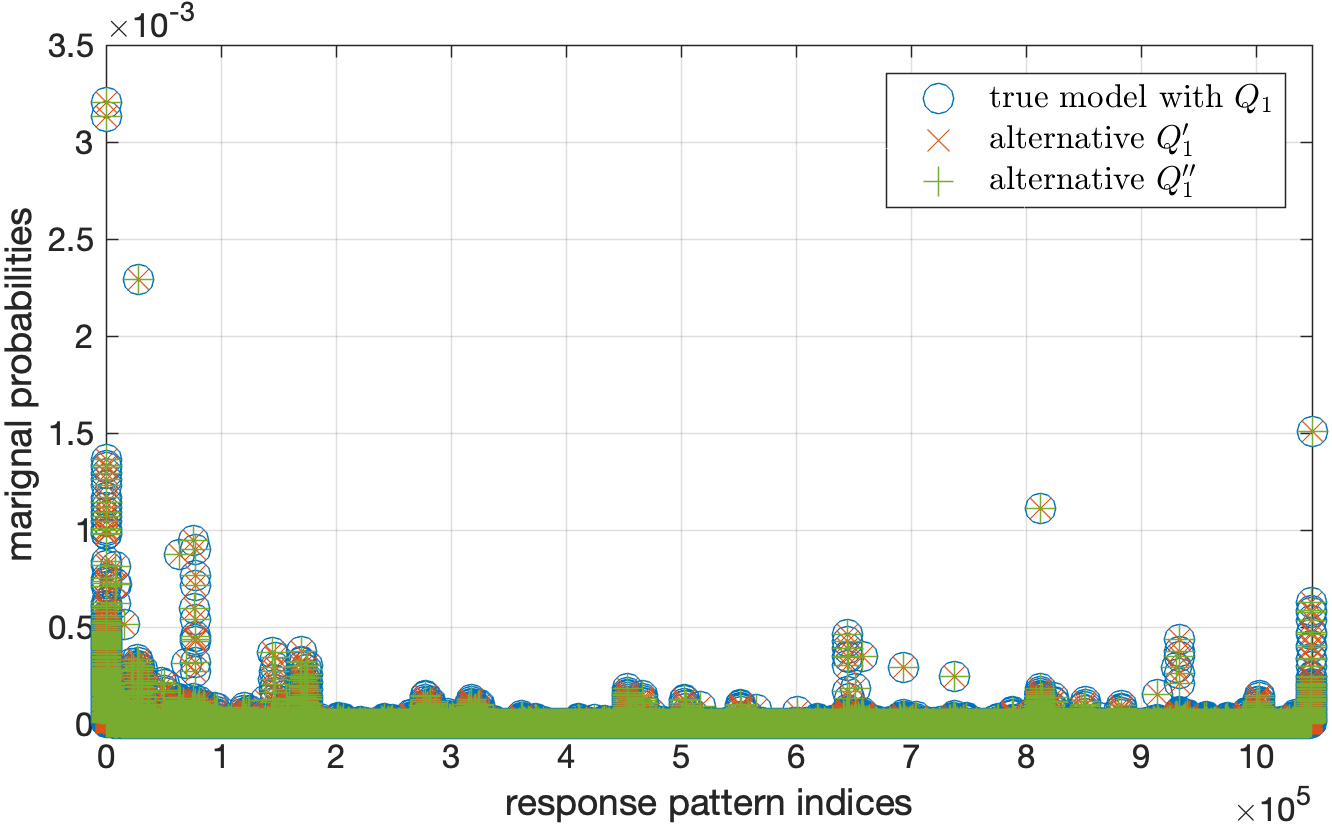}
\caption{$K=3$, $J=20$, $|\{0,1\}^{20}|=2^{20}$ response probabilities}
\end{subfigure}

\begin{subfigure}{0.58\textwidth}
\includegraphics[width=\linewidth]
{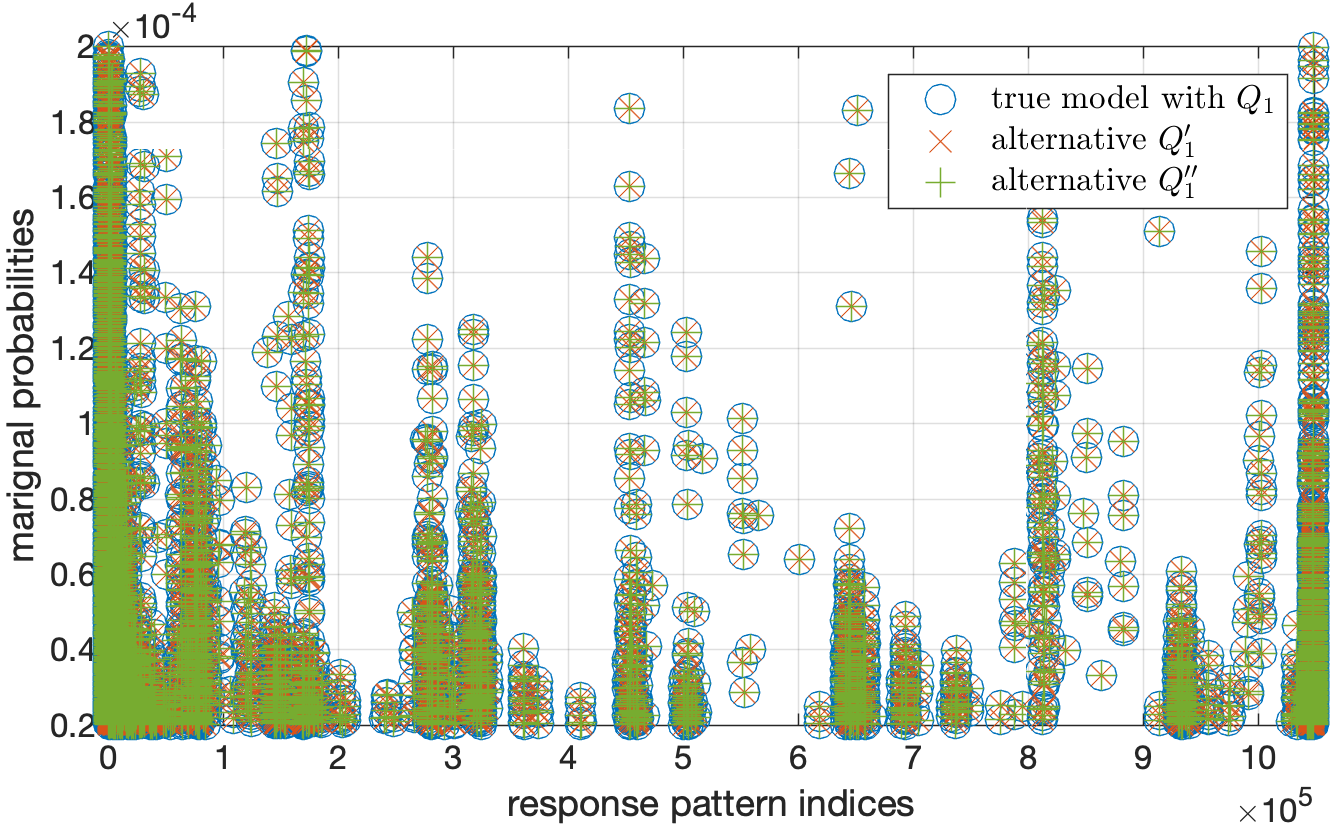}
\caption{$K=3$, $J=20$, response probabilities zoomed in}
\end{subfigure}

\end{figure}

In Figure \ref{fig-basis-k3}, we visualize the non-identifiability phenomenon of $Q_1$. In Figure \ref{fig-basis-k3}(a), we plot the differences of proportions parameters under the alternative models and the true model with $Q_1$. The red dotted line with ``$\times$" plots the values
$\pp'-\pp=(p'_{000}-p_{000},~p'_{001}-p_{001},~p'_{010}-p_{010},~p'_{011}-p_{011},
~p'_{100}-p_{100},~p'_{101}-p_{101},~p'_{110}-p_{110},~p'_{111}-p_{111})$ correspondent to the 8 attribute patterns; and the green dotted line with ``$+$" plots $\pp''-\pp$. Despite these three sets of parameters are quite different, the $2^{20}$-dimensional vector of marginal probabilities of $\RR$ are exactly the same, as shown in plots (b) and (c) of Figure \ref{fig-basis-k3}. In particular, in plot (b), the $x$-axis presents the indices of the response patterns in $\rr\in\{0,1\}^J$, the $y$-axis presents the values of $\mathbb P(\RR=\rr\mid Q,\cc,\cg,\pp)$, where the blue circles denote those under $(Q_1,\pp)$, red ``$\times$" for $(Q'_1,\pp')$, and green ``$+$'' for $(Q''_1,\pp'')$. Plot (c) of Figure \ref{fig-basis-k3} is a zoomed-in version of plot (b), by only showing those marginal probabilities in $[0.2\times 10^{-4}, 2\times 10^{-4}]$, which contains around $7\times 10^3$ response patterns. One can roughly see from both plots (b) and (c) that the three underlying parameters yield identical distribution of the response vector.
Indeed, the computation carried out using Matlab yields
\begin{align*}
\max_{\rr\in\{0,1\}^{20}}\left|\mathbb P(\RR=\rr\mid Q_1,\cc,\cg,\pp)-\mathbb P(\RR=\rr\mid Q'_1,\cc,\cg,\pp')\right|
&= 2.17\times 10^{-19},\\
\max_{\rr\in\{0,1\}^{20}}\left|\mathbb P(\RR=\rr\mid Q_1,\cc,\cg,\pp)-\mathbb P(\RR=\rr\mid Q''_1,\cc,\cg,\pp'')\right|
&= 4.34\times 10^{-19},
\end{align*}
which are both smaller than the machine epsilon (machine error) of Matlab $2.22\times 10^{-16}$. This confirms that $Q_1$ defined in \eqref{eq-inc-k3} is not identifiable.

Figure \ref{fig-basis-k5} shows the analogous results for $Q_2$ of size $20\times 5$. Plot (a) in Figure \ref{fig-basis-k5} shows the difference of the $2^5=32$-dimensional proportion parameters under alternative and true $Q$-matrices, and plots (b) and (c) give marginal probabilities of $\RR$. The computation using Matlab gives
\begin{align*}
\max_{\rr\in\{0,1\}^{20}}\left|\mathbb P(\RR=\rr\mid Q_2,\cc,\cg,\pp)-\mathbb P(\RR=\rr\mid Q'_2,\cc,\cg,\pp')\right|
&= 2.17\times 10^{-19},\\
\max_{\rr\in\{0,1\}^{20}}\left|\mathbb P(\RR=\rr\mid Q_2,\cc,\cg,\pp)-\mathbb P(\RR=\rr\mid Q''_2,\cc,\cg,\pp'')\right|
&= 6.51\times 10^{-19},
\end{align*}
which are also both smaller than the machine error $2.22\times 10^{-16}$ of Matlab. This verifies the non-identifiability of $Q_2$ defined in \eqref{eq-inc-k5}.

\newpage

\subsection{General RLCM: GDINA Model}\label{sec-simu-gdina}
In this section, we design simulation studies to verify the proposed identifiability conditions under the GDINA model introduced in Example \ref{exp-gdina}. In Study V, we use exhaustive search within $5\times 2$ $Q$-matrices to verify the sufficient conditions in Theorem \ref{thm-multi}.
In Study VI and Study VII, we verify the necessary conditions in Theorem \ref{prop-gdina-nece}.

\smallskip
\noindent
\textbf{Study V: When $Q$-matrix satisfies Conditions $D$, $E$ for generic identifiability.}

Within the set of $5\times 2$ $Q$-matrices $\text{Exhaus}(Q_{5\times 2})$, if $Q$ satisfies the sufficient conditions $D$ and $E$ for generic identifiability under the GDINA model, then other than the all-one $Q$-matrix $Q^{121}$ which corresponds to the unrestricted latent class model, $Q$ can take the forms of $Q^{15}$, $Q^{18}$, $Q^{27}$, $Q^{54}$, and $Q^{81}$ (up to rearrangement of rows and columns). 
When using some $Q$-matrix to generate data, we also set the sample size to $N=10^5$ and randomly set the true parameters which satisfy the monotonicity constraint \eqref{eq-mono} in the main text.
In plots (a), (b), (c), (d) and (e) in Figure \ref{fig-gdina}, we present the exhaustive search results when the true data-generating $Q$-matrix is $Q^{15}$, $Q^{18}$, $Q^{27}$, $Q^{54}$, or $Q^{81}$. We point out that for GDINA model, in each scenario, we did not plot all the 121 $Q$-matrices' log-likelihood values, although we fit all the 121 ones to the simulated data. Instead, we only plot those $Q$-matrices under which the estimated parameters satisfies the stringent monotonicity constraint
\begin{equation}
	\label{eq-strmono}
	\theta_{j,\aaa}>\theta_{j,\aaa'}~~\text{if}~~\aaa\odot \qq_j\succ \aaa'\odot \qq_j.
\end{equation}
This constraint is stronger than requiring merely \eqref{eq-mono}, and it is often imposed in practice when fitting the general RLCM that models the main and interaction effects of the latent attributes; for example, see the LCDM proposed in \cite{HensonTemplin09}. So each blue triangle in each plot of Figure \ref{fig-gdina} corresponds to a $Q$-matrix with estimated $\TT$ satisfying \eqref{eq-strmono}.
We can see from the five plots in Figure \ref{fig-gdina} that when the generic identifiability conditions $D$ and $E$ are satisfied, the true data-generating $Q$-matrix achieves the maximum of the data likelihood compared to all the candidate $Q$-matrices of the same size.

\newpage
\begin{figure}[H]
\caption{GDINA: exhaustive search in $5\times 2$ $Q$-matrices with a true $Q$ satisfying Conditions $D$ and $E$.}
\label{fig-gdina}
\centering
\begin{subfigure}{0.7\textwidth}
\includegraphics[width=\linewidth]
{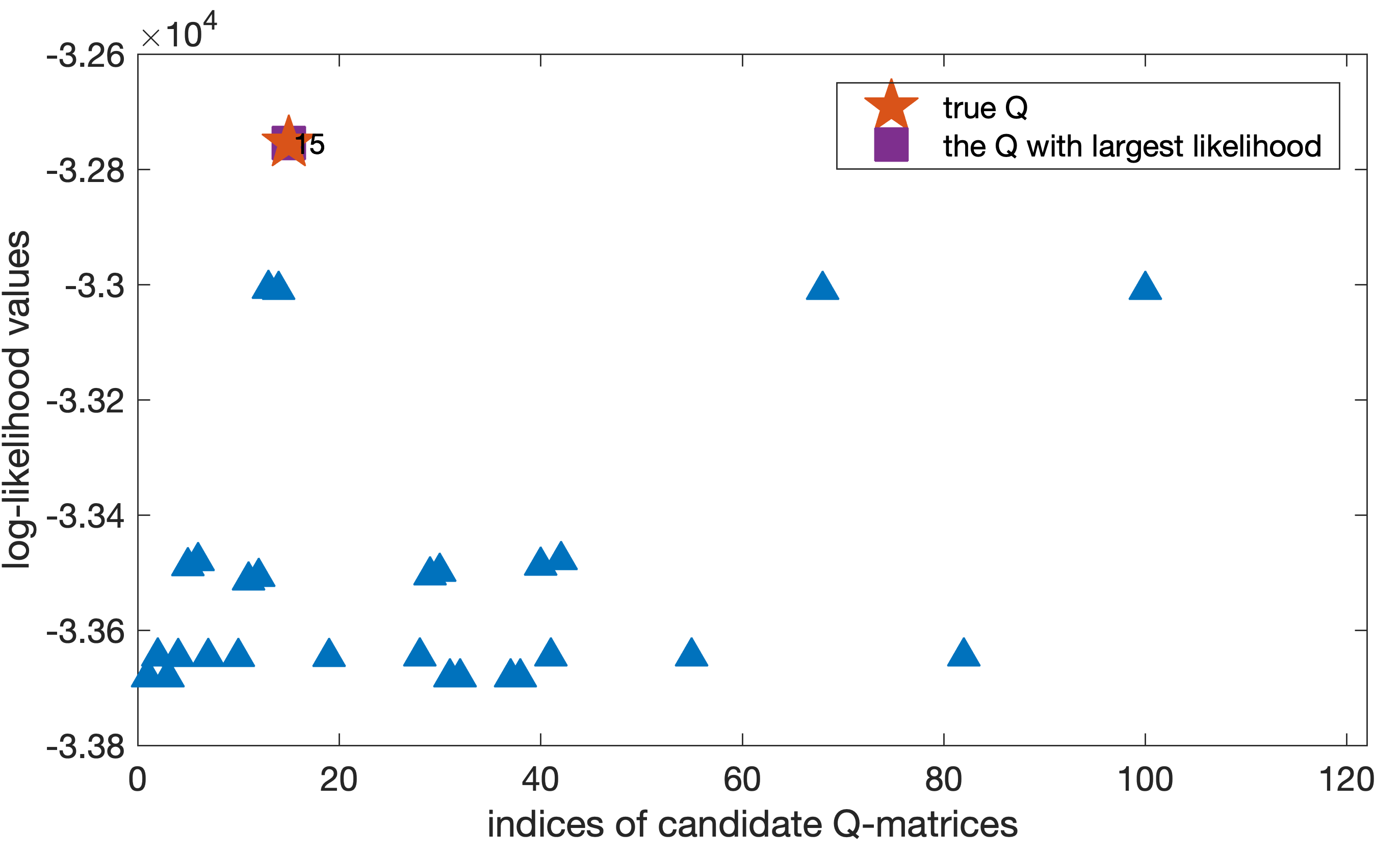}
\caption{GDINA: generically identifiable:
$
\small
Q^{15} = \begin{pmatrix}
  0   &  1  &   1   &  1   &  0\\
     1 &    1   &  0  &   0  &   1
\end{pmatrix}^\top
$} 
\end{subfigure}
\end{figure}

\begin{figure}[H]\ContinuedFloat
\centering
\begin{subfigure}{0.7\textwidth}
\includegraphics[width=\linewidth]
{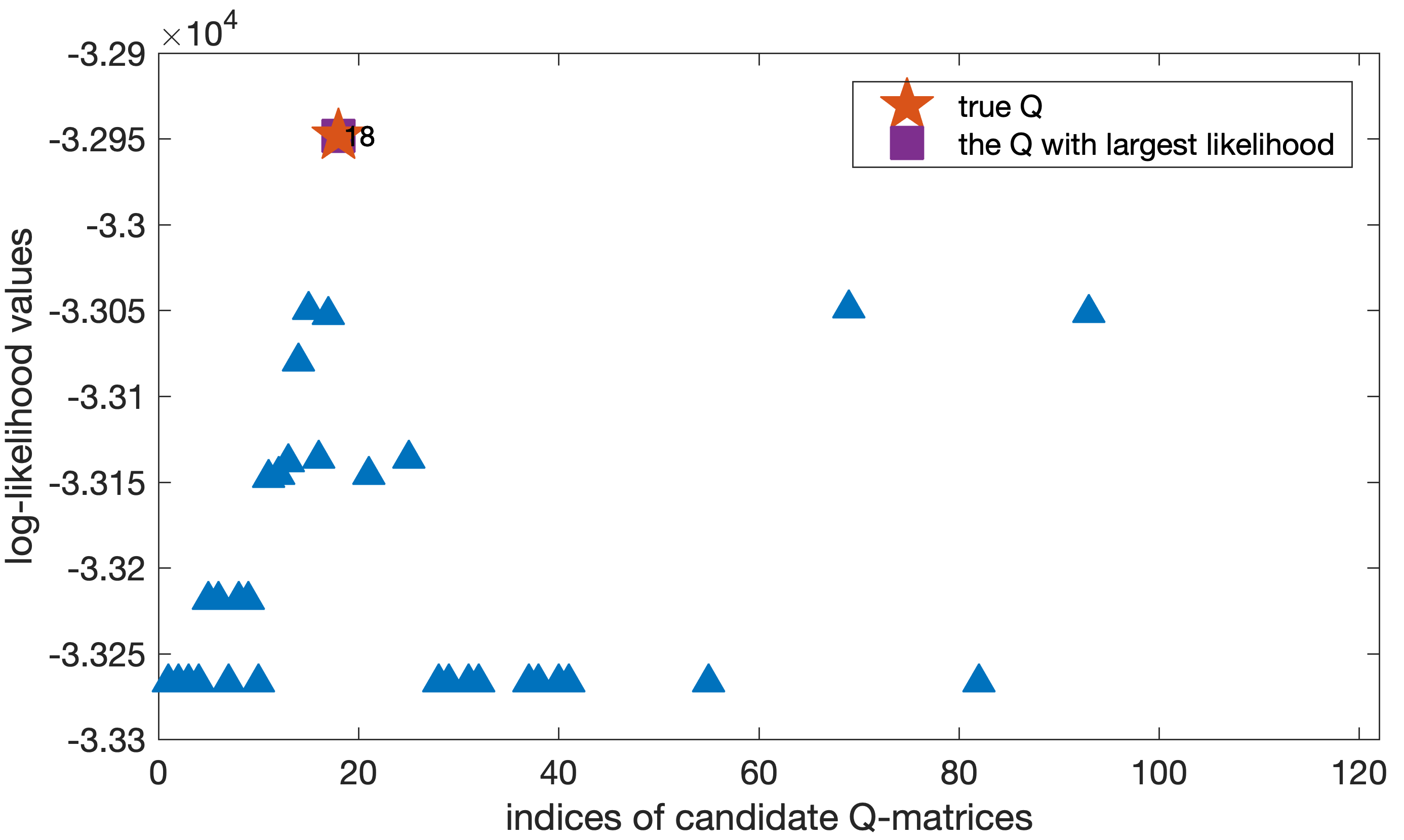}
\caption{GDINA: generically identifiable:
$\small
Q^{18} = \begin{pmatrix}
  0   &  1  &   1   &  1   &  0\\
     1 &    1   &  1  &   0  &   1
     \end{pmatrix}^\top
$}
\end{subfigure}
\end{figure}

\begin{figure}[H]\ContinuedFloat
\centering
\begin{subfigure}{0.7\textwidth}
\includegraphics[width=\linewidth]
{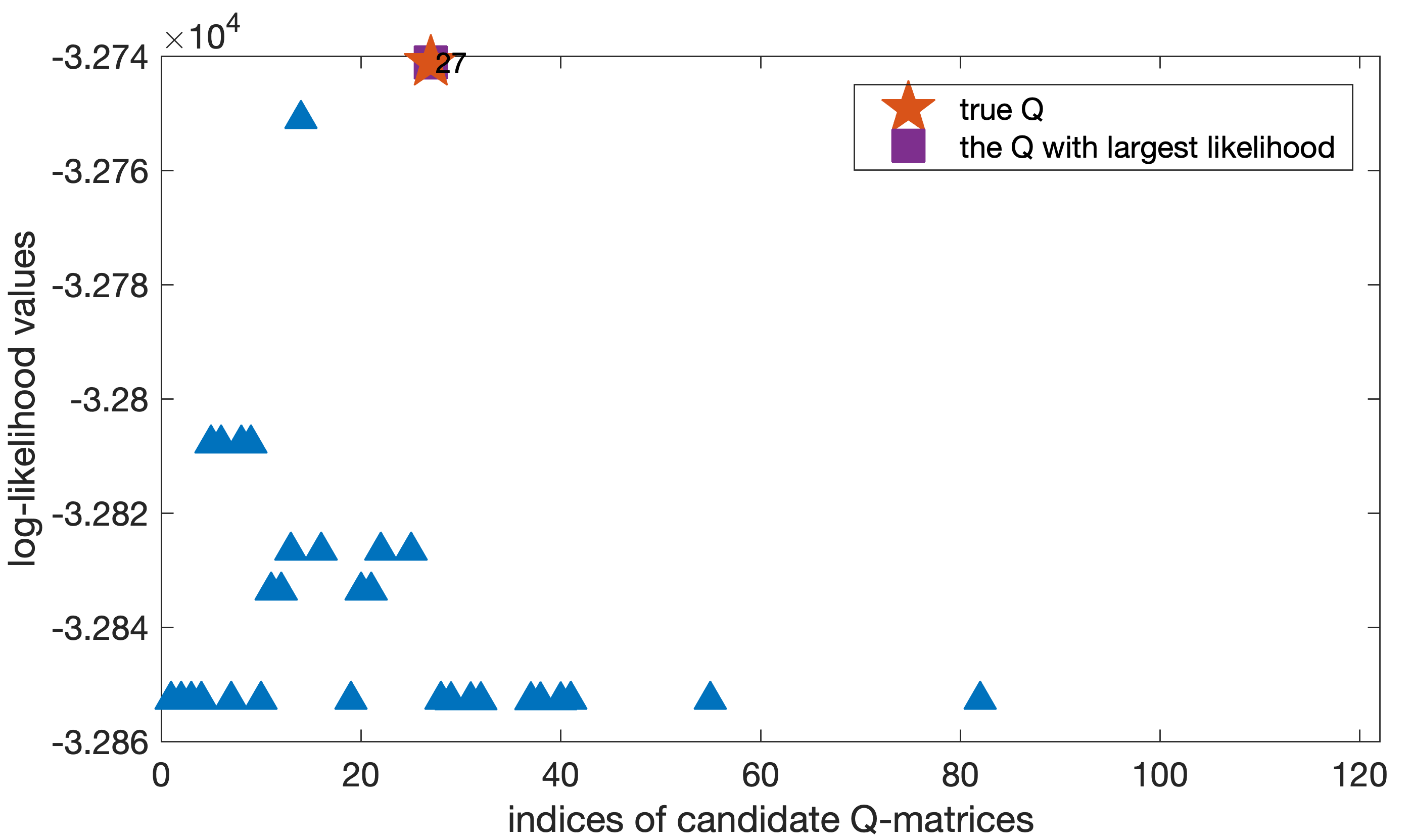}
\caption{GDINA: generically identifiable:
$\small
Q^{27} = \begin{pmatrix}
  0   &  1  &   1   &  1   &  0\\
  1 &    1   &  1  &   1  &   1
\end{pmatrix}^\top
$}
\end{subfigure}
\end{figure}

\begin{figure}[H]\ContinuedFloat
\centering
\begin{subfigure}{0.7\textwidth}
\includegraphics[width=\linewidth]
{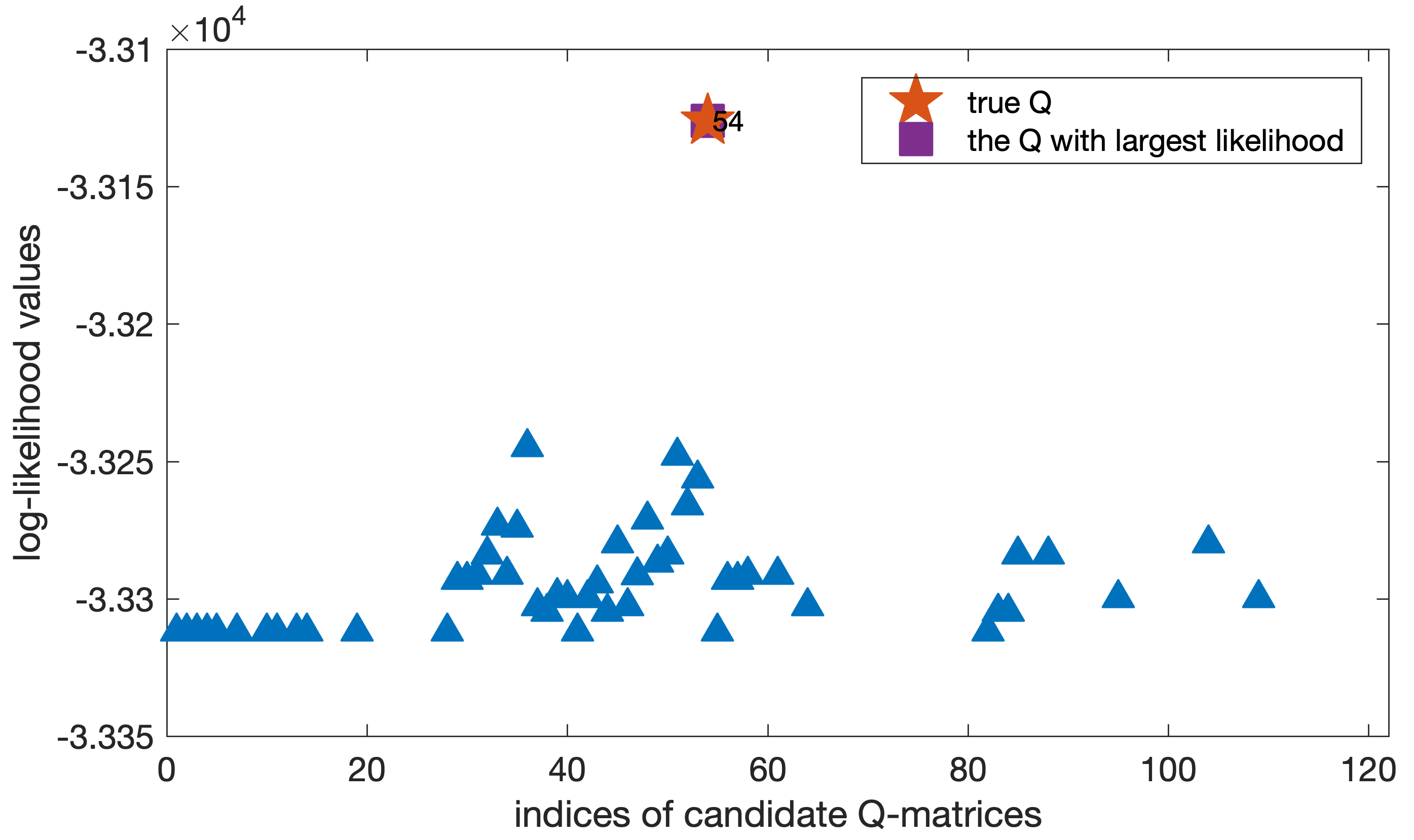}
\caption{GDINA: generically identifiable:
$\small
Q^{54} = \begin{pmatrix}
  0   &  1  &   1   &  1   &  1\\
  1 &    1   &  1  &   1  &   0
\end{pmatrix}^\top
$}
\end{subfigure}
\end{figure}

\newpage
\begin{figure}[H]\ContinuedFloat
\centering
\begin{subfigure}{0.7\textwidth}
\includegraphics[width=\linewidth]
{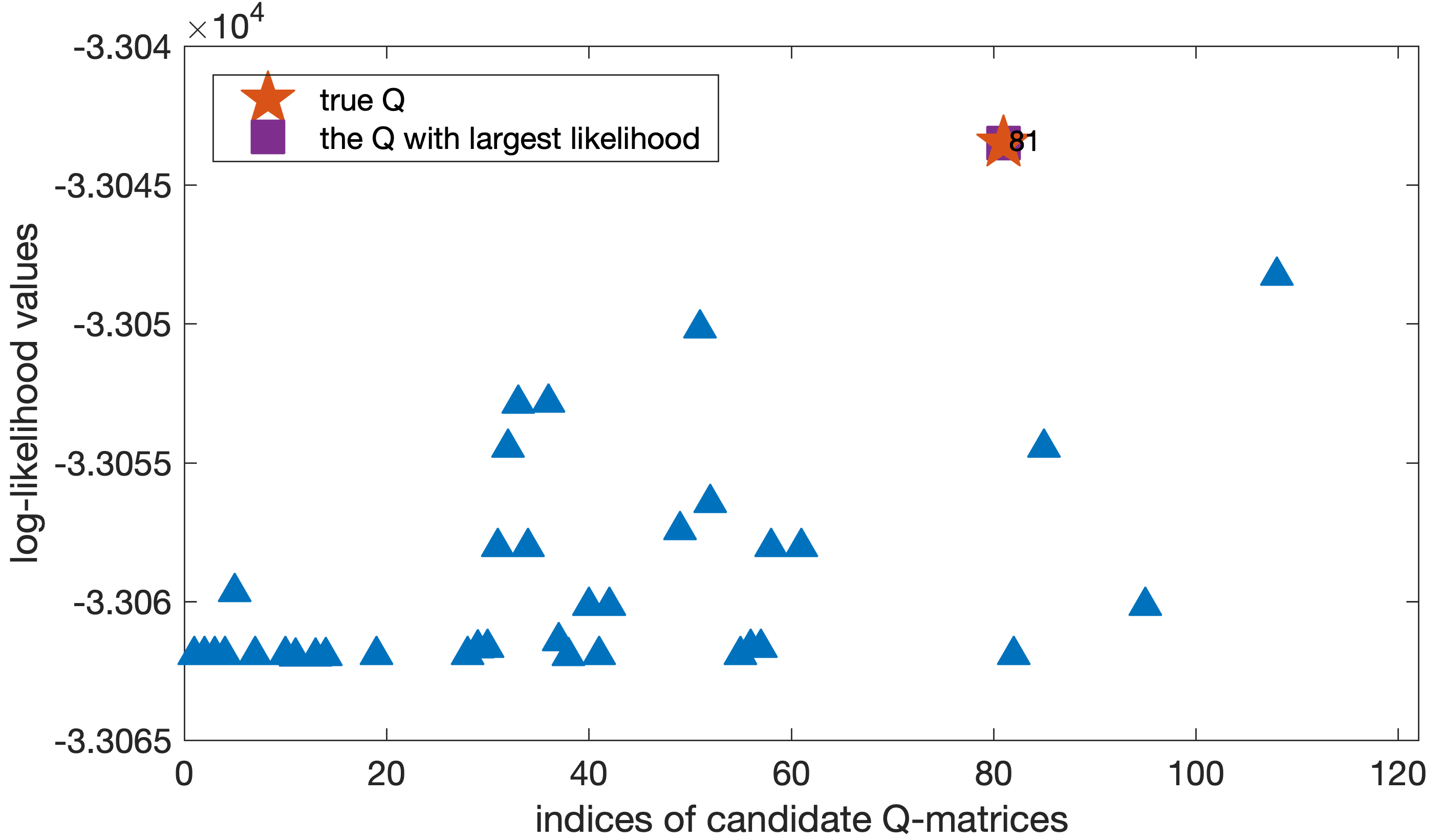}
\caption{GDINA: generically identifiable:
$\small
Q^{81} = \begin{pmatrix}
  0   &  1  &   1   &  1   &  1\\
  1 &    1   &  1  &   1  &   1
\end{pmatrix}^\top
$}
\end{subfigure}
\end{figure}

%\begin{figure}[H]\ContinuedFloat
%\centering
%\begin{subfigure}{0.8\textwidth}
%\includegraphics[width=\linewidth]
%{exhaus_gdina_onlymono_Q121}
%\caption{GDINA: generically identifiable:
%$\small
%Q^{121} = \begin{pmatrix}
%  1   &  1  &   1   &  1   &  1\\
%  1 &    1   &  1  &   1  &   1
%\end{pmatrix}^\top
%$}
%\end{subfigure}
%\end{figure}

\noindent
\textbf{Study VI: When $Q$-matrix does not even lead to local generic identifiability.}

We now use the not even locally generically identifiable $Q$-matrices $Q^1$, $Q^2$, or $Q^3$ to generate the data, and perform the exhaustive search among $\text{Exhaus}(Q_{5\times 2})$. The log-likelihood plots along with the forms of the data generating matrices $Q^1$, $Q^2$, $Q^3$ are presented in Figure \ref{fig-gdina-search}. Similar to the previous Study V, here in each scenario we only plot those $Q$-matrix whose estimated $\TT$ parameters satisfy the stringent monotonicity constraint \eqref{eq-strmono}.
One can see from the plots in Figure \ref{fig-gdina-search} that these $Q^1$, $Q^2$, $Q^3$ do not maximize the data likelihood, implying severe non-identifiability. Note that for Figure \ref{fig-gdina-search}(b) corresponding to $Q^2$, there are only two $Q$-matrices satisfying the constraint \eqref{eq-strmono} among the 121 $Q$-matrices fitted to the data; these two $Q$-matrices are the true $Q$-matrix $Q^2$ and another $Q$-matrix $Q^{56}$,
\begin{equation*}%\small
Q^{2}=\begin{pmatrix}
     0  &   1\\
     1  &   0\\
     0  &   1\\
     0  &   1\\
     0  &   1
     \end{pmatrix},\qquad
Q^{56}=\begin{pmatrix}
     0  &   1\\
     1  &   0\\
     0  &   1\\
     0  &   1\\
     1  &   1
     \end{pmatrix}.
\end{equation*}
Note that even there are only two $Q$-matrices satisfying the monotonicity constraint \eqref{eq-strmono}, the true $Q^{2}$ used to generate the data is not the one that has the larger likelihood, according to Figure \ref{fig-gdina-search}(b). This illustrates the non-identifiability of $Q_2$.

%\newpage
\begin{figure}[H]
\caption{GDINA: exhaustive search in $5\times 2$ $Q$-matrices with a true $Q$-matrix which leads to a not even locally generically identifiable model.}
\label{fig-gdina-search}
\centering

\begin{subfigure}{0.75\textwidth}
\includegraphics[width=\linewidth]
{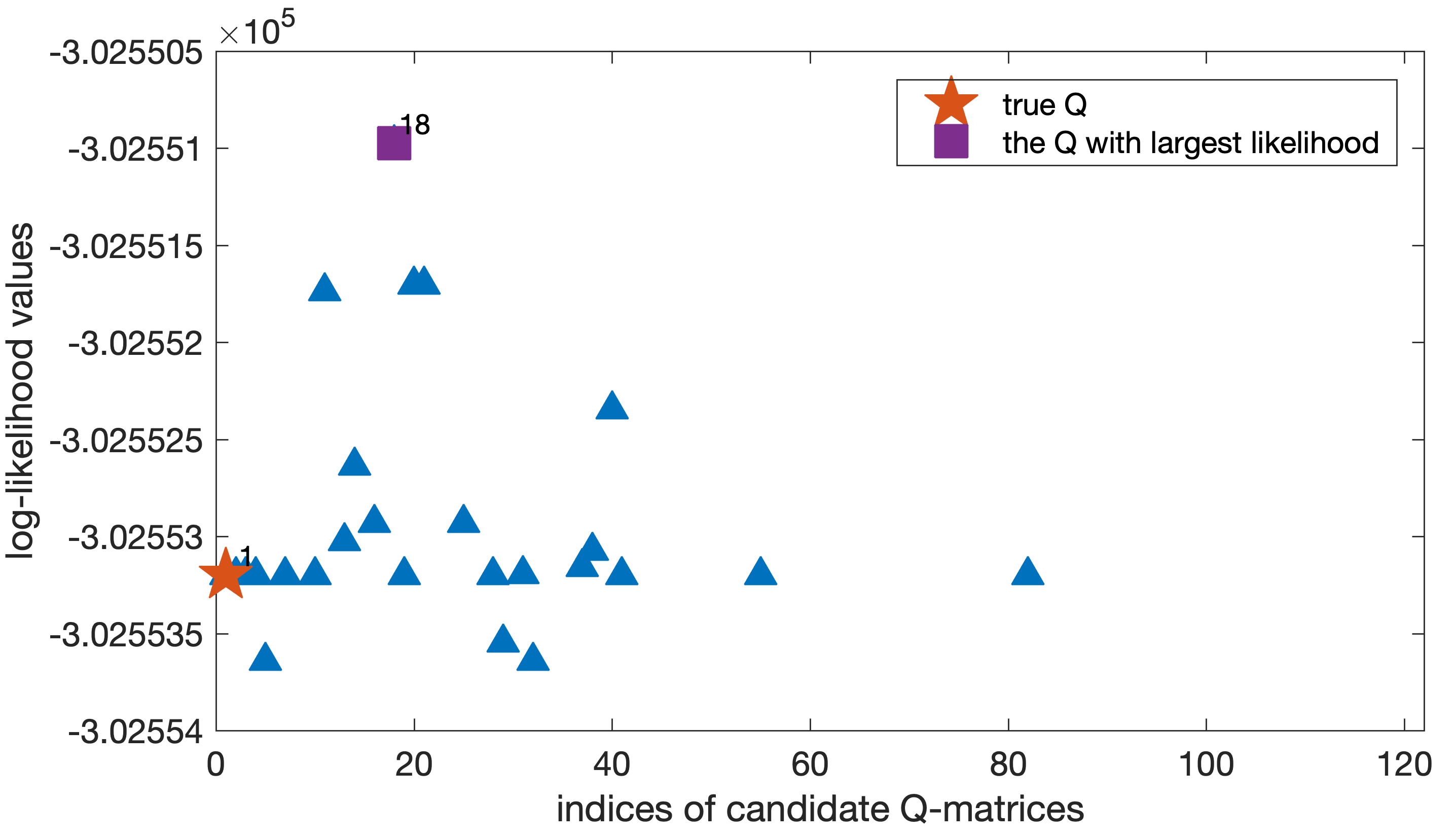}
\caption{GDINA: true $Q$ not even locally identifiable: $\small
Q^{1} = \begin{pmatrix}
  0   &  0  &   0   &  0   &  0\\
  1 &    1   &  1  &   1  &   1
\end{pmatrix}^\top
$} 
\end{subfigure}
\end{figure}

\begin{figure}[H]\ContinuedFloat
\centering

\begin{subfigure}{0.75\textwidth}
\includegraphics[width=\linewidth]
{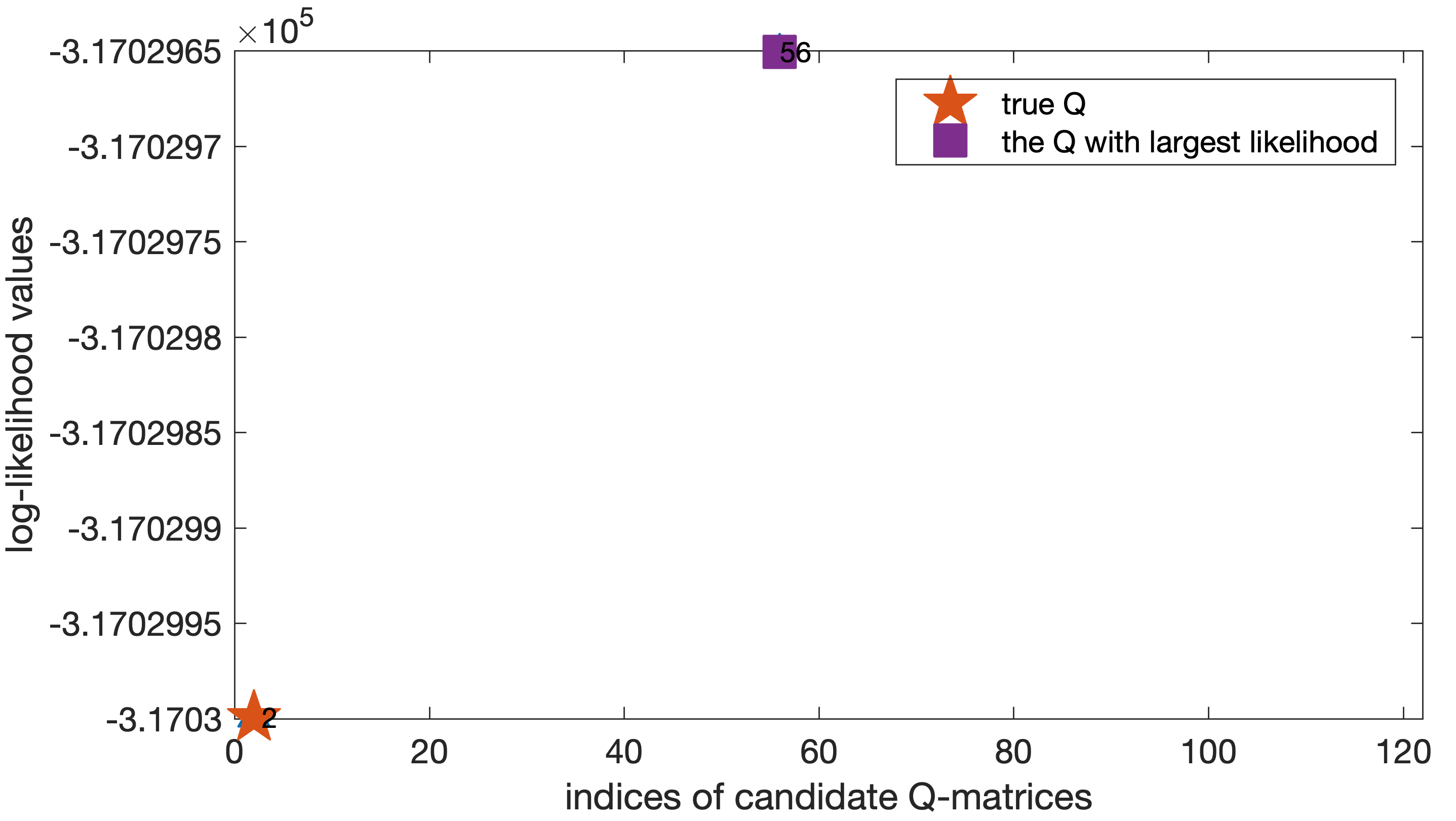}
\caption{GDINA: true $Q$ not even locally identifiable: $\small
Q^{2} = \begin{pmatrix}
  0   &  1  &   0   &  0   &  0\\
  1 &    0   &  1  &   1  &   1
\end{pmatrix}^\top
$}
\end{subfigure}
\end{figure}

\begin{figure}[H]\ContinuedFloat
\centering

\begin{subfigure}{0.75\textwidth}
\includegraphics[width=\linewidth]
{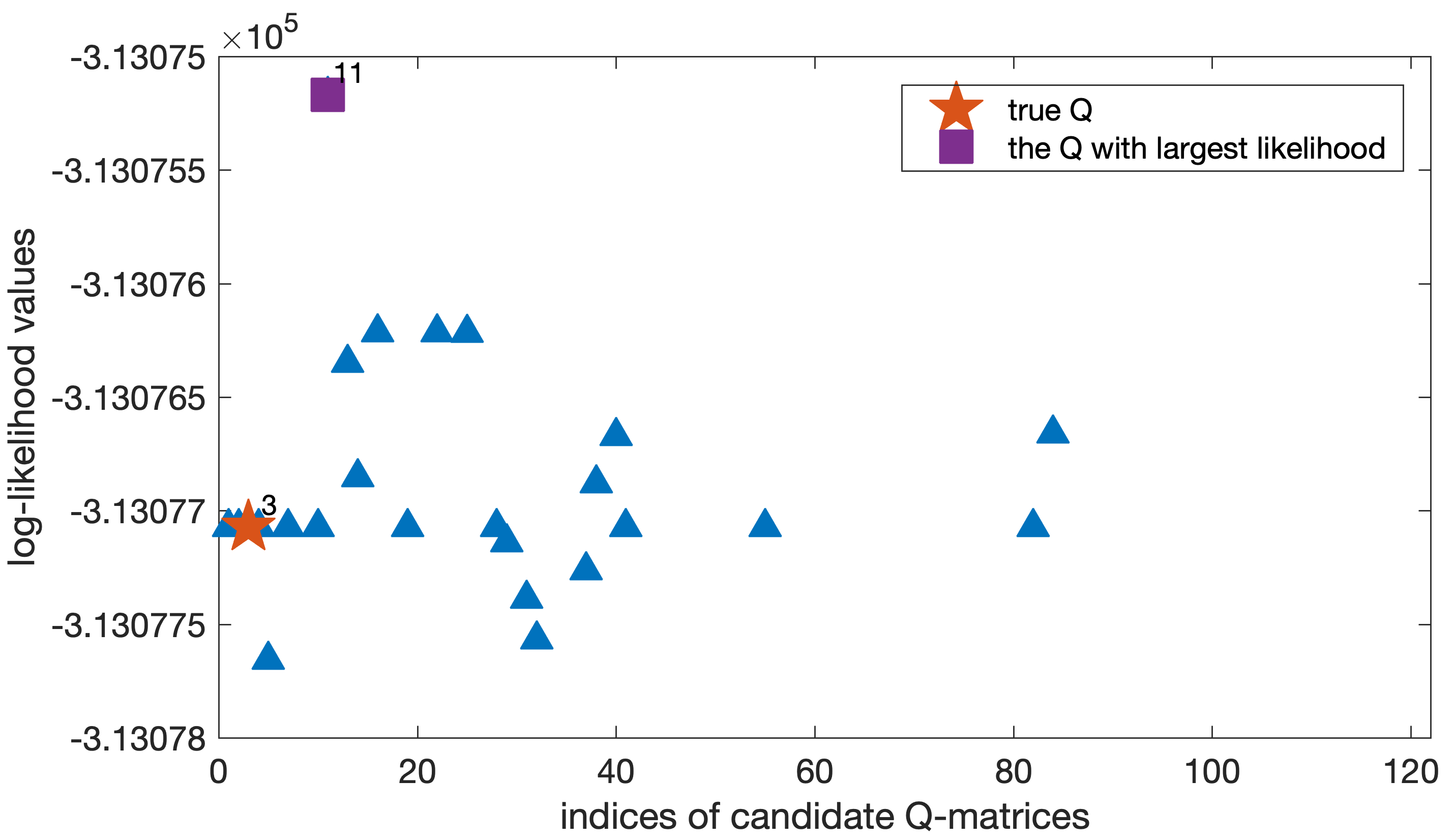}
\caption{GDINA: true $Q$ not even locally identifiable: $\small
Q^{3} = \begin{pmatrix}
  0   &  1  &   0   &  0   &  0\\
  1 &    1   &  1  &   1  &   1
\end{pmatrix}^\top
$}
\end{subfigure}

\end{figure}

\newpage
\noindent
\textbf{Study VII: Construction of many alternative sets of parameters when true $Q$-matrix violates the necessary condition for generic identifiability.}

In this study, we verify Theorem \ref{prop-gdina-nece}, i.e., verify the necessity of Condition $C$ that \textit{each attribute is required by at least two items in the $Q$-matrix} for joint generic identifiability.
We consider two cases with $(K, J)=(3,20)$ and $(K,J)=(5,20)$.

First, for $(K,J)=(3,20)$, consider the following $Q$-matrix $Q_3$ and an alternative $\bar Q_3$. 
\begin{equation}\label{eq-q3-attr2}
\small
Q_3 = \begin{pmatrix}
	 \mathbf 1  &   \mathbf 1  &   \mathbf 0\\
     \mathbf 1  &   \mathbf 0  &   \mathbf 1\\
     0  &   1  &   0\\
     0  &   0  &   1\\
     0  &   1  &   1\\
     0  &   1  &   0\\
     0  &   0  &   1\\
     0  &   1  &   1\\
     0  &   1  &   0\\
     0  &   0  &   1\\
     0  &   1  &   1\\
     0  &   1  &   0\\
     0  &   0  &   1\\
     0  &   1  &   1\\
     0  &   1  &   0\\
     0  &   0  &   1\\
     0  &   1  &   1\\
     0  &   1  &   0\\
     0  &   0  &   1\\
     0  &   1  &   1
\end{pmatrix}_{20\times 3}\quad
\bar Q_3 = \begin{pmatrix}
	 \mathbf 1  &   \mathbf 1  &   \mathbf 1\\
     \mathbf 1  &   \mathbf 1  &   \mathbf 1\\
     0  &   1  &   0\\
     0  &   0  &   1\\
     0  &   1  &   1\\
     0  &   1  &   0\\
     0  &   0  &   1\\
     0  &   1  &   1\\
     0  &   1  &   0\\
     0  &   0  &   1\\
     0  &   1  &   1\\
     0  &   1  &   0\\
     0  &   0  &   1\\
     0  &   1  &   1\\
     0  &   1  &   0\\
     0  &   0  &   1\\
     0  &   1  &   1\\
     0  &   1  &   0\\
     0  &   0  &   1\\
     0  &   1  &   1
\end{pmatrix}_{20\times 3}
\end{equation}
We first construct true parameters $(\TT,\pp)$ under $Q_3$. For each attribute pattern $\aaa$, we set its population proportion $p_{\aaa}$ to be $1/2^K$. For each item, set the baseline probability, the positive response probability of the all-zero attribute profile $\aaa=\mathbf 0^\top$, to be $0.2$ and the positive response probability of $\aaa=\mathbf 1^\top$ to be $0.8$. And we take all the main effects and interaction effects parameters to be equal.

%In summary, the item parameter matrix $\TT$ of size $J\times 2^K$ and the proportion vector $\pp$ of length $2^K$ can be written as
%$$
%\TT=
%\begin{pmatrix}
%0.2  &  0.2  &  0.4 &  0.4  &  0.4 &   0.4 &   0.8  &  0.8 \\
%0.2  &  0.4  &  0.2 &  0.4  &  0.4 &  0.8  &  0.4  &  0.8\\
%0.2  &  0.2  &  0.8 &  0.8 &  0.2  &  0.2   & 0.8   & 0.8\\
%0.2  &  0.8  &  0.2 &  0.8 &   0.2 &   0.8  & 0.2   & 0.8\\
%0.2  &  0.2  &  0.8 &  0.8 &   0.2 &   0.2   & 0.8  &  0.8\\
%0.2  &  0.8  &  0.2 &  0.8 &   0.2  &  0.8   & 0.2   & 0.8\\
%0.2  &  0.4  &  0.4 &  0.8 &  0.2  & 0.4   & 0.4   & 0.8\\
%\end{pmatrix},~
%%\end{block}
%%\end{blockarray}~,
%\quad
%\pp = \begin{pmatrix}
%0.125\\
%0.125\\
%0.125\\
%0.125\\
%0.125\\
%0.125\\
%0.125\\
%0.125\\
%\end{pmatrix}
%$$
%where the columns of $\TT$ are indexed by the eight attribute profiles. %$000,~001,~010,~011,~100,~101,\allowbreak~110,~111$.
%Under $\TT$ and $\pp$, following the G-DINA model assumption we can get the  Probability Mass Function (PMF) for the response pattern vector $\RR$, which is a $2^J = 128$ dimensional vector $\PP(Q,\TT,\pp)$ storing  $\{P(\RR=\rr\mid Q,\TT,\pp)$ for all $\rr\in\{0,1\}^J$.  Here 
%\begin{align*}
%\PP(Q,\TT,\pp)&=(\PP_{\rr_1}(Q,\TT,\pp),\,\ldots,\,\PP_{\rr_{2^J}}(Q,\TT,\pp))\\
%&=( 0.0433, ~   0.0183  ,~  0.0183  , ~ 0.0152,\allowbreak\ldots,~0.0338).
%\end{align*}
 
For the defined true parameters $(\TT,\pp)$ under $Q_3$, we next construct 70 alternative sets of parameters $(\bar\TT^\ell,\bar\pp^\ell)$ for $\ell=1,2,\ldots,70$, all under the alternative $Q$-matrix $\bar Q_3$, that are non-distinguishable from the true parameters.
Following the proof of Theorem \ref{prop-gdina-nece}, 
we first set $\bar\theta_{j,\aaa}=\theta_{j,\aaa}$ for any $j>2$ and any $\aaa$. Then we \textit{randomly generate} the values of the $\bar\TT_{1:2,\,1:4}$ (the first four elements of the first two rows of $\bar\TT$) from the neighborhood of their true values, and enforce the monotonicity constraint \eqref{eq-mono}. Specifically, for each alternative set (the $\ell$-th set) of parameters, there is
$$
\bar\TT^{\ell}_{i, j} = \TT_{i, j} + \mathcal U(-0.1, 0.1), \quad i=1,2;~ j=1,2,3,4;~\ell=1,2,\ldots,70.
$$
where $\mathcal U(-0.1, 0.1)$ denotes a uniformly distributed random variable between $-0.1$ and $0.1$. Next we just use Equation \eqref{eq-sol} to get the remaining item parameters $\bar\TT^{\ell}_{1:2,\,5:8}$ and $\bar\pp^{\ell}$. 

\begin{figure}[H]
\caption{GDINA: true $Q$ is $Q_3$ with $(K,J)=(3,20)$; each of the 70 colored line corresponds to one set of alternative parameters under $\bar Q_3$; all sets  non-distinguishable.}
\label{fig-gdina_attr2_K3}

%\vspace*{-3mm}
\centering
\begin{subfigure}{0.72\textwidth}
\includegraphics[width=\linewidth]
{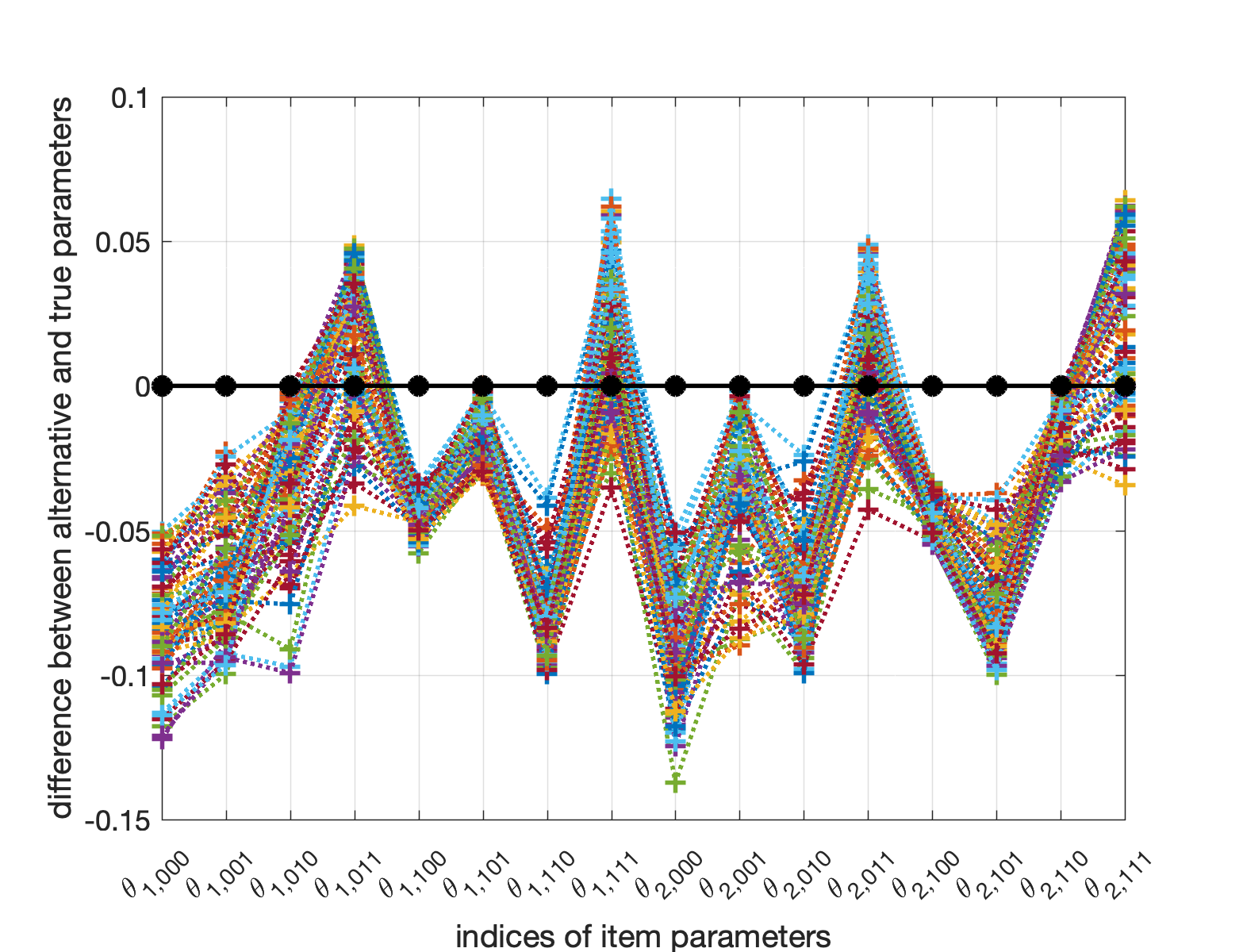}
\caption{$K=3$ and  $J=20$, $70$ alternative sets of parameters} 
\end{subfigure}

%\vspace{-3mm}
\begin{subfigure}{0.72\textwidth}
\includegraphics[width=\linewidth]
{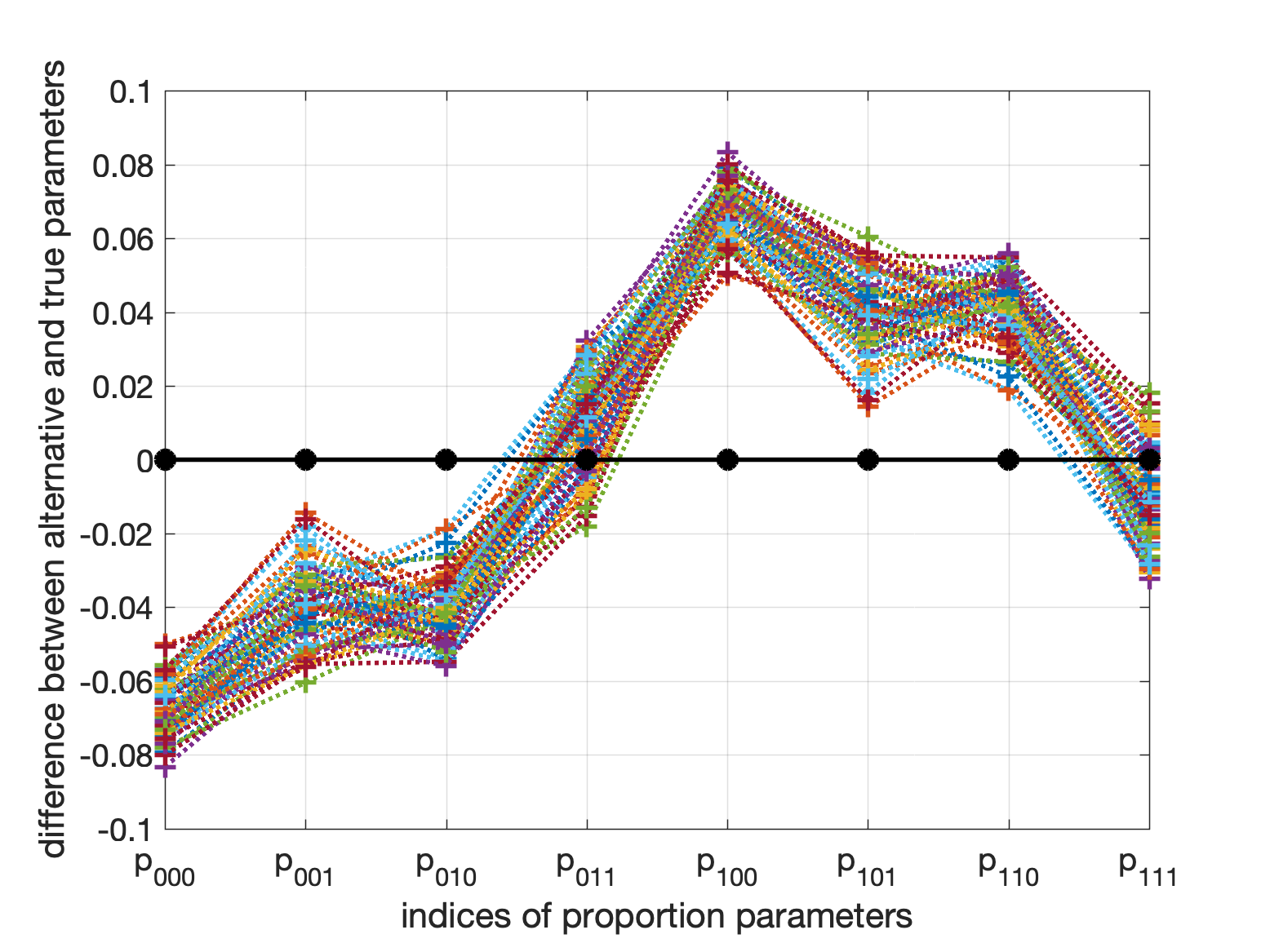}
\caption{$K=3$ and  $J=20$, $70$ sets of parameters} 
\end{subfigure}
\end{figure}

%Along with the true parameters $(\TT,\pp)$ under $Q_3$, 
Figure \ref{fig-gdina_attr2_K3} presents the constructed 70 other parameters sets $(\bar\TT^\ell,\bar\pp^\ell)$ under the alternative $\bar Q_3$, by plotting the values of difference between the alternative parameters and the true parameters.
In particular, In Figure \ref{fig-gdina_attr2_K3}(a), the black solid line with dots is the reference line at zero, and each of the 70 colored dotted line with ``$+$'''s represents one particular set of alternative parameters. For each colored line corresponding to the $\ell$th set of parameters, the following $16$-dimensional vector of parameter difference is plotted,
\begin{align*}
(
&\bar\theta^\ell_{1,\,000}-\theta_{1,\,000},
~\bar\theta^\ell_{1,\,001}-\theta_{1,\,010},
~\bar\theta^\ell_{1,\,010}-\theta_{1,\,010},
~\bar\theta^\ell_{1,\,011}-\theta_{1,\,011},\\
~&\bar\theta^\ell_{1,\,100}-\theta_{1,\,100},
~\bar\theta^\ell_{1,\,101}-\theta_{1,\,110},
~\bar\theta^\ell_{1,\,110}-\theta_{1,\,110},
~\bar\theta^\ell_{1,\,111}-\theta_{1,\,111},\\
~&\bar\theta^\ell_{2,\,000}-\theta_{2,\,000},
~\bar\theta^\ell_{2,\,001}-\theta_{2,\,010},
~\bar\theta^\ell_{2,\,010}-\theta_{2,\,010},
~\bar\theta^\ell_{2,\,011}-\theta_{2,\,011},\\
~&\bar\theta^\ell_{2,\,100}-\theta_{2,\,100},
~\bar\theta^\ell_{2,\,101}-\theta_{2,\,110},
~\bar\theta^\ell_{2,\,110}-\theta_{2,\,110},
~\bar\theta^\ell_{2,\,111}-\theta_{2,\,111}
).
\end{align*}
Similarly, in Figure \ref{fig-gdina_attr2_K3}(b), for each colored line corresponding to the $\ell$th set of parameters, the following $8$-dimensional vector of parameter difference is plotted,
%\begin{align*}
$(
\bar p^\ell_{000}- p_{000},
~\bar p^\ell_{001}- p_{010},
~\bar p^\ell_{010}- p_{010},
~\bar p^\ell_{011}- p_{011},
~\bar p^\ell_{100}- p_{100},
~\bar p^\ell_{101}- p_{110},
~\bar p^\ell_{110}- p_{110},
~\bar p^\ell_{111}- p_{111}
).$
%\end{align*}
In summary, a total number of 70 colored lines corresponding to 70 alternative sets of parameters are plotted in Figure \ref{fig-gdina_attr2_K3}.

%The $(\bar\TT^\ell,\bar\pp^\ell)$ are constructed based on randomly setting several parameters and then determining other parameters following our theoretical derivation non-identifiability.
The $(\TT,\pp)$ and all the $(\bar\TT^\ell,\bar\pp^\ell)$, $\ell=1,\ldots,70$ give the identical distribution of $\RR$. Specifically, from the computation in Matlab, we have
\begin{align*}
\max_{1\leq\ell\leq 70}\max_{\rr\in\{0,1\}^{20}}\left|\mathbb P(\RR=\rr\mid Q_3,\TT,\pp)-\mathbb P(\RR=\rr\mid \bar Q_3,\TT^\ell,\pp^\ell)\right|
&= 1.30\times 10^{-18},
\end{align*}
which is smaller than the Matlab machine error $2.22\times 10^{-16}$. This verifies that despite the underlying parameters are different from the truth, they all lead to the identical distribution of responses. So $(Q_3,\TT,\pp)$ are not identifiable.
We emphasize that under $Q_3$, for any true parameters, one can  construct arbitrarily many such alternative parameter sets under $\bar Q_3$.

\newpage

\begin{figure}[H]
\caption{GDINA: true $Q$ is $Q_4$ with $(K,J)=(5,20)$; each of the 70 colored line corresponds to one set of alternative parameters under $\bar Q_4$; all sets  non-distinguishable.}
\label{fig-gdina_attr2_K5}

%\vspace*{-3mm}
\centering
\begin{subfigure}{0.72\textwidth}
\includegraphics[width=\linewidth]
{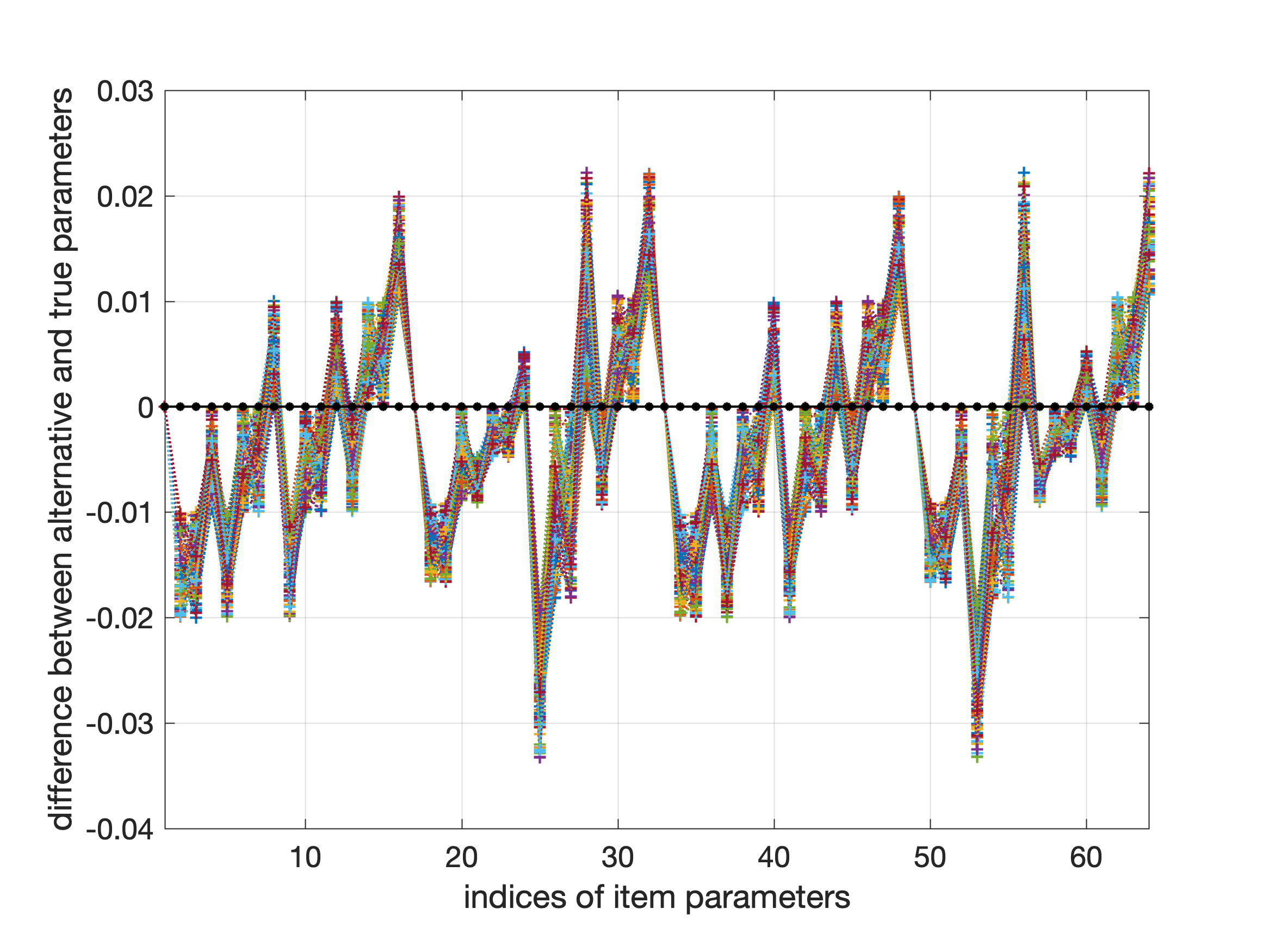}
%\vspace*{-12mm}
\caption{$K=5$ and  $J=20$, $70$ alternative sets of parameters}  %marginal prob. for response patterns}
\end{subfigure}

%\hspace*{\fill}
%\vspace{-3mm}
\begin{subfigure}{0.72\textwidth}
\includegraphics[width=\linewidth]
{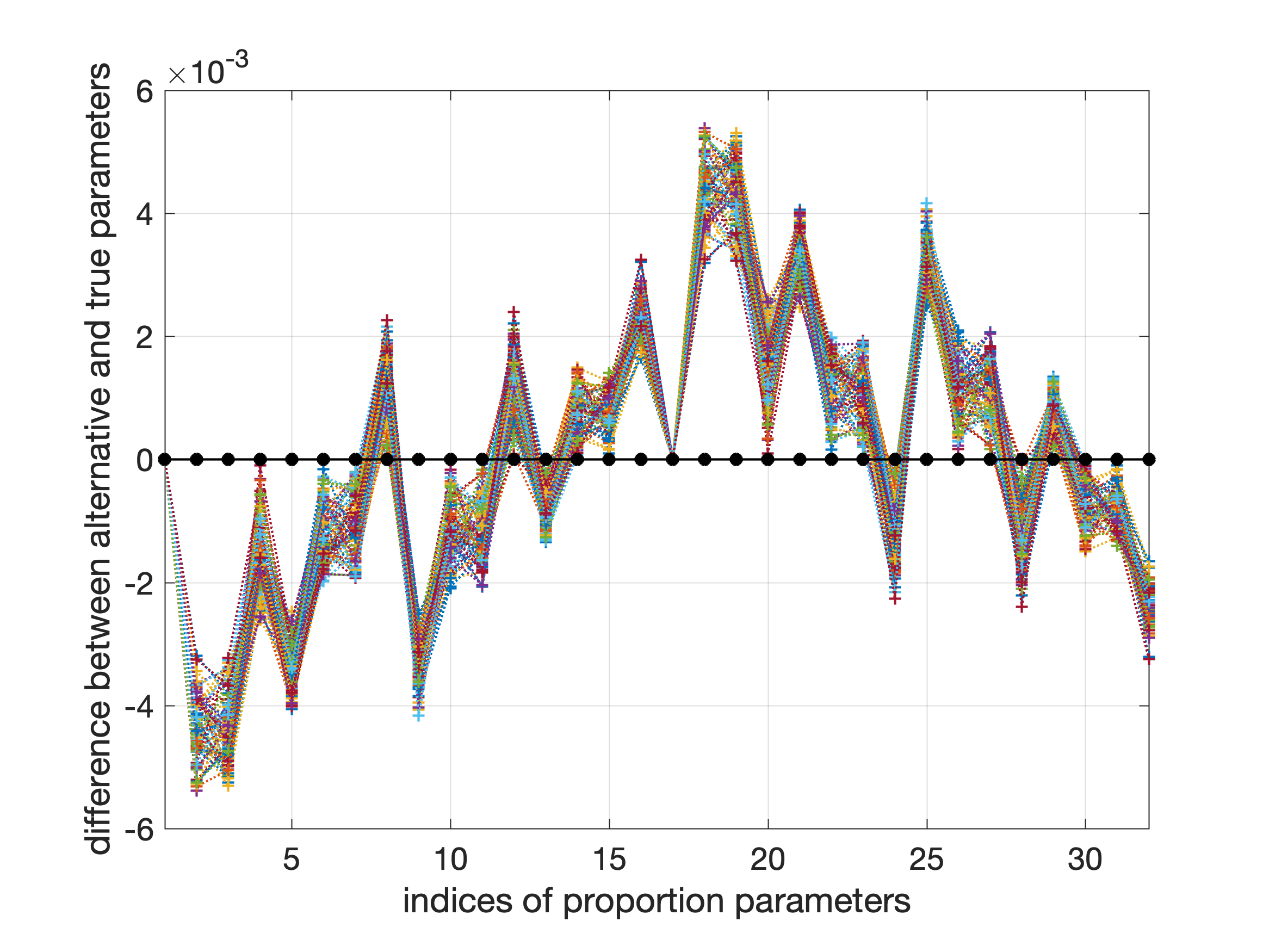}
%\vspace*{-12mm}
\caption{$K=5$ and  $J=20$, $70$ alternative sets of parameters} 
\end{subfigure}

\end{figure}

For $(K,J)=(5,20)$, consider the following $Q_4$ and an alternative $\bar Q_4$,
\begin{equation}\small
Q_4 = \begin{pmatrix}
     \mathbf 1  &  \mathbf  1  &   \mathbf 0  &  \mathbf  0  &   \mathbf 0\\
     \mathbf 1  &  \mathbf  0  &  \mathbf  1  & \mathbf   0  &   \mathbf 0\\
     0  &   1  &   0  &   0  &   0\\
     0  &   0  &   1  &   0  &   0\\
     0  &   0  &   0  &   1  &   0\\
     0  &   0  &   0  &   0  &   1\\
     0  &   1  &   0  &   0  &   0\\
     0  &   0  &   1  &   0  &   0\\
     0  &   0  &   0  &   1  &   0\\
     0  &   0  &   0  &   0  &   1\\
     0  &   1  &   0  &   0  &   0\\
     0  &   0  &   1  &   0  &   0\\
     0  &   0  &   0  &   1  &   0\\
     0  &   0  &   0  &   0  &   1\\
     0  &   1  &   1  &   0  &   0\\
     0  &   1  &   0  &   1  &   0\\
     0  &   1  &   0  &   0  &   1\\
     0  &   0  &   1  &   1  &   0\\
     0  &   0  &   1  &   0  &   1\\
     0  &   0  &   0  &   1  &   1
\end{pmatrix}_{20\times 5}\quad
\bar Q_4 = \begin{pmatrix}
     \mathbf 1  &  \mathbf  1  &   \mathbf 1  &  \mathbf  1  &   \mathbf 1\\
     \mathbf 1  &  \mathbf  1  &  \mathbf  1  & \mathbf   1  &   \mathbf 1\\
     0  &   1  &   0  &   0  &   0\\
     0  &   0  &   1  &   0  &   0\\
     0  &   0  &   0  &   1  &   0\\
     0  &   0  &   0  &   0  &   1\\
     0  &   1  &   0  &   0  &   0\\
     0  &   0  &   1  &   0  &   0\\
     0  &   0  &   0  &   1  &   0\\
     0  &   0  &   0  &   0  &   1\\
     0  &   1  &   0  &   0  &   0\\
     0  &   0  &   1  &   0  &   0\\
     0  &   0  &   0  &   1  &   0\\
     0  &   0  &   0  &   0  &   1\\
     0  &   1  &   1  &   0  &   0\\
     0  &   1  &   0  &   1  &   0\\
     0  &   1  &   0  &   0  &   1\\
     0  &   0  &   1  &   1  &   0\\
     0  &   0  &   1  &   0  &   1\\
     0  &   0  &   0  &   1  &   1
\end{pmatrix}_{20\times 5}.
\end{equation}
We set the true parameters under $Q_4$ similarly as those under $Q_3$, and also use \eqref{eq-sol} in the proof of Theorem \ref{prop-gdina-nece} to randomly construct 70  sets of parameters under the $\bar Q_4$. Figure \ref{fig-gdina_attr2_K5} (a) and (b) plot the values of difference between alternative and true item parameters (of the first two items), and that between alternative and true proportion parameters, respectively. Despite the differences in parameter values, our computation in Matlab shows the maximum difference between marginal response probabilities is
\begin{align*}
\max_{1\leq\ell\leq 70}\max_{\rr\in\{0,1\}^{20}}\left|\mathbb P(\RR=\rr\mid Q_4,\TT,\pp)-\mathbb P(\RR=\rr\mid \bar Q_4,\TT^\ell,\pp^\ell)\right|
&= 5.42\times 10^{-19},
\end{align*}
also smaller than the Matlab machine error $2.22\times 10^{-16}$. This illustrates the non-identifiability of $Q_4$.

\newpage

\bibliographystyle{apa}
\bibliography{ref}

\vskip .65cm
\noindent
Yuqi Gu and Gongjun Xu 
\vskip 2pt
\noindent
Department of Statistics, University of Michigan 
\vskip 2pt
\noindent
E-mail: yuqigu, gongjun@umich.edu
\vskip 2pt

%\noindent
%second author affiliation
%\vskip 2pt
%\noindent
%E-mail: (second author email)

\end{document}